\documentclass[leqno,a4,11pt]{article}
\usepackage{amsmath,amssymb,amsthm,longtable}
\numberwithin{equation}{section}
\newtheorem{thm}{Theorem}[section]
\newtheorem{lem}{Lemma}[section]
\newtheorem{prop}{Proposition}[section]
\newcommand{\lw}[1]{\smash{\lower2.0ex\hbox{#1}}}
\makeatletter
\def\@seccntformat#1{\csname the#1\endcsname.\quad}
\def\section{\@startsection {section}{1}{\z@}{ 2.3ex plus 2ex minus
 -.2ex}{2.3 ex plus .2ex}{\Large\bf}}
\def\subsection{\@startsection {subsection}{1}{\z@}{ 2.3ex plus 2ex minus
 -.2ex}{2.3 ex plus .2ex}{\large\bf}}
\makeatother

\begin{document}
\title{\bf Involutions of compact Riemannian \\ $4$-symmetric spaces}
\author{Hiroyuki Kurihara and Koji Tojo}
\date{}
\maketitle
\begin{abstract}
Let $G/H$ be a compact $4$-symmetric space of inner type such that 
the dimension of the center $Z(H)$ of $H$ is at most one. 
In this paper we shall classify involutions of $G$ preserving $H$ 
for the case where $\dim Z(H)=0$, or $H$ is a centralizer of a toral subgroup 
of $G$. 
\end{abstract}


\section{Introduction}

 It is known that Riemannian $k$-symmetric spaces is a generalizations 
of Riemannian symmetric spaces. The definition is as follows:\par
  Let $G$ be a Lie group and $H$ a compact subgroup of $G$. A homogeneous
 space $(G/H, \langle,\rangle)$ with $G$-invariant Riemannian metric
 $\langle,\rangle$ is called 
 a {\it Riemannian $k$-symmetric space} if there
 exists an automorphism $\sigma$ of on $G$ such that 
\begin{enumerate}
\item $G^{\sigma}_o \subset H \subset G^{\sigma}$, where $G^{\sigma}$
      and $G^{\sigma}_o$ is the set of fixed points of $\sigma$ and its
      identity component, respectively,
\item $\sigma^k={\rm Id}$ and $\sigma^l\ne{\rm Id}$ for any $l<k$,
\item The transformation of $G/H$ induced by $\sigma$ is an isometry.
\end{enumerate}
 We denote by $(G/H, \langle,\rangle, \sigma)$ a Riemannian $k$-symmetric
 space with an automorphism $\sigma$. Gray \cite{G} classified
 Riemannian $3$-symmetric spaces (see also Wolf and Gray~\cite{WG}). 
Moreover compact Riemannian $4$-symmetric spaces is classified by 
Jem\'{e}nez~\cite{J}. The structure of Riemannian $k$-symmetric
 spaces is closely rerated to the study of finite order automorphisms of
 Lie groups. Such automorphisms of compact simple Lie groups were classified
 (cf. Kac~\cite{K} and Helgason~\cite{Helg}).  
\par
It is known that involutions on $k$-symmetric spaces are important. 
  For example, the classifications of affine symmetric spaces by 
Berger~\cite{Be} are, in essence, the classification of involutions on 
compact symmetric spaces $G/H$ preserving $H$. Similarly, 
such involutions play an important role in the classification of symmetric 
submanifolds on compact symmetric spaces (cf. Naitoh \cite{N1} and 
\cite{N2}). \par
On a compact $3$-symmetric space $(G/H,\langle,\rangle,\sigma)$, an 
involution $\tau$ preserving $H$ satisfies 
$\tau\circ\sigma=\sigma\circ\tau$ or $\tau\circ\sigma=\sigma^{-1}\circ\tau$. 
The classification of affine $3$-symmetric spaces (\cite{WG})
was made by classifying involutions $\tau$ satisfying $\tau\circ\sigma=
\sigma\circ\tau$. Moreover, \cite{T1} and \cite{T2} classify half-dimensional,
 totally real and totally geodesic submanifold (with respect to the canonical 
almost complex structures) of compact Riemannian $3$-symmetric spaces 
$(G/H,\langle,\rangle,\sigma)$ by classifying involutions $\tau$ on $G$ 
satisfying $\tau\circ\sigma=\sigma^{-1}\circ\tau$.
 \par
In general, there exists an involution $\tau$ such that
 $\tau\circ\sigma\circ\tau^{-1}\ne\sigma$ or $\sigma^{-1}$ for Reimannian 
$4$-symmetric spaces. These automorphisms do not appear in Riemannian 
symmetric spaces and $3$-symmetric spaces. However, if the dimension of the 
center of $H$ is at most one, each involution $\tau$ preserving $H$ satisfies 
$\tau\circ\sigma\circ\tau^{-1}=\sigma$ or $\sigma^{-1}$. 
\par
According to \cite{J}, a compact simply connected Riemannian
$4$-symmetric space decomposes as a product $M_1\times\cdots\times
M_r$, where $M_i(1\leq i\leq r)$ is compact, irreducible Riemannian 
$4$-symmetric space. In this paper we treat a compact, irreducible Riemannian 
$4$-symmetric space $(G/H,\langle,\rangle,\sigma)$ such that the dimension of 
the center of $H$ is at most one. In purticular we classify involutions of $G$
 preserving $H$ for the case where $\dim Z(H)=0$, or $\dim Z(H)=1$ and $H$ is 
a centralizer of a toral subgroup of $G$. 
More precisely, let $\mathfrak{g}$ and $\mathfrak{h}$ be the Lie algebras of 
$G$ and $H$, respectively. Then we first prove that there exists a maximal 
abelian subalgebra $\mathfrak{t}$ of $\mathfrak{g}$ contained in 
$\mathfrak{h}$ such that $\tau(\mathfrak{t})=\mathfrak{t}$ for any involution 
$\tau$ preserving $\mathfrak{h}$. Except for the case where $\dim Z(H)=1$ and 
$\tau\circ\sigma\circ\tau^{-1}=\sigma^{-1}$, we classify involutions 
$\bar{\tau}$ of the root system of $\mathfrak{h}$ with respect to 
$\mathfrak{t}$. Moreover, for each involution $\bar{\tau}$ 
($\bar{\tau}\not={\rm Id}$) of the root system of $\mathfrak{h}$, we prove 
that there exists an involution $\tau_0$ preserving $\mathfrak{h}$ such that 
$\tau_0|_{\mathfrak{h}}=\bar{\tau}$. Then each involution $\tau$ can be 
written as $\tau=\tau_0\circ{\rm Ad}(\exp\sqrt{-1}h)$ or 
$\tau={\rm Ad}(\exp\sqrt{-1}h)$ for some $\sqrt{-1}h\in\mathfrak{t}$ since 
$\tau|_{\mathfrak{t}}$ is an involution of the root system of $\mathfrak{h}$, 
and we obtain all $\tau$ by considering conjugations within automorphisms 
preserving $\mathfrak{h}$. For the case where $\dim Z(H)=1$ and 
$\tau\circ\sigma\circ\tau^{-1}=\sigma^{-1}$, using graded Lie algebras, 
we classify all $\tau$ by an argument similar to that in \cite{T1}. 
\par
According to \cite{T2}, for $3$-symmetric spaces 
$(G/H,\langle,\rangle,\sigma)$ with 
$\dim Z(H)=0$, each involution $\tau$ with $\tau\circ\sigma\circ\tau^{-1}=
\sigma^{-1}$ preserving $H$ is obtained from a grade-reversing Cartan 
involution of some graded Lie algebra of the third kind. 
In the case where $(G/H,\langle,\rangle,\sigma)$ is $4$-symmetric with $\dim Z(H)=0$ 
and $\tau\circ\sigma\circ\tau^{-1}=\sigma^{-1}$, we can see that there exists 
$\tau$ which is not obtained from a grade-reversing Cartan involution of any 
graded Lie algebra of the fourth kind.

\par 
 The organization  of this paper is as follows:\par 
  In Section 2, we recall the notions of root systems and graded Lie algebras 
needed for the remaining part of this paper. Moreover we recall some results 
on automorphisms of order $k$ ($k\leq 4$). \par
  In Section 3, we remark on some relation between involutions of 
$4$-symmetric space $(G/H,\langle,\rangle,\sigma)$ reserving $H$ 
and root systems of the Lie algebra of $G$. \par
In Section 4, by using the results in Section 3, we describe the restrictions 
of involutions to the root systems for the case where the dimension of the 
center is zero.  
\par  
In Section 5-8, we enumerate all involutions $\tau$ of compact $4$-symmetric 
spaces such that $\tau(H)=H$ and the dimension of the center of $H$ is
 zero, or $H$ is a centralizer of a toral subgroup of $G$. \par
In Section 9, we describe some conjugations between involutions. \par
In Section 10, by making use of the results in Section 5-8 together with 
conjugations in Section 9, we give the classification theorem of the
 equivalence classes of involutions.

\section{Preliminaries}
\subsection{Root systems}
  Let $\mathfrak{g}$ and $\mathfrak{t}$ be a compact semisimple Lie algebra
 and a maximal abelian subalgebra of $\mathfrak{g}$, respectively. We
 denote by $\mathfrak{g}_{\mathbb{C}}$ and $\mathfrak{t}_{\mathbb{C}}$ the
 complexifications of $\mathfrak{g}$ and $\mathfrak{t}$, respectively. Let
 $\varDelta (\mathfrak{g}_{\mathbb{C}},\mathfrak{t}_{\mathbb{C}})$ be the
 root system of $\mathfrak{g}_{\mathbb{C}}$ with respect to $\mathfrak
 t_{\mathbb{C}}$ and $\varPi (\mathfrak{g}_{\mathbb{C}},\mathfrak
  t_{\mathbb{C}}) = \{\alpha_1, \ldots ,\alpha_n\}$ the set of 
 fundamental roots of $\varDelta (\mathfrak{g}_{\mathbb{C}},\mathfrak
  t_{\mathbb{C}})$ with respect to a lexicographic order. 
For $\alpha\in\varDelta (\mathfrak{g}_{\mathbb{C}},\mathfrak{t}_{\mathbb{C}})$,
put 
\begin{equation}
  \mathfrak{g}_{\alpha} := \{X\in\mathfrak{g}_{\mathbb{C}} \ ; \ [H,X] =
  \alpha(H)X \ {\rm for \ any} \ H\in \mathfrak{t}_{\mathbb{C}}\}.
   \label{eqn:21}
\end{equation}
Since the Killing form $B$ of $\mathfrak{g}_{\mathbb{C}}$ is nondegenerate, 
we can define $H_{\alpha}\in\mathfrak{t}_{\mathbb{C}}$ 
$(\alpha\in\varDelta(\mathfrak{g}_{\mathbb{C}}, \mathfrak{t}_{\mathbb{C}}))$ 
 by $\alpha (H) = B(H_{\alpha}, H)$ for any $H\in \mathfrak{t}
_{\mathbb{C}}$. As in \cite{Helg}, we take the Weyl basis
 $\{E_{\alpha} \in \mathfrak{g}_{\alpha} \ ; \ \alpha \in \varDelta
 (\mathfrak{g}_{\mathbb{C}},\mathfrak{t}_{\mathbb{C}})\}$ of
 $\mathfrak{g}_{\mathbb{C}}$ so that 
  
$$\begin{array}{ll}
  &[E_{\alpha},E_{-\alpha}]  = H_{\alpha}, \\
  &[E_{\alpha}, E_{\beta}] = N_{\alpha,\beta} E_{\alpha + \beta}, \ \
   N_{\alpha,\beta}\in\mathbb{R}, \\ 
  &N_{\alpha,\beta} = -N_{-\alpha,-\beta},\\
  &A_{\alpha}:= E_{\alpha} - E_{-\alpha}, \ \ B_{\alpha}:=
   \sqrt{-1}(E_{\alpha} + E_{-\alpha})\in \mathfrak{g}.
\end{array}$$
 We denote by $\varDelta^+(\mathfrak{g}_{\mathbb{C}}, \mathfrak{t}_{\mathbb{C}})$
 the set of positive roots of $\varDelta
(\mathfrak {g_{\mathbb{C}}, t_{\mathbb{C}}})$ with respect to the order. 
Then it follows that 

\begin{equation}
  \mathfrak{g} = \mathfrak{t}
+\sum_{\alpha\in\varDelta^+(\mathfrak{g}_{\mathbb{C}},\mathfrak{t}_{\mathbb{C}})}
({\mathbb{R}}  A_{\alpha}+{\mathbb{R}} B_{\alpha}),\ \ 
\mathfrak{t}=\sum_{i=1}^{n}{\mathbb{R}}\sqrt{-1}H_{\alpha_i}.
\end{equation}
For $\alpha\in\varDelta(\mathfrak{g}_{\mathbb{C}}, \mathfrak{t}_{\mathbb{C}})$, 
define a Lie subalgebra $\mathfrak{su}_{\alpha}(2)$ of $\mathfrak{g}$ by 

\begin{equation}
   \mathfrak {su} _{\alpha}(2):= \mathbb{R} \sqrt{-1}H_{\alpha}+\mathbb{R}
   A_{\alpha}+\mathbb{R} B_{\alpha}\cong \mathfrak {su} (2).
\label{eqn:2.3}
\end{equation}
We denote by $t_{\alpha}$ the root reflection for 
$\alpha\in \varDelta (\mathfrak{g}_{\mathbb{C}},\mathfrak{t}_{\mathbb{C}})$. 
Then there exists an extension of $t_{\alpha}$ to an element of 
the group ${\rm Int}(\mathfrak{g})$ of inner automorphisms of $\mathfrak{g}$, 
which is denoted by the same symbol as $t_{\alpha}$. 
Since the root reflection of $\mathfrak{su}_{\alpha}(2)$ for $\alpha$ 
coincides with the restriction of $t_{\alpha}$ to 
${\mathbb{R}}\sqrt{-1}H_{\alpha}$ and $t_{\alpha}$ is the identical 
transformation on the orthogonal complement of 
${\mathbb{R}}\sqrt{-1}H_{\alpha}$ in $\mathfrak{t}$, the 
following lemma holds. 

\begin{lem}
There exists an element $\phi\in{\rm Int}(\mathfrak{su}_{\alpha}(2))
(\subset{\rm Int}(\mathfrak{g}))$
 such that $\phi|_{\mathfrak{t}}=t_{\alpha}|_{\mathfrak{t}}$. 
 \label{lem:21}
\end{lem}

Define $K_j\in\mathfrak{t}_{\mathbb{C}} \ \ (j = 1,\ldots, n) $ by
$$\alpha_i (K_j) = \delta_{ij}, \ \ \ i,j = 1, \ldots ,n,$$
and denote the highest root $\delta$ by
$$\delta := \sum_{j=1}^n m_j \alpha_j, \ \ m_j\in\mathbb Z.$$
We set
$$\tau_H := {\rm Ad}(\exp\pi\sqrt{-1}H),\ \ H\in\mathfrak
t_{\mathbb{C}}.$$
Then from (\ref{eqn:21}) we have
\begin{equation}
  \tau_H(E_{\alpha})=e^{\pi\sqrt{-1}\alpha(H)}E_{\alpha}, \ 
\alpha\in\varDelta(\mathfrak{g}_{\mathbb{C}},\mathfrak{t}_{\mathbb{C}}).
    \label{eqn:20}
\end{equation}

Assume that $\mathfrak{g}$ is simple. Then the following is known. 

\begin{lem}[\cite{M}]
  Any inner automorphism of order 2 on $\mathfrak{g}$ is conjugate within
  ${\rm Int} (\mathfrak{g})$ to some $\tau_{K_i}$ with $m_i=1$ or $2$.
      \label{lem:22}
\end{lem}

If $h-h^{\prime} = \sum_{i=1}^n a_i K_i, \ a_i\in 2\mathbb Z$ for
$h,h^{\prime}\in \mathfrak{t}_{\mathbb{C}}$, we say that $h$ is congruent
to $h^{\prime}$ modulo 
$2\varPi(\mathfrak{g}_{\mathbb{C}},\mathfrak{t}_{\mathbb{C}})$ and it is denoted by 
$h\equiv h^{\prime}({\rm mod} \ 
2\varPi(\mathfrak{g}_{\mathbb{C}},\mathfrak{t}_{\mathbb{C}}))$. It follows from 
(\ref{eqn:20}) that $\tau_h=\tau_{h'}$ if $h\equiv h'({\rm mod}\ 
2\varPi(\mathfrak{g}_{\mathbb{C}},\mathfrak{t}_{\mathbb{C}}))$. 
\vspace{0.5cm}\par

{\sc Remark 2.1.}
   According to Lemma \ref{lem:22}, for any inner automorphism $\tau_H$
   of order $2$ on $\mathfrak{g}$, there exists an inner automorphism
   $\nu$ of $\mathfrak{g}$ such that $\nu(H) \equiv K_i
   \ (m_i=1\ {\rm or} \ 2) \ ({\rm mod} \ 
   2\varPi(\mathfrak{g}_{\mathbb{C}},\mathfrak{t}_{\mathbb{C}}))$.
\vspace{0.5cm}\par

We write $h\sim k$ if $\tau_h$ is conjugate to
$\tau_k$ within the group of inner automorphism
of $\mathfrak{g}$.

\begin{lem}
  $(A_n)$ If $\mathfrak{g}$ is of type $A_n$, then $K_i\sim K_{n+1-i}$.\\
  $(D_n)$ If $\mathfrak{g}$ is of type $D_n$, then $K_i\sim K_{n-i} \
 (1\leq i\leq[n/2])$. 
  In particular if $n$ is odd, then $K_{n-1}\sim K_{n}$.\\
  $(E_6)$ If $\mathfrak{g}$ is of type $E_6$, then $K_1\sim K_6, \
 K_2\sim K_3\sim K_5$. 
   \label{lem:23}
\end{lem}

 Proof. $(A_n):$ We identify $\varDelta(\mathfrak{g}_{\mathbb{C}}, \mathfrak
h_{\mathbb{C}})$ with 
$$\{e_i - e_j \ ; \ 1\leq i\ne j\leq n+1\}$$
(for example, see \cite{Helg}), where $\{e_1 , \ldots ,e_{n+1}\}$ is an 
orthonormal basis of ${\mathbb{R}}^{n+1}$. From \cite{B} there exists 
an element $w$ of the Weyl group $W(\mathfrak{g}, \mathfrak{t})$ of 
$\mathfrak{g}$ with respect to $\mathfrak{t}$ such that 
$w(e_j)=e_{n-j+2} \ (1\leq j\leq n+1)$. Set $\alpha_i=e_i-e_{i+1}$. 
Then we have
$$w(\alpha_i)=w(e_i-e_{i+1})=e_{n-i+2}-e_{n-i+1}=-\alpha_{n-i+1}.$$
It is easy to see that $w^{-1}(K_i)=-K_{n+1-i}\equiv K_{n+1-i}\ 
(\mbox{mod}\ 2\varPi(\mathfrak{g}_{\mathbb{C}}, \mathfrak
h_{\mathbb{C}}))$. Hence $\tau_{w^{-1}(K_i)}=w^{-1}\circ\tau_{K_i} \circ w
=\tau_{K_{n+1-i}}$. \par 
$(D_n):$ 
$$\varDelta(\mathfrak{g}_{\mathbb{C}}, \mathfrak
h_{\mathbb{C}})=\{\pm e_i \pm e_j \ ; \ 1\leq i\ne j\leq n\}.$$
Set 
$$\alpha_i=e_i-e_{i+1}(1\leq i\leq n-1), \ \alpha_n=e_{n-1}+e_n.$$
Since there exists $w\in W(\mathfrak{g}, \mathfrak{t})$ such that 
$w(e_j)=e_{n-j+1} \ (1\leq j\leq n)$, we have  
$$w(\alpha_i)=w(e_i-e_{i+1})=e_{n-i+1}-e_{n-i}=-\alpha_{n-i}.$$
Hence we get $w^{-1}(K_i)=-K_{n-i}\equiv K_{n-i}\ 
({\rm mod}\ 2\varPi(\mathfrak{g}_{\mathbb{C}}, \mathfrak{h}_{\mathbb{C}}))$.  
In particular, if $n$ is odd, then there exists a unique 
$\bar w\in W(\mathfrak{g}, \mathfrak{t})$ such that $\{\alpha_1,\ldots
,\alpha_n\}\to\{-\alpha_1,\ldots,-\alpha_n\}$. If
$\bar w(\alpha_i)=-\alpha_i$ for $1\leq i\leq n$, then $\bar w=-{\rm Id}$, 
which is a contradiction (cf, \cite{B}). Thus we get
$$\bar w(\alpha_i)=-\alpha_i \ (1\leq i\leq n-2), \ \bar
w(\alpha_{n-1})=-\alpha_n, \ \bar w(\alpha_n)=-\alpha_{n-1}.$$ 
Hence we obtain $\bar w^{-1}(K_{n-1})=-K_n\equiv K_n \ ({\rm mod}\ 
2\varPi(\mathfrak{g}_{\mathbb{C}}, \mathfrak{h}_{\mathbb{C}}))$. 
 \par
$(E_6):$ There exists a unique $w\in W(\mathfrak{g}, \mathfrak{t})$ 
such that $\{\alpha_1,\ldots,\alpha_6\}\to\{-\alpha_1,\ldots,-\alpha_6\}$. 
Similarly as in the proof of $(D_n)$, we have
$$w(\alpha_1)=-\alpha_6, \ w(\alpha_2)=-\alpha_2, \
w(\alpha_3)=-\alpha_5, \ w(\alpha_4)=-\alpha_4.$$
Hence we obtain $w^{-1}(K_1)=-K_6$ and $w(K_3)=-K_5$. On the other hand, 
it is easy to see that
$t_{\alpha_1+\alpha_2+2\alpha_3+2\alpha_4+\alpha_5}\circ
t_{\alpha_2+\alpha_4+\alpha_5}(K_2)=-K_5+2K_6\equiv K_5\ 
({\rm mod}\ 2\varPi(\mathfrak{g}_{\mathbb{C}}, \mathfrak{h}_{\mathbb{C}}))$. 
Thus we have $K_2\sim K_5$.\hfill$\Box$\vspace{0.5cm}\par

 Let $(G/H,\langle,\rangle,\sigma )$ be a compact Riemannian 4-symmetric 
space such that $\sigma$ is inner. Then the following holds.
\begin{lem}[\cite{J}]
  $\sigma$ is conjugate within ${\rm Int} (\mathfrak{g})$ to some
  ${\rm Ad}(\exp(\pi/2)\sqrt{-1}h_a)$ where either   
   $$\begin{array}{lll}
    &h_0 = K_i, & m_i = 4, \\
    &h_1 = K_i \ \mbox{or}\ K_j+K_k,& m_i = 3, \ m_j=m_k=2, \\
    &h_2 = K_i + K_j, &  m_i = 1, m_j = 2,\\
    &h_3 = K_i + K_j + K_k, \ \ &  m_i = m_j = m_k =1, \\
    &h_4 = K_i, & m_i=1, \\
    &h_5 = K_i,\ K_j +K_k \ \mbox{or}\ 2K_p +K_q,  & m_i=2, 
    m_j=m_k=m_p=m_q=1.
 \end{array}$$
     \label{lem:24}
\end{lem}
\vspace{0.5cm}\par

{\sc Remark 2.2.} $(1)$ If $\sigma$ is conjugate to $\tau_{(1/2)h_4}$, then a 
pair $(\mathfrak{g},\mathfrak{g}^{\sigma})$ is symmetric. 
Indeed, for $\alpha=\sum_{r=1}^{n}k_r\alpha_r\in
\varDelta(\mathfrak{g}_{\mathbb{C}},\mathfrak{t}_{\mathbb{C}})$, we have 
$\alpha(h_4)=k_i$ and 
$$\alpha(h_4)\equiv0\!\!\!\pmod{4}\Longleftrightarrow
\alpha(h_4)\equiv0\!\!\!\pmod{2}\Longleftrightarrow k_i=0$$
since $m_i=1$. Therefore it follows that $\mathfrak{g}^{\tau_{(1/2){K_i}}}
=\mathfrak{g}^{\tau_{K_i}}$. Hence $(\mathfrak{g},\mathfrak{g}^{\sigma})$
is a symmetric pair, because $\tau_{K_i}$ is an involution. 
\par
If $\sigma$ is conjugate to $\tau_{(1/2)h_5}$, then a pair $(\mathfrak{g},
\mathfrak{g}^{\sigma})$ is $3$-symmetric. Indeed, for example, if 
$h_5=2K_p +K_q$, then we have 
$$\alpha(h_5)\equiv0\!\!\!\pmod{4}\Longleftrightarrow
\alpha(K_p+K_q)\equiv0\!\!\!\pmod{3}
\Longleftrightarrow k_p=k_q=0,
$$
for $\alpha=\sum_{r=1}^nk_r\alpha_r\in
\varDelta(\mathfrak{g}_{\mathbb{C}},\mathfrak{t}_{\mathbb{C}})$. Therefore, 
we obtain $\mathfrak{g}^{\tau_{(1/2)h_5}}=
\mathfrak{g}^{\tau_{(2/3)(K_p +K_q)}}$, and hence $(\mathfrak{g},
\mathfrak{g}^{\sigma})$ is a $3$-symmetric pair because 
$\tau_{(2/3)(K_p +K_q)}$ is of order $3$. 
\smallskip\par
$(2)$ Let $\mathfrak{z}$ be the center of $\mathfrak{h}$. If 
$\sigma={\rm Ad}(\exp(\pi/2)\sqrt{-1}h_a) \ (a=0,1,2,3)$, then the dimension 
of $\mathfrak{z}$ is equal to $a$ (\cite{J}).

\subsection{Graded Lie algebras}
 In this subsection we recall notions and some results on graded Lie 
algebras. 
\par
Let $\mathfrak{g}^*$ be a noncompact semisimple Lie algebra over $\mathbb
 R$. Let $\tau$ be a Cartan involution of $\mathfrak{g}^*$ and 
\begin{equation}
  \mathfrak{g^* =k + p^*} ,  \ \ \tau |_{\mathfrak k} = {\rm Id}_{\mathfrak k}
  ,  \ \tau |_{\mathfrak p^*} = -{\rm Id}_{\mathfrak p^*}
    \label{eqn:25}
\end{equation} 
the Cartan decomposition of $\mathfrak{g}^*$
 corresponding to $\tau$. Let $\mathfrak{a}$ be a maximal abelian
 subspace of $\mathfrak p^*$ and $\varDelta$ the set of restricted
 roots of $\mathfrak{g}^*$ with respect to $\mathfrak{a}$. We denote by
 $\varPi = \{\lambda_1, \ldots ,\lambda_l\}$ the set of fundamental roots of
 $\varDelta$ with respect to a lexicographic ordering of $\mathfrak{a}$. We
 call a collection of subsets $\{\varPi_0,\varPi_1,\ldots,\varPi_n\}$ of $\varPi$ a {\it
 partition} of $\varPi$ if $\varPi_1\ne\emptyset, \varPi_n\ne\emptyset$
 and

$$\varPi = \varPi_0\cup\varPi_1\cup\cdots\cup\varPi_n \ \ {\rm (disjoint
 \  union).}$$ 
Let $\varPi$ and $\bar{\varPi}$ be fundamental root systems of
 noncompact semisimple Lie algebras $\mathfrak{g}^*$ and $\bar{\mathfrak
 g^*}$ respectively. Partitions $\{\varPi_0,\varPi_1,\ldots,\varPi_m\}$
 of $\varPi$ and $\{\bar{\varPi_0},\bar{\varPi_1},\ldots,\bar{\varPi_n}\}$
 of $\bar{\varPi}$ are said to be {\it equivalent} if there exists an
 isomorphism $\phi$ from Dynkin diagram of $\varPi$ to that of $\bar{\varPi}$
 such that $m=n$ and $\phi(\varPi_i)=\bar{\varPi_i} \ (i=0,1,\ldots
 ,n)$.\par
  Take a gradation 
\begin{equation}
 \begin{array}{l}
   \mathfrak{g}^* = \mathfrak{g}_{-\nu}^* + \cdots + \mathfrak
   {g}_0^* + \cdots + \mathfrak{g}_{\nu}^* , \ \ \ 
 \vspace{0.1cm}\\
   \left[\mathfrak {g}_p^*,\mathfrak{g}_q^*\right] \subset\mathfrak
   {g}_{p+q}^*, \ \ \ 
   \tau (\mathfrak {g}_p^*) = \mathfrak {g}_{-p}^*, \ \ \ -\nu \leq p, \ q
   \leq \nu, 
 \end{array}
\end{equation}
of $\nu$-th kind on $\mathfrak{g}^*$ so that $\mathfrak {g}_1^* \ne \{0\}$. 
We denote by $Z$ the charactiristic element of the gradation, i.e. $Z$
 is a unique element in $\mathfrak {p}^* \cap \mathfrak {g}_0^*$ such
 that 
$$\mathfrak{g}_p^* = \{ X\in\mathfrak{g}^* \ ; \ [Z,X] = pX\}, \ \ \ -\nu\leq
 p \leq\nu .$$
Let
$$\mathfrak{g}^* = \sum _{i=-\nu}^\nu \mathfrak {g}_i^* , \ \ \
 \bar{\mathfrak{g}}^* = \sum _{i=-\bar{\nu}}^{\bar{\nu}} \bar{\mathfrak
 {g}}_i ^* $$
be two graded Lie algebras. These gradations are said to be {\it isomorphic}
 if $\nu = \bar{\nu}$ and there exists an isomorphism $\phi : \mathfrak
 {g}^* \to \bar{\mathfrak {g}}^*$ such that $\phi(\mathfrak {g}_i^*) =
 \bar{\mathfrak {g}}_i^*(-\nu\leq i \leq\nu)$. Then the following holds.

\begin{thm}[Kaneyuki and Asano {\cite{KA}}]
   Let $\mathfrak{g}^*$ be a noncompact semisimple Lie algebra over
   $\mathbb{R}$ and $\varPi$ a fundamental root system of $\mathfrak
   g^*$. Then there exists a bijection between the set of equivalent
   classes of partitions of $\varPi$ and set of isomorphic classes of
   gradations on $\mathfrak{g}^*$. 
    \label{thm:21}
\end{thm}
\noindent
The bijection in the theorem is constructed as follows: Let
$\{\varPi_0,\varPi_1,\ldots,\varPi_n\}$ be a partition of
$\varPi$. Define $h_{\varPi} : \varDelta \to \mathbb Z$ by 
$$h_{\varPi} (\lambda) := \sum_{\lambda_i \in \varPi_1}m_i +
2\sum_{\lambda_j \in \varPi_2}m_j + \cdots + n\sum_{\lambda_k \in
\varPi_n}m_k, \ \ \ \lambda = \sum_{i=1}^l m_i \lambda_i \in \varDelta.$$
Then there is a unique $Z$ in $\mathfrak{a}$ such that $\lambda (Z) =
h_{\varPi} (\lambda)$ for all $\lambda \in \varDelta$. For a partition
$\{\varPi_0,\varPi_1,\ldots,\varPi_n\}$ we obtain a gradation
$\mathfrak{g}^* = \sum_{i=-\nu}^{\nu} \mathfrak{g}_i^*$ whose
charactiristic element equals $Z$. This correspondence induces a
bijection mentioned in the theorem.\par
Define $h_i\in\mathfrak{a}$ ($i=1,2,\cdots,l$) by 
$$\lambda_j(h_i)=\delta_{ij}.$$
 Let $\mathfrak{t}^*$ be a Cartan subalgebra 
of $\mathfrak{g}^*$ such that $\mathfrak{a}\subset\mathfrak{t}^*$. Take 
compatible orderings on $\mathfrak{t}^*$ and $\mathfrak{a}$. 
We clarify the relation between $K_i$ and $h_j$.
\begin{lem} Let $\lambda_i$ be any root in $\varPi$.
\par
$(1)$ If there exists a unique $\alpha_j\in\varPi(\mathfrak{g}^*_{\mathbb{C}},
\mathfrak{t}^*_{\mathbb{C}})$ such 
that $\alpha_j|_{\mathfrak{a}}=\lambda_i$, then $h_i=K_j$. 
\par
$(2)$ If there exist two fundamental roots 
$\alpha_j$, $\alpha_k\in\varPi(\mathfrak{g}^*_{\mathbb{C}},
\mathfrak{t}^*_{\mathbb{C}})$ such that 
$\alpha_j|_{\mathfrak{a}}=\alpha_k|_{\mathfrak{a}}=\lambda_i$, then 
$h_i=K_j+K_k$.
\label{lem:25} 
\end{lem}

 Proof. (1): Considering the classification of the Satake diagrams, 
for $\alpha_p\in\varPi(\mathfrak{g}^*_{\mathbb{C}},
\mathfrak{t}^*_{\mathbb{C}})$, 
$p\not=j$, it follows that $\alpha_p|_{\mathfrak{a}}=0$ or 
$\alpha_p|_{\mathfrak{a}}=\lambda_q$ for some $q$ ($q\not=i$). 
Thus we have 
$$\alpha_p(h_i)=\alpha_p|_{\mathfrak{a}}(h_i)=0,\ \ \alpha_j(h_i)
=\lambda_i(h_i)=1,$$
which implies $h_i=K_j$. 
\par
(2): Similarly as above, for $\alpha_p\in\varPi(\mathfrak{g}^*_{\mathbb{C}},
\mathfrak{t}^*_{\mathbb{C}})$, $p\not=j,k$, 
it follows that $\alpha_p|_{\mathfrak{a}}=0$ or $\alpha_p|_{\mathfrak{a}}
=\lambda_q$ for some $q$ ($q\not=i$). 
Therefore 
$$\alpha_p(h_i)=\alpha_p|_{\mathfrak{a}}(h_i)=0,\ \ \alpha_m(h_i)
=\alpha_m|_{\mathfrak{a}}(h_i)=\lambda_i(h_i)=1, \ \ m=j,k,
$$
which implies $h_i=K_j+K_k$. \hfill $\Box$


\section{Riemannian 4-symmetric spaces}
 In this section we use the same notation as in Section 2. Let $(G/H,
\langle,\rangle, \sigma)$ be a Riemannian 4-symmetric space with an inner 
automorphism $\sigma$ of order $4$.  Let $\mathfrak{g}$ and
$\mathfrak{h}$ be the Lie algebras of $G$ and $H$, respectively. Note
that $\mathfrak{h}$ coincides with the set $\mathfrak{g}^{\sigma}$ of
fixed points of $\sigma$. Choose a subspace $\mathfrak{m}$ of
$\mathfrak{g}$ so that $\mathfrak{g = h + m}$ is an ${\rm Ad}(H)$- and
$\sigma$-invariant decomposition.  Let $\mathfrak{t}$ be a maximal
abelian subalgeba of $\mathfrak{g}$  contained in $\mathfrak h$, and
$\mathfrak{z}$ the center of $\mathfrak{h}$. \par  
Suppose that $\mathfrak{g}$ is a compact simple Lie algebra. 
Let ${\rm Aut}_{\mathfrak{h}}(\mathfrak{g})$ be the set of 
automorphisms of $\mathfrak{g}$ preserving $\mathfrak{h}$.  
\begin{lem}
  Assume $\sigma={\rm Ad}(\exp({\pi}/2)\sqrt{-1}K_i), 
  \ m_i =3 \ or \ 4$, where $\delta = \sum_{j=1}^n m_j \alpha_j$ 
is the highest root of $\varDelta(\mathfrak{g}_{\mathbb{C}},
\mathfrak{t}_{\mathbb{C}})$ as in Section 2. Then for each 
$\mu\in{\rm Aut}_{\mathfrak{h}}(\mathfrak{g})$, we have 
$\mu\circ\sigma\circ\mu^{-1}=\sigma$ or $\sigma^{-1}$. 
    \label{lem:31}
\end{lem}

 Proof. 
Since $\mu(\mathfrak h)=\mathfrak h$, we obtain
$\mathfrak{g} ^{\tilde{\sigma}} =\mathfrak{h}$, where $\tilde{\sigma} :=
\mu\circ\sigma\circ\mu^{-1}$. In particular, we have 
$\tilde{\sigma}|_{\mathfrak{t}}={\rm Id}_{\mathfrak{t}}$. Therefore, it follows 
from Proposition 5.3 of Chapter IX of \cite{Helg} that there is 
$\sqrt{-1}Z\in \mathfrak{t}$ such that 
\begin{equation}
  \tilde{\sigma} = {\rm Ad}(\exp\frac{\pi}{2}\sqrt{-1}Z).
    \label{eqn:31}
\end{equation}
Since $\sigma={\rm Ad}(\exp(\pi/2)\sqrt{-1}K_i)$ with $m_i=3$ or $4$, we obtain 
$E_{\alpha_j}\in{\mathfrak{h}}_{\mathbb{C}}$ ($j\not=i$) and 
$E_{\alpha_i}\not\in{\mathfrak{h}}_{\mathbb{C}}$. Moreover, since 
$\mathfrak{g}^{\sigma}=\mathfrak{g}^{\tilde{\sigma}}=\mathfrak h$, it follows 
from (\ref{eqn:31}) that
\begin{equation}
\tilde{\sigma}(E_{\alpha_j}) = E_{\alpha_j}, \ \ 
\tilde{\sigma}(E_{\alpha_i}) = c E_{\alpha_i}, 
\label{eqn32}\end{equation}
for some $c\in\mathbb{C}$ with $|c|=1$. Then $c^4=1$ and $c^2\not=1$, because 
$\tilde{\sigma}^4={\rm Id}$ and $\tilde{\sigma}^2\not={\rm Id}$. From 
(\ref{eqn32}), we can see that if $c=\sqrt{-1}$, then $\tilde{\sigma}=\sigma$,
and if $c=-\sqrt{-1}$, then $\tilde{\sigma}=\sigma^{-1}$.\hfill$\Box$
\vspace{0.5cm}\par

{\sc Remark 3.1.} 
Lemma \ref{lem:31} dose not hold in general. If $\sigma$ is
 conjugate to ${\rm Ad}(\exp(\pi/2)\sqrt{-1}(K_i +
 K_j)) (m_i=m_j=2)$, then Lemma \ref{lem:31} holds. However in other cases, 
Lemma \ref{lem:31} dose not hold.  
\vspace{0.5cm}\par

{\sc Remark 3.2.}
If $\sigma$ is an automorphism of order 2 or 3, then by an argument similar to 
the proof of Lemma \ref{lem:31}, it follows that
 $\mu\circ\sigma\circ\mu^{-1}=\sigma$ for any $\mu\in{\rm
 Aut}_{\mathfrak{h}}(\mathfrak{g})$. 
\vspace{0.5cm}\par

\begin{lem}
  Suppose that $\sigma={\rm Ad}(\exp(\pi/2)\sqrt{-1}K_i)$ with $m_i
 =3$. Let $\tau$ be an involutive automorphism of $\mathfrak{g}$ such that
 $\tau (\mathfrak h) = \mathfrak h$. Then
   \begin{itemize}
    \item[{\rm (i)}] $\tau\circ\sigma = \sigma\circ\tau$ if and only if the
		 coefficient of $\alpha_i$ in $\tau(\delta)$ is equal to
		 $3$. 
    \item[{\rm (ii)}] $\tau\circ\sigma = \sigma^{-1}\circ\tau$ if and only if the
		 coefficient of $\alpha_i$ in $\tau(\delta)$ is equal to
		 $-3$. 
   \end{itemize}
     \label{lem:32}
\end{lem}

 Proof. It is easy to see that $\mathfrak{z} = \mathbb{R}\sqrt{-1}K_i$
 for some $i$ with $m_i =3$. Since $\tau 
(\mathfrak{h}) =\mathfrak{h}$, we have $\tau(\sqrt{-1} K_i) =
\pm\sqrt{-1}K_i$, and therefore
$$\tau\circ\sigma\circ\tau^{-1} = {\rm Ad}(\exp\frac{\pi}{2}\tau(\sqrt{-1}K_i))
=\left\{
\begin{array}{ll}
\sigma & \ \ {\rm if} \ \tau(\sqrt{-1}K_i)=\sqrt{-1}K_i, \\
\sigma^{-1} & \ \ {\rm if} \ \tau(\sqrt{-1}K_i)=-\sqrt{-1}K_i,
\end{array}
\right.$$ 
and 
$$\tau(\delta)(K_i) = \delta(\tau(K_i)) = \delta(\pm K_i) =
\left\{
\begin{array}{rl}
3 & \ \ {\rm if} \ \tau(\sqrt{-1}K_i)=\sqrt{-1}K_i, \\
-3 & \ \ {\rm if} \ \tau(\sqrt{-1}K_i)=-\sqrt{-1}K_i.
\end{array}
\right.$$
This completes the proof of tha lemma. \hfill$\Box$
\vspace{0.5cm}\par

\begin{lem} Suppose that $\sigma={\rm Ad}(\exp(\pi/2)\sqrt{-1}K_i)$ with 
 $m_i=3$ or $4$. \par
\noindent
  {\rm (i)} \  Let $\tau_1, \tau_2$ be involutive automorphisms of 
  $\mathfrak{g}$ such that $\tau_i (\mathfrak h) = \mathfrak h, (i = 1,2)$. 
  If there exists $\mu\in {\rm Aut_{\mathfrak h} (\mathfrak{g})}$ such that
 $\mu\circ\tau_1\circ\mu^{-1} = \tau_2$. Then 
     \begin{equation}
       \mathfrak{g} ^{\tau_1} \cong \mathfrak{g} ^{\tau_2}, \ \ \mathfrak
       {h\cap g}^{\tau_1}\cong \mathfrak {h\cap g}^{\tau_2}.
     \end{equation}
  {\rm (ii)} \ Put $\tau' := \mu\circ\tau\circ\mu^{-1}, \ \mu\in{\rm
 Aut_{\mathfrak h}(\mathfrak{g})}$. If $\tau\circ\sigma = \sigma^{\pm
 1}\circ\tau$, then $\tau'\circ\sigma = \sigma^{\pm 1}\circ\tau'$,
 respectively. 
    \label{lem:34}
\end{lem}

 Proof. (i) is trivial. \\
(ii) We have
 $$\begin{array}{lll}
   \tau\circ\sigma = \sigma^{\pm
   1}\circ\tau&\Longleftrightarrow&\mu\circ\tau\circ\sigma\circ\mu^{-1} =
   \mu\circ\sigma^{\pm 1}\circ\tau\circ\mu^{-1}\\
   &\Longleftrightarrow&\tau'\circ\mu\circ\sigma\circ\mu^{-1} =
   \mu\circ\sigma^{\pm 1}\circ\mu^{-1}\circ\tau'
\end{array}$$
Hence, it follows from Lemma 3.1 that if 
$\tau\circ\sigma=\sigma^{\pm 1}\circ\tau$, then $
\tau'\circ\sigma=\sigma^{\pm 1}\circ\tau'$. \hfill$\Box$ 
\vspace{0.5cm}\par
In the remaining part of this paper, we suppose that 
$\sigma={\rm Ad}(\exp({\pi}/2)\sqrt{-1}K_i)$ for some $\alpha_i\in
\varPi(\mathfrak{g}_{\mathbb{C}},\mathfrak{t}_{\mathbb{C}})$ with $m_i =3$ 
or $4$. If $m_i =  3$, the Dynkin diagram of
 $\mathfrak h$ is isomorphic to the extended Dynkin diagram of $\varPi
 (\mathfrak{g}_{\mathbb{C}},\mathfrak{t}_{\mathbb{C}})$ except $\alpha_i$
 and $\alpha_0$, and if $m_i =4$, it is isomorphic to that of $\varPi
 (\mathfrak{g}_{\mathbb{C}},\mathfrak{t}_{\mathbb{C}})$ except
 $\alpha_i$ (cf. Theorem 5.15 of Chapter X of \cite{Helg}). 
We denote by $\varPi(\mathfrak h)$ the fundamental root system of $\mathfrak{h}$ 
corresponding to the Dynkin diagram of $\mathfrak h$. 
\begin{lem}
   For any involutive
   automorphism $\tau$ of $\mathfrak{g}$ satisfying $\tau(\mathfrak h) =
   \mathfrak h$, there exists $\mu\in{\rm Int}(\mathfrak h)$
   such that $\mu\circ\tau\circ\mu^{-1}(\varPi(\mathfrak h)) =
 \varPi(\mathfrak{h})$.
     \label{lem:35}
\end{lem}

 Proof.  
Put $\tilde{\tau}: = \tau |_{\mathfrak h}$. Then $\tilde{\tau}$ is an
involution of $\mathfrak h = \mathfrak z \oplus\mathfrak h _s$, where
$\mathfrak h _s : = [\mathfrak{h , h}]$. 
It is obvious that $\tilde{\tau}(\mathfrak z) = \mathfrak z$ and
$\tilde{\tau}(\mathfrak h _s) = \mathfrak h _s$. 
Decompose $\mathfrak{h}_s$ into $\mathfrak{h}_s = \mathfrak{h}_1
 \oplus \cdots \oplus \mathfrak{h}_m$ where $\mathfrak{h}_1, \ldots ,
 \mathfrak{h}_m$ are simple ideals. From the classification of compact $4$-symmetric pairs 
$(\mathfrak{g},\mathfrak{h})$ (\cite{J}), it follows that 
(i) $\mathfrak{h} _i \ncong \mathfrak{h} _j$ for 
any $i,j\in\{1,\ldots ,m\} \ (i\not=j)$, or (ii) there exists only one
 pair $(p,q)$ such that $\mathfrak{h}_p \cong \mathfrak{h}_q$. 
\vspace{0.2cm}
\par

{\it The case} (i): Since $\tau(\mathfrak{h}_i)$ ($1\leq i\leq m$) is 
a simple ideal of $\mathfrak{h}_s$ and $\mathfrak{h}_i\not\cong\mathfrak{h}_j$ 
($i\not=j$), it follows that $\tilde{\tau}(\mathfrak{h}_i)=\mathfrak{h}_i$ 
($1\leq i\leq m$). 
Therefore we have a direct sum decomposition $\mathfrak h _i =
\mathfrak k _i \oplus \mathfrak p _i$. Let $\mathfrak{a} _i$ be a maximal
abelian subspace of $\mathfrak p _i$ and $\mathfrak{t} _i$ be a maximal
abelian subalgebra of $\mathfrak h _i$ containing $\mathfrak{a} _i$. We
take a fundamental root system $\varPi_i = \{\lambda_1, \ldots ,
\lambda_{n_i}\}$ for the set of nonzero roots with respect to
$({\mathfrak h _i}_{\mathbb{C}}, {\mathfrak{t} _i}_{\mathbb{C}})$. From 
Theorem 5.15 of \cite{Helg}, there exists $\mu_i\in{\rm Aut}
(\mathfrak{h}_i)$ such that $\mu_i\circ\tau|_{\mathfrak{h}_i}\circ\mu_i^{-1}$ 
is an automorphism of $\varPi_i$ of order $1$ or $2$. Hence we have 
\begin{equation}
\tau(\mu_i^{-1}(\mathfrak{t}_i))=\mu_i^{-1}(\mathfrak{t}_i),\ \ 
\tau(\mu_i^{-1}(\varPi_i))=\mu_i^{-1}(\varPi_i).
\label{eqn:34}
\end{equation}
Set $\tilde{\mathfrak{t}}:=\mu_1^{-1}(\mathfrak{t}_1)\oplus\cdots
\oplus\mu_m^{-1}(\mathfrak{t}_m)\oplus\mathfrak{z}$ and 
$\tilde{\varPi}:=\mu_1^{-1}(\varPi_1)\cup\cdots\cup\mu_m^{-1}(\varPi_m)$. 
Then by (\ref{eqn:34}), we have $\tau(\tilde{\mathfrak{t}})
=\tilde{\mathfrak{t}}$ and $\tau(\tilde{\varPi})=\tilde{\varPi}$. Since there 
exist $\mu\in{\rm Int}(\mathfrak{h})$ and $w\in W(\mathfrak{h},\mathfrak{t})$ 
such that $\mu(\tilde{\mathfrak{t}})=\mathfrak{t}$ and 
$w(\mu(\tilde{\varPi}))=\varPi(\mathfrak{h})$, we obtain 
$$
(w\mu)\circ\tau\circ(w\mu)^{-1}(\varPi(\mathfrak{h}))=\varPi(\mathfrak{h}),
$$
which completes the proof of the lemma for the case (i). 
\vspace{0.2cm}
\par 
{\it The case} (ii): If $\tau(\mathfrak{h} _i) = \mathfrak{h} _i$ for $i =
1,\ldots ,m$, then by the same argument as in the case (i), we can prove the
claim. Hence we assume that and $\tau(\mathfrak{h} _p) = \mathfrak{h}
_q$ and $\tau(\mathfrak{h}_i) = \mathfrak{h}_i$ for $i\ne p , q$. 
Define isomorphisms $\tau_1:\mathfrak{h}_q\longrightarrow\mathfrak{h}_p$ and 
$\tau_2:\mathfrak{h}_p\longrightarrow\mathfrak{h}_q$ by
$$
\tau (X, Y) = (\tau_1 (Y), \tau_2(X)), \ \ X\in\mathfrak{h}_p,\ Y\in
\mathfrak{h}_q.
$$
Since $\tau$ is an involution, it follows that 
$\tau_1 \circ \tau_2 =\tau_2 \circ \tau_1 = {\rm Id}$. 
Hence we have $\tau(X, Y) = (\tau_1(Y), \tau_1^{-1}(X))$. \par

Put $\mathfrak{b}:=\mathfrak{h}_p$ and define an isomorphism 
$\nu:\mathfrak{h}_p\oplus\mathfrak{h}_q\longrightarrow
\mathfrak{b}\oplus\mathfrak{b}$ by 
$$\nu (X,Y):=(X,\tau_1(Y)).$$
Then it is easy to see that $\nu\circ\tau\circ\nu^{-1}(X,Y) = (Y,X)$.
Therefore, considering a symmetric pair $(\mathfrak{b}\oplus\mathfrak{b},
\Delta\mathfrak{b})$ ($\Delta\mathfrak{b}:=\{(X,X) \ ; \ X\in\mathfrak{b}\}$), 
we can see that there exist a fundamental root system of 
$\mathfrak{h}_p\oplus\mathfrak{h}_q$ preserved by $\tau|_{\mathfrak{h}_p
\oplus\mathfrak{h}_q}$. Hence, by an argument similar to (i), there exists 
$\nu\in{\rm Int}(\mathfrak{h})$ such that $\nu\circ\tau\circ\nu^{-1}
(\varPi(\mathfrak{h}))=\varPi(\mathfrak{h})$. 
This completes the proof of the lemma for the case
(ii). \hfill$\Box$\vspace{0.3cm}\par 

In the following sections, we shall classify the equivalence classes of 
involutive automorphisms $\tau$ within ${\rm Aut}_{\mathfrak h}(\mathfrak{g})$
of $\mathfrak{g}$  such that $\tau(\mathfrak h) = \mathfrak h$. From 
Remark 2.2 and Lemma \ref{lem:31} we have the following 
four type:
$$
\begin{array}{l}
  \dim\mathfrak z =0, \ \tau\circ\sigma = \sigma^{\pm 1}\circ\tau, \\
  \dim\mathfrak z =1, \ \tau\circ\sigma = \sigma^{\pm 1}\circ\tau. \\
\end{array}
$$


\section{The case where $\dim\mathfrak{z}=0$}
 In the remaining part of this paper we use the same notation as in 
Section 2 and Section 3. Let $(G/H,\langle,\rangle,\sigma)$ be a Riemannian 
$4$-symmetric space such that $\sigma$ is inner and $\dim
\mathfrak{z}=0$. From Lemma \ref{lem:24} together with Remark 2.2 we may
suppose that  
$$\sigma = {\rm Ad}(\exp\frac{\pi}{2}\sqrt{-1}K_i) \ {\rm  \ for  \
some} \ i \ {\rm with \ the \ property} \ m_i = 4.$$
According to Section 3 and Jim\'{e}nez \cite {J},
$4$-symmetric pairs $(\mathfrak {g,h})$ satisfying the condition $\dim
\mathfrak z =0$ are given by
\begin{equation}
 \begin{array}{l}
   (\mathfrak{e} _7, \mathfrak {so}(6) \oplus \mathfrak {so}(6)  \oplus
   \mathfrak{su}(2)), \ (\mathfrak{e} _8, \mathfrak {su}(8)  \oplus
   \mathfrak{su}(2)), \\
   (\mathfrak{e} _8, \mathfrak {so}(10)  \oplus \mathfrak{so}(6)), \
   (\mathfrak{f} _4, \mathfrak {so}(6)  \oplus \mathfrak {so}(3)).
 \end{array}
\end{equation}
\par
Let $\tau$ be an involution of $\mathfrak{g}$ preserving $\mathfrak h$. 
By Lemma \ref{lem:35}, we may assume $\tau(\mathfrak{t})=\mathfrak{t}$ and 
$\tau(\varPi(\mathfrak h))=\varPi(\mathfrak h)$. 
If $\tau|_{\mathfrak{t}} = {\rm Id}_{\mathfrak{t}}$, then there exists 
$\sqrt{-1}H\in\mathfrak{t}$ such that $\tau={\rm Ad}(\exp\pi\sqrt{-1}H)$ 
and $\tau\circ\sigma=\sigma\circ\tau$. 
\par
Now, we assume $\tau |_{\mathfrak{t}}\ne{\rm Id}_{\mathfrak
t}$. Suppose that $\mathfrak{g}$ is of type $\mathfrak{e} _8$. From
Section 3, the Dynkin diagram of $\mathfrak{h}$ coincides with the 
extended Dynkin diagram of $\mathfrak{e}_8$ except $\oplus$ as follows: 
\begin{equation}
\setlength{\unitlength}{0.7pt}
\begin{picture}(270,50)
 \put(10,13){(i)}
 \put(55,15){\circle{8}}
 \put(52,0){\small $\alpha_0$}
 \put(59,15){\line(1,0){20}}
 \put(83,15){\circle{8}}
 \put(80,0){\small $\alpha_8$}
 \put(87,15){\line(1,0){20}}
 \put(111,15){\circle{8}}
 \put(108,0){\small $\alpha_7$}
 \put(115,15){\line(1,0){20}}
 \put(139,15){\circle{8}}
 \put(136,0){\small $\alpha_6$}
 \put(143,15){\line(1,0){20}}
 \put(167,15){\circle{8}}
 \put(164,0){\small $\alpha_5$}
 \put(171,15){\line(1,0){20}}
 \put(195,15){\circle{8}}
 \put(192,0){\small $\alpha_4$}
 \put(195,19){\line(0,1){20}}
 \put(195,43){\circle{8}}
 \put(204,43){\small $\alpha_2$}
 \put(199,15){\line(1,0){20}}
 \put(223,15){\circle{8}}
 \put(219,15){\line(1,0){8}}
 \put(223,11){\line(0,1){8}}
 \put(220,0){\small $\alpha_3$}
 \put(227,15){\line(1,0){20}}
 \put(251,15){\circle{8}}
 \put(248,0){\small $\alpha_1$}
\end{picture} 
\setlength{\unitlength}{0.7pt}
\begin{picture}(300,50)
 \put(10,15){(ii)}
 \put(55,15){\circle{8}}
 \put(52,0){\small $\alpha_0$}
 \put(59,15){\line(1,0){20}}
 \put(83,15){\circle{8}}
 \put(80,0){\small $\alpha_8$}
 \put(87,15){\line(1,0){20}}
 \put(111,15){\circle{8}}
 \put(108,0){\small $\alpha_7$}
 \put(115,15){\line(1,0){20}}
 \put(139,15){\circle{8}}
 \put(135,15){\line(1,0){8}}
 \put(139,11){\line(0,1){8}}
 \put(136,0){\small $\alpha_6$}
 \put(143,15){\line(1,0){20}}
 \put(167,15){\circle{8}}
 \put(164,0){\small $\alpha_5$}
 \put(171,15){\line(1,0){20}}
 \put(195,15){\circle{8}}
 \put(192,0){\small $\alpha_4$}
 \put(195,19){\line(0,1){20}}
 \put(195,43){\circle{8}}
 \put(204,43){\small $\alpha_2$}
 \put(199,15){\line(1,0){20}}
 \put(223,15){\circle{8}}
 \put(220,0){\small $\alpha_3$}
 \put(227,15){\line(1,0){20}}
 \put(251,15){\circle{8}}
 \put(248,0){\small $\alpha_1$}
\end{picture}
\label{eqn:42}
\end{equation}
We denote $\sum_{i=1}^8k_i\alpha_i$ by
$$\left(\begin{array}{lllllll}
    &     &     &     & k_2 &     &     \\
k_8 & k_7 & k_6 & k_5 & k_4 & k_3 & k_1 \\
\end{array}\right).$$
In the above case (i), since $\tau|_{\mathfrak{t}}\ne{\rm Id}_{\mathfrak{t}}$ 
and $\tau(\varPi(\mathfrak h))=\varPi(\mathfrak h)$, the possibility of 
$\tau|_{\mathfrak{t}}$ is as follows: 
$$\tau(\alpha_1) = \alpha_1, \ \tau(\alpha_2) = \alpha_0, \
\tau(\alpha_4) = \alpha_8, \ \tau(\alpha_5) = \alpha_7, \ 
\tau(\alpha_6) = \alpha_6.$$
Then we get
$$\alpha_2 = \tau(\alpha_0)=-4\tau(\alpha_3)-3\tau(\alpha_2)-
\left(\begin{array}{lllllll}
  &   &   &   & 0 &   &   \\
6 & 5 & 4 & 3 & 2 & 0 & 2 \\
\end{array}\right),$$
and hence 
$$\tau(\alpha_3)=
\left(\begin{array}{lllllll}
  &   &   &   & 2 &   &   \\
0 & 1 & 2 & 3 & 4 & 3 & 1 \\
\end{array}\right)\in\varDelta(\mathfrak{g}_{\mathbb{C}},\mathfrak
t_{\mathbb{C}}).$$\par
By a similar argument as above, we obtain the following proposition. 
\begin{prop} Suppose that $\dim\mathfrak{z}=0$. Let $\tau$ be an involution 
of $\mathfrak{g}$ such that 
$\tau(\mathfrak{h})=\mathfrak{h}$ and $\tau(\varPi(\mathfrak{h}))=
\varPi(\mathfrak{h})$.  Then all the possibilities of $\tau|_{\mathfrak{t}}$ 
such that $\tau|_{\mathfrak{t}}\not={\rm Id}_{\mathfrak{t}}$ 
are given by Table 1. 
\label{prop:41}
\end{prop}
{\small
\begin{center}
\renewcommand{\arraystretch}{1.2}
\begin{longtable}{c|c|c|c}
\caption{The possibilities of $\tau|_{\mathfrak{t}}$ such that 
$\tau|_{\mathfrak{t}}\ne{\rm Id}$ 
($\sigma =\tau_{(1/2)H})$).}\\ 
\hline
Type & $\mathfrak{g}$ & $H$ & $\tau|_{\mathfrak{t}}$ 
 \\ 
\hline
\hline
I & $\mathfrak{e}_7$ & $K_4$ &
 $\alpha_1\mapsto \alpha_6, \ \alpha_2\mapsto \alpha_2, \ \alpha_3\mapsto
 \alpha_5, \ \alpha_4\mapsto\alpha_4, \ \alpha_7\mapsto \alpha_0$  \\   
\hline
\lw{II} & \lw{$\mathfrak{e}_7$} & \lw{$K_4$} &
 $\alpha_1\mapsto \alpha_1, \ \alpha_2\mapsto \alpha_2, \ \alpha_3\mapsto
 \alpha_0, \ \alpha_5\mapsto\alpha_7, \  \alpha_6\mapsto\alpha_6$  \\   
&&&$\alpha_4\mapsto\alpha_1 + \alpha_2 + 2\alpha_3 + 3\alpha_4 +
 2\alpha_5 + \alpha_6$ \\
\hline
\lw{III} & \lw{$\mathfrak{e}_7$} & \lw{$K_4$} &
 $\alpha_1\mapsto \alpha_6, \ \alpha_2\mapsto \alpha_2, \ \alpha_3\mapsto
 \alpha_7, \ \alpha_5\mapsto \alpha_0$  \\
&&&$\alpha_4\mapsto\alpha_1 + \alpha_2 + 2\alpha_3 + 3\alpha_4 +
 2\alpha_5 + \alpha_6$ \\
\hline
\lw{IV} & \lw{$\mathfrak{e} _8$} & \lw{$K_3$} &
 $\alpha_1\mapsto \alpha_1, \ \alpha_2\mapsto \alpha_0, \ \alpha_4\mapsto
 \alpha_8, \ \alpha_5\mapsto\alpha_7, \ \alpha_6\mapsto\alpha_6$ \\ 
&&&$\alpha_3\mapsto\alpha_1 + 2\alpha_2 + 3\alpha_3 + 4\alpha_4 +
 3\alpha_5 + 2\alpha_6 + \alpha_7$ \\
\hline
\lw{V} & \lw{$\mathfrak{e}_8$} & \lw{$K_6$} &
 $\alpha_1\mapsto \alpha_1, \ \alpha_2\mapsto \alpha_5, \ \alpha_3\mapsto
 \alpha_3, \ \alpha_4\mapsto\alpha_4, \ \alpha_7\mapsto \alpha_0 \
 \alpha_8\mapsto\alpha_8$  \\   
&&&$\alpha_6\mapsto\alpha_1 + \alpha_2 + 2\alpha_3 + 3\alpha_4 +
 3\alpha_5 + 3\alpha_6 + 2\alpha_7 +\alpha_8$ \\
\hline
\lw{VI} & \lw{$\mathfrak{f}_4$} & \lw{$K_3$} &
 $\alpha_1\mapsto \alpha_1, \ \alpha_2\mapsto \alpha_0, \ \alpha_4\mapsto
 \alpha_4$ \\
&&&$\alpha_3\mapsto\alpha_1 + 2\alpha_2 + 3\alpha_3 + \alpha_4$ \\
\hline
\end{longtable}
\end{center}}
For Type IV in Table 1, it is easy to see $\tau(K_3)=-4K_2 + 3K_3\equiv 
-K_3\pmod 4$. Hence we have $\tau\circ\sigma =\sigma^{-1}\circ\tau$. 
Similarly, for Type I we get $\tau\circ\sigma=\sigma\circ\tau$ and 
for the other types, we have 
$\tau\circ\sigma = \sigma^{-1}\circ\tau$.\par
Finally, in order to compute the dimension of $\mathfrak{g}^{\tau}$, we
prove the following Lemma.
\begin{lem} Let $\mathfrak{t}_+$ be the $(+1)$-eigenspace of 
$\tau|_{\mathfrak{t}}$. Then 
 \begin{eqnarray*}
\dim\mathfrak{g}^{\tau}&\!\!\!=&\!\!\!\dim\mathfrak{t}_+ 
+\#\varDelta^+(\mathfrak{g}_{\mathbb{C}},\mathfrak{t}_{\mathbb{C}}) +
     2\#\{\alpha\in\varDelta^+(\mathfrak{g}_{\mathbb{C}},
\mathfrak{t}_{\mathbb{C}}) \ ; \ \tau(E_{\alpha}) = E_{\alpha} \}\\
&\!\!\! &\!\!\! -\#\{\alpha\in\varDelta^+(\mathfrak{g}_{\mathbb{C}},
\mathfrak{t}_{\mathbb{C}}) \ ; \ \tau(\alpha) = \alpha \}.
\end{eqnarray*}
         \label{lem:41}
\end{lem}

 Proof. If $\tau(\alpha)=\beta, \ (\beta\ne\pm\alpha)$, 
we can put $\tau(E_{\alpha}) = cE_{\beta}$ for some $c$. Since $\tau$ is 
involutive and $\tau(H_{\alpha})=H_{\beta}$, it is easy to
see that $E_{\alpha}+ cE_{\beta}$ and $E_{-\alpha}+c^{-1}E_{-\beta}$ are 
$(+1)$-eigenvectors of $\tau$. If $\tau(\alpha)=\alpha$, we get
$\tau(E_{\alpha})=E_{\alpha}$ or $\tau(E_{\alpha})=-E_{\alpha}$. Furthermore, 
if $\tau(\alpha)=-\alpha$, we can put $\tau(E_{\alpha})=cE_{-\alpha}$ for
some $c$. Then we have $\tau(E_{\alpha}\pm
cE_{-\alpha})=\pm(E_{\alpha}\pm cE_{-\alpha})$. Therefore we obtain
$$\begin{array}{lll}
  \dim\mathfrak{g}^{\tau} &\!\!\! = &\!\!\!\dim\mathfrak{t}_+ +
   \#\{\alpha\in\varDelta^+(\mathfrak{g}_{\mathbb{C}},
   \mathfrak{t}_{\mathbb{C}}) \ ; \ \tau(\alpha) \ne\pm\alpha \}\\ 
  &\!\!\!&\!\!\! + 2\#\{\alpha\in\varDelta^+(\mathfrak{g}_{\mathbb{C}},
\mathfrak{t}_{\mathbb{C}}) \ ; \ \tau(E_{\alpha}) = E_{\alpha} \} 
  + \#\{\alpha\in\varDelta^+(\mathfrak{g}_{\mathbb{C}},
\mathfrak{t}_{\mathbb{C}})\ ; \ \tau(\alpha) = -\alpha \} \\
  &\!\!\!=&\!\!\!\dim\mathfrak{t}_+ + \#\varDelta^+(\mathfrak{g}_{\mathbb{C}},
\mathfrak{t}_{\mathbb{C}}) + 2\#\{\alpha\in\varDelta^+(\mathfrak{g}_{\mathbb{C}},
\mathfrak{t}_{\mathbb{C}}) \ ; \ \tau(E_{\alpha}) = E_{\alpha} \} \\
&\!\!\!&\!\!\!- \#\{\alpha\in\varDelta^+(\mathfrak{g}_{\mathbb{C}},\mathfrak{t}_{\mathbb{C}}) \
  ; \ \tau(\alpha) = \alpha \}.
\end{array}$$\hfill$\Box$


\section{The case where $\dim\mathfrak z = 0$ and $\tau\circ\sigma =
 \sigma^{-1}\circ\tau$} 
 We consider the cases of Type II, III, IV, V and VI in Table 1. 
First we construct $\tau$ by using graded Lie algebras. Let $\mathfrak 
g^*$ be a normal real form of a complex simple Lie algebra 
$\mathfrak{g}_{\mathbb{C}}$. 
Let $\mathfrak{t}^*$ be a Cartan subalgebra of $\mathfrak{g}^*$. Then 
we have a Cartan decomposition $\mathfrak{g}^*=\mathfrak{k}+\mathfrak{p}^*$ 
with 
\begin{equation}
  \mathfrak k := \sum_{\alpha\in\varDelta^+ (\mathfrak{g}^*_{\mathbb{C}},
   \mathfrak{t}^*_{\mathbb{C}}) } \mathbb{R} A_{\alpha}, \ \ \ 
  \ \ \mathfrak{p}^* := \mathfrak{t}^* + \sum_{\alpha\in\varDelta^+
  (\mathfrak{g}^*_{\mathbb{C}}, \mathfrak{t}^*_{\mathbb{C}}) }  
{\mathbb{R}}\sqrt{-1} B_{\alpha}.
   \label{eqn:51}
\end{equation}
We take a gradation $\mathfrak{g}^*
=\sum_{p=-4}^4\mathfrak{g}_p^*$ of the fourth kind on $\mathfrak{g}^*$
corresponding to a partition  
$$\varPi = \varPi_0\cup\varPi_1,\quad\varPi_1=\{\alpha_i\},\quad m_i =4.$$
Then the characteristic element of the gradation coincides with $K_i$. 
\par
Let $\tau^*$ be the Cartan involution defined by (\ref{eqn:25}). 
 Put $\sigma := {\rm Ad}(\exp(\pi/2)\sqrt{-1}K_i)$. Then $\sigma$ is an 
automorphism of order $4$ on the compact dual 
$\mathfrak{g}:=\mathfrak{k}+\sqrt{-1}\mathfrak{p}^*$ of $\mathfrak{g}^*$. 
Since $\tau^*(K_i)=-K_i$, it is obvious that
\begin{equation}
\tau^*\circ\sigma\circ(\tau^*)^{-1}=\sigma^{-1}.
\label{eqn5-2}
\end{equation} 
By Lemma \ref{lem:35} and Proposition \ref{prop:41}, $\tau^*$ is conjugate 
within ${\rm Int}(\mathfrak{h})$ to an involutive automorphism $\tau^{\varPi}$ 
of Type II, III, IV, V or VI in Table 1, that is, there exists
$\mu\in{\rm Int}(\mathfrak{h})$ such that 
$\tau^{\varPi}|_{\mathfrak{t}}=(\mu\circ\tau^*\circ\mu^{-1})|_{\mathfrak{t}}$.  Note that
$\dim \mathfrak{z} = 0$ by Theorem 5.15 of Chapter X of \cite{Helg}, 
and it follows from (\ref{eqn:51}) that 
\begin{equation}
   \mathfrak {h\cap k} = \sum_{\alpha\in\varDelta^+
    (\mathfrak{g}^*_{\mathbb{C}}, \mathfrak{t}^*_{\mathbb{C}})
    \atop\alpha(K_i)\equiv 0 
   (\bmod 4)} \mathbb{R} A_{\alpha}.
     \label{eqn:52}
\end{equation}
Now we prove the following Lemma.
\begin{lem} Let $\mu$ be in ${\rm Int}(\mathfrak{h})$. Then 
   $\mu\circ\tau^*\circ\mu^{-1}$ is conjugate within ${\rm Int}(\mathfrak h)$ to
   $\mu\circ\tau^*\circ\mu^{-1}\circ\sigma$. 
     \label{lem:51}
\end{lem}

 Proof. Put $\nu:=\tau_{(-1/4)K_i}$. 
Then we have $\nu\circ\tau^*\circ\nu^{-1}=\tau^*\circ\nu^{-1}\circ\nu^{-1}
=\tau^*\circ\sigma$. Since $\nu\in{\rm Int}(\mathfrak h)$, it follows that 
$$\mu\circ\tau^*\circ\mu^{-1}\circ\sigma =
\mu\circ\tau^*\circ\sigma\circ\mu^{-1} =
\mu\circ\nu\circ\tau^*\circ\nu^{-1}\circ\mu^{-1} =
\nu\circ\mu\circ\tau^*\circ\mu^{-1}\circ\nu^{-1},$$ 
and hence $\mu\circ\tau^*\circ\mu^{-1}$ is conjugate within 
${\rm Int}(\mathfrak h)$ to
$\mu\circ\tau^*\circ\mu^{-1}\circ\sigma$. \hfill$\Box$\vspace{0.5cm}\par

In the remaining part of this section, we shall determine all involutions for
each type. Furthermore for each involution $\tau$, 
we shall determine $\mathfrak{h\cap g}^{\tau}$ and
$\mathfrak{g}^{\tau}$. \par
Let $\tau^{\varPi}_1, \tau^{\varPi}_2, \tau^{\varPi}_3,
\tau^{\varPi}_4$ be the involutive automorphisms which conjugate 
within ${\rm Int}(\mathfrak{h})$ to the Cartan involutions $\tau^*$ with
respect to Type IV, V, III, VI, respectively. We denote by $\tau$ any
involution of each types.
\vspace{0.3cm}\par

 {\it Type} IV : Now, we investigate involutions of Type IV in Table 1. 
Since $\mathfrak{g}^*$ is a normal real form and of type $\mathfrak{e}_8$, 
the pair $(\mathfrak{g^*,k})$ is given by $(\mathfrak{e}_{8(8)},
\mathfrak{so}(16))$. Note that $\dim\mathfrak k = 120$. 
Set $\sigma={\rm Ad}(\exp(\pi/2)\sqrt{-1}K_3)$. From (\ref{eqn:51}) and
(\ref{eqn:52}), considering the number of roots 
$\alpha\in\varDelta^+(\mathfrak{g}^*_{\mathbb{C}},
\mathfrak{t}^*_{\mathbb{C}})$ such that $\alpha(K_3)=0$ or $4$ (for
example see Freudenthal and Vries~\cite{FV}), we get
$\dim(\mathfrak{h\cap k}) = 29$. 
Then it follows from (\ref{eqn5-2}) and 
Proposition \ref{prop:41} that $\tau^{\varPi}_1$ is of Type IV in Table 1.
\par
Let $\mathfrak{t}_{\pm}$ be the $(\pm 1)$-eigenspaces of 
$\tau^{\varPi}_1|_{\mathfrak{t}}$, respectively. 
Since $\alpha_i(\tau^{\varPi}_1(K_j))=\tau^{\varPi}_1(\alpha_i)(K_j)$, we have
$$\begin{array}{l}
   \mathfrak{t}_+ = {\rm span}\{2K_1 - K_2, \ 2K_1 - K_3, \ 2K_1 + K_6, \
    4K_1 - K_5 - K_7, \ 4K_1 + K_4 + K_8\}, \\
   \mathfrak{t}_- = {\rm span}\{2K_2 - K_3 , \ 2K_2 - K_4 + K_8, \ K_2
    - K_5 +K_7 \}.
\end{array}$$
For $\tau_{h_-} \ (h_-\in\mathfrak{t}_-)$, we have
$\tau^{\varPi}_1\circ\tau_{h_-}\circ\tau^{\varPi}_1 = \tau_{\tau^{\varPi}_1(h_-)}=\tau_{-h_-}$. Thus we
get
\begin{equation}
   (\tau_{h_-})^{-1}\circ\tau^{\varPi}_1\circ\tau_{h_-}=\tau^{\varPi}_1\circ\tau_{2h_-}.
     \label{eqn:53}  
\end{equation}
Then using $h_-:=t(K_3-K_4+K_8)\in{\mathfrak{t}}_-$, we may assume 
$\tau^{\varPi}_1(E_{\alpha_4})=E_{\alpha_8}$. Indeed, if
$\tau^{\varPi}_1(E_{\alpha_4}) =b_4E_{\alpha_8} \ (b_4\in\mathbb{C}, \ |b_4|=1)$, then it follows from
(\ref{eqn:53}) that  
$$(\tau_{h_-})^{-1}\circ\tau^{\varPi}_1\circ\tau_{h_-}(E_{\alpha_4}) =
b_4e^{-2t\pi\sqrt{-1}}E_{\alpha_8},$$
Taking $t$ so that $b_4=e^{2t\pi\sqrt{-1}}$, we may assume
$\tau^{\varPi}_1(E_{\alpha_4})=E_{\alpha_8}$.
Similarly, using $h_-=t(2K_2 - K_3)$ or $h_-=2(K_2 -K_5 +K_7)-(2K_2 -K_3)$, 
we may assume $\tau^{\varPi}_1(E_{\alpha_2})=E_{\alpha_0}$ and 
$\tau^{\varPi}_1(E_{\alpha_5})=E_{\alpha_7}$.\par 
 On the other hand, for any involution $\tau$ of Type IV, 
the number of the subsets $\{\alpha , \beta\}$ such that
$\alpha\in\varDelta^+ (\mathfrak{g}_{\mathbb{C}},
\mathfrak{t}_{\mathbb{C}}) , \ \tau(\alpha) = \beta, \
\alpha\ne\pm\beta$ and $\alpha(K_3)\equiv 0\pmod 4$ is $12$. Since
$\dim\mathfrak{t}_+ = 5$, by an argument similar to the proof of Lemma
\ref{lem:41} we obtain 
$$\dim(\mathfrak{h\cap g}^{\tau})\geq 5 + (12 \times 2) = 29.$$ 
Because $\dim(\mathfrak{h}\cap\mathfrak{g}^{\tau^{\varPi}_1})=
\dim(\mathfrak{h}\cap\mathfrak{k})=29$, we obtain
\begin{equation}
  \begin{array}{l}
    \tau^{\varPi}_1(E_{\alpha_1}) = -E_{\alpha_1}, \ \tau^{\varPi}_1(E_{\alpha_2}) =
    E_{\alpha_0}, \ \tau^{\varPi}_1(E_{\alpha_3}) = c_1E_{\beta_1}, \\
    \tau^{\varPi}_1(E_{\alpha_4}) = E_{\alpha_8}, \ \tau^{\varPi}_1(E_{\alpha_5}) =
    E_{\alpha_7}, \ \tau^{\varPi}_1(E_{\alpha_6}) = -E_{\alpha_6},
  \end{array}
\label{eqn5-5}
\end{equation}
where $\beta_1 =\alpha_1+2\alpha_2+3\alpha_3+4\alpha_4+3\alpha_5+2\alpha_6
+\alpha_7$ (see Table 1) and $c_1\in{\mathbb{C}}$ with $|c_1|=1$ (cf. 
Corollary 5.2 of Chapter IX of \cite{Helg}). 
\vspace{0.5cm}\par

{\sc Remark 5.1.}
Except for conjugations within ${\rm Aut}_{\mathfrak h}(\mathfrak{g})$, 
we can determine the constant $c_1$ uniquely. Indeed, from the proof of 
Theorem 5.1 of Chapter IX of \cite{Helg}, there exists 
$\mu\in{\rm Aut}(\mathfrak{g})$ such that
\begin{equation}
  \begin{array}{l}
    \mu(E_{\alpha_1}) = E_{\alpha_1}, \ \mu(E_{\alpha_2}) =
    E_{\alpha_0}, \ \mu(E_{\alpha_3}) = E_{\beta_1}, \ \mu(E_{\alpha_4}) =
    E_{\alpha_8}, \\
    \mu(E_{\alpha_5}) = E_{\alpha_7}, \ \mu(E_{\alpha_6}) = E_{\alpha_6}, \
    \mu(E_{\alpha_7}) = E_{\alpha_5}, \ \mu(E_{\alpha_8}) = E_{\alpha_4},
    \\
    \mu(E_{\alpha_0}) = \epsilon_0
    E_{\alpha_2}, \ \mu(E_{\beta_1}) = \epsilon_{\beta_1} E_{\alpha_3}, 
  \end{array}
\label{eqn:56}
\end{equation}
where $\epsilon_0 =\pm 1$ and $\epsilon_{\beta_1} = \pm 1$. 
Note that $\epsilon_0$ and $\epsilon_{\beta_1}$ are uniquely determined since 
$E_{\pm\alpha_i}$ ($1\leq i\leq 8$) generate ${\mathfrak{g}}_{\mathbb{C}}$. 
Since $((\tau^{\varPi}_1)^{-1}\circ\mu)|_{\mathfrak{t}}={\rm Id}_{\mathfrak{t}}$,  
it follows from Proposition 5.3 of Chapter IX of \cite{Helg} that 
there exists $\sqrt{-1}H\in\mathfrak{t}$ such that 
$(\tau^{\varPi}_1)^{-1}\circ\mu= \tau_H$, and therefore $\mu=\tau^{\varPi}_1\circ\tau_H$.
Put $H=\sum_{i=1}^8 a_i K_i, \ a_i\in\mathbb{R}$. Then from
(\ref{eqn:20}), we have
$$
E_{\alpha_1}=\mu(E_{\alpha_1}) =
 \tau^{\varPi}_1\circ\tau_H(E_{\alpha_1})=e^{\pi\sqrt{-1}\alpha_1(H)}\tau^{\varPi}_1(E_{\alpha_1})
 =-e^{\pi\sqrt{-1}\alpha_1(H)}E_{\alpha_1}.
$$
Thus we get $a_1=\alpha_1(H)\equiv 1\pmod2$. Similarly as above, we
obtain $a_2\equiv 0, \ a_4\equiv 0, \ a_5\equiv 0, \ a_6\equiv 1, \
a_7\equiv 0, \ a_8\equiv 0\pmod 2$. Moreover, since 
$$
E_{\beta_1}=\mu(E_{\alpha_3}) =
 e^{\pi\sqrt{-1}\alpha_3(H)}\tau^{\varPi}_1(E_{\alpha_3})
 =c_1e^{\pi\sqrt{-1}a_3}E_{\beta_1},
$$
we have 
\begin{equation}
  c_1=e^{-\pi\sqrt{-1}a_3}.
   \label{eqn:618}
\end{equation}
Then by (\ref{eqn5-5}) and (\ref{eqn:56}) we have 
$$\epsilon_0E_{\alpha_2}=\mu(E_{\alpha_0})=\tau^{\varPi}_1\circ\tau_H(E_{\alpha_0})
=e^{-4\pi\sqrt{-1}a_3}E_{\alpha_2},
$$
and it follows from (\ref{eqn:618}) that $c_1^4=\epsilon_0$. 
\par
If $\epsilon_0=1$, then $c_1=\pm 1$ or $\pm\sqrt{-1}$. Considering  
(\ref{eqn:53}) for $h_-=2K_2-K_3\in{\mathfrak{t}}_{-}$, we may assume that 
$c_1=1$ or $\sqrt{-1}$. Moreover, by Lemma \ref{lem:51} we may assume that 
$c_1=1$. If $\epsilon_0=-1$, then by the same argument as above, we may assume 
$c_1=e^{(\pi/4)\sqrt{-1}}$. 
Consequently, $c_1$ is uniquely determined except for 
conjugations within ${\rm Aut}_{\mathfrak h}(\mathfrak{g})$. 
\vspace{0.5cm}\par

By an argument similar to (\ref{eqn5-5}), we may assume
\begin{equation}
  \begin{array}{l}
    \tau(E_{\alpha_1}) = \pm E_{\alpha_1}, \ \tau(E_{\alpha_2}) =
    E_{\alpha_0}, \ \tau(E_{\alpha_3}) = \tilde{c}_1E_{\beta_1}, \\
    \tau(E_{\alpha_4}) = E_{\alpha_8}, \ \tau(E_{\alpha_5}) =
    E_{\alpha_7}, \ \tau(E_{\alpha_6}) = \pm E_{\alpha_6},
  \end{array}
\label{eqn5-8}
\end{equation}
where $\tilde{c_1}\in\mathbb{C}$ and $|\tilde{c}|=1$. Then, by
Proposiion 5.3 of Chapter IX of \cite{Helg}, there exists
$\sqrt{-1}h\in\mathfrak{t}$ such that
$\tau =\tau^{\varPi}_1\circ\tau_{h}$. Put
$$
\begin{array}{ll}
h  := &\!\!\! h_+ + h_-, \\
h_+ := &\!\!\! k_1(2K_1 - K_2) + k_2(2K_1 - K_3) + k_3(2K_1 + K_6) \\
       &\!\!\! + k_4(4K_1 - K_5 - K_7) + k_5(4K_1 + K_4 + K_8) \in
	\sqrt{-1}\mathfrak{t}_+, \\ 
h_- :=  &\!\!\! k_6(2K_2 - K_3) +k_7(2K_2 - K_4 +
        K_8) +k_8(K_2 - K_5 +K_7)\in\sqrt{-1}\mathfrak{t}_-,\\
\end{array}$$
where $k_1,\ldots ,k_8\in\mathbb{R}.$
Then since $\tau^2={\rm Id}$ and $\tau^{\varPi}_1(h)=h_+-h_-$, we have
$\tau_{2h_+}=\rm{Id}$ and hence $2h_+\equiv
0\pmod{2\varPi({\mathfrak{g}}_{\mathbb{C}},
{\mathfrak{t}}_{\mathbb{C}})}$. Therefore we get $k_1,\ldots , k_5
\in\mathbb Z$. Then we have
$$
\begin{array}{ll}
   h \equiv & \!\!\! k_1 K_2 + k_2 K_3 +k_5 K_4 +k_4 K_5 +k_3 K_6 +k_4 K_7 +
    k_5 K_8 \\
  &\!\!\! +k_6(2K_2 - K_3) +k_7(2K_2 - K_4 + K_8) +k_8(K_2 - K_5 +K_7)
\pmod{2\varPi({\mathfrak{g}}_{\mathbb{C}},
{\mathfrak{t}}_{\mathbb{C}})}.
\end{array}
$$
Considering (\ref{eqn5-5}) and (\ref{eqn5-8}) together with 
(\ref{eqn:20}), we obtain 
$$\alpha_2(h)\equiv \alpha_4(h)\equiv\alpha_8(h)\equiv\alpha_5(h)\equiv
\alpha_7(h)\equiv 0 \pmod{2},$$ 
and therefore 
$$h\equiv (k_2-k_6)K_3+k_3K_6 \pmod{2\varPi({\mathfrak{g}}_{\mathbb{C}},
{\mathfrak{t}}_{\mathbb{C}})}.$$
Furthermore, since $\tau(E_{\alpha_0})=\tau^{\varPi}_1(E_{\alpha_0})
=E_{\alpha_2}$, it follows that $\alpha_0(h)\equiv0 \pmod{2}$, and therefore 
$2k_6\in\mathbb{Z}$. Hence we may assume that $\tau$ is one of the following:
$$\tau^{\varPi}_1, \ \ \tau^{\varPi}_1\circ\tau_{K_j},\ \ 
\tau^{\varPi}_1\circ\tau_{K_j}\circ\sigma,\ \ 
\tau^{\varPi}_1\circ\tau_{K_3+K_6},\ \ \tau^{\varPi}_1\circ
\tau_{K_3+K_6}\circ\sigma, \ \ j=3,6.$$
Indeed, $\tau^{\varPi}_1\circ\tau_{-k_6K_3}$ is conjugate within 
${\rm Int}(\mathfrak{t})$ to one of $\tau^{\varPi}_1$ and 
$\tau^{\varPi}_1\circ\sigma$ since 
$$\tau^{\varPi}_1\circ\tau_{-(1/2)K_3}=\tau^{\varPi}_1\circ\sigma^{-1}
=\sigma\circ(\tau^{\varPi}_1\circ\sigma)\circ\sigma^{-1}.$$
Moreover, since $\tau_{K_3}=\sigma^2$ and $\tau^{\varPi}_1\circ\sigma=
\sigma^{-1}\circ\tau^{\varPi}_1$, it follows that 
$\tau^{\varPi}_1\circ\tau_{K_3}$ and $\tau^{\varPi}_1\circ\tau_{K_3+K_6}$ are 
conjugate within ${\rm Int}(\mathfrak{t})$ to $\tau^{\varPi}_1$ and 
$\tau^{\varPi}_1\circ\tau_{K_6}$, respectively. Consequently, $\tau$ is 
conjugate within ${\rm Aut}_{\mathfrak{h}}(\mathfrak{g})$ to one of following: 
$$
  \tau^{\varPi}_1, \qquad  \ \  \tau^{\varPi}_1\circ\tau_{K_6},
  \qquad  \ \ \tau^{\varPi}_1\circ\tau_{K_6}\circ\sigma. 
$$

Now we shall compute the dimension of $\mathfrak{h\cap g}^{\theta}$ and
$\mathfrak{g}^{\theta}$, where $\theta$ is one of $\tau^{\varPi}_1$, 
$\tau^{\varPi}_1\circ\tau_{K_6}$ and $\tau^{\varPi}_1\circ\tau_{K_6}\circ\sigma$. 
Since $\tau^{\varPi}_1\circ\tau_{K_6}(E_{\alpha_6})=E_{\alpha_6}$ and
$\dim(\mathfrak{h}\cap\mathfrak{g}^{\tau^{\varPi}_1})=29$, we have 
$\dim(\mathfrak{h}\cap\mathfrak{g}^{\tau^{\varPi}_1
\circ\tau_{K_6}})=36$. Therefore we get
$\mathfrak{h}\cap\mathfrak{g}^{\tau^{\varPi}_1}\cong D_4\oplus D_1$ and
$\mathfrak{h}\cap\mathfrak{g}^{\tau^{\varPi}_1 \circ\tau_{K_6}}\cong C_4 \oplus
D_1$. Put $\nu:=\tau|_{\mathfrak{t}}$. 
It is easy to see that positive roots $\alpha$ such that
$\nu(\alpha)=\alpha$ are
$$\varDelta_{\nu}^+ :=\left\{
\begin{tabular}{p{300pt}}
 $\alpha_1, \ \alpha_6,\ \alpha_5+ \alpha_6+ \alpha_7,\ \alpha_4+
  \alpha_5+ \alpha_6 + \alpha_7 + \alpha_8, \ \alpha_1 + \alpha_2
 +2\alpha_3 + 2\alpha_4 + 
 \alpha_5, \ \alpha_1+ \alpha_2+ 2\alpha_3+ 2\alpha_4+ \alpha_5+
 \alpha_6, \ \alpha_1+ \alpha_2+ 2\alpha_3+ 2\alpha_4+ 2\alpha_5+
 \alpha_6+ \alpha_7, \ \alpha_1+ \alpha_2+ 2\alpha_3+ 2\alpha_4+
 2\alpha_5+ 2\alpha_6+ \alpha_7, \ 
  \alpha_1 + \alpha_2 + 2\alpha_3 + 3\alpha_4 + 2\alpha_5 + \alpha_6 +
  \alpha_7 +\alpha_8, \
  \alpha_1 + \alpha_2 + 2\alpha_3 + 3\alpha_4 + 2\alpha_5 + 2\alpha_6 +
  \alpha_7 +\alpha_8, \
  \alpha_1 + \alpha_2 + 2\alpha_3 + 3\alpha_4 + 3\alpha_5 + 2\alpha_6 +
  2\alpha_7 +\alpha_8, \
  \alpha_1 + \alpha_2 + 2\alpha_3 + 3\alpha_4 + 3\alpha_5 + 3\alpha_6 +
  2\alpha_7 +\alpha_8, \ 2\alpha_1+ 2\alpha_2+ 4\alpha_3+ 5\alpha_4+
 4\alpha_5+ 3\alpha_6+ 2\alpha_7+ \alpha_8$
\end{tabular}
\right\}.$$\vspace{0.2cm}\par
We consider the case of $\tau^{\varPi}_1\circ\tau_{K_6}$. 
Put $\gamma:=\alpha_1 + \alpha_2 +2\alpha_3 + 2\alpha_4 + \alpha_5.$
Take a Weyl basis so that $\tau^{\varPi}_1(E_{\gamma})=E_{\gamma}$ 
(cf. see Gilkey and Seitz \cite{GS}). Then it is easy to see that 
$\tau^{\varPi}_1(E_{\alpha}) = E_{\alpha}$ for any
$\alpha\in\varDelta_{\nu}^+\setminus\{\alpha_1, \ \alpha_6,\ \alpha_5+ \alpha_6+
\alpha_7,\ \alpha_4+ \alpha_5+ \alpha_7, \ 2\alpha_1+ 2\alpha_2+
4\alpha_3+ 5\alpha_4+ 4\alpha_5+ 3\alpha_6+ 2\alpha_7+ \alpha_8\}$ and
therefore $\tau^{\varPi}_1\circ\tau_{K_6}(E_{\alpha}) = E_{\alpha}$ for any
$\alpha\in\varDelta_{\nu}^+\setminus\{\alpha_1\}$. It follows from Lemma
4.1 that $\dim\mathfrak{g}^{\tau^{\varPi}_1\circ\tau_{K_6}}= 136$. By using the
classification of symmetric spaces, we get
$\mathfrak{g}^{\tau^{\varPi}_1\circ\tau_{K_6}}\cong E_7\oplus A_1$.\par

Similarly as above we can obtain $\mathfrak{h\cap g}^{\theta}$ and
$\mathfrak{g}^{\theta}$ for $\theta=\tau^{\varPi}_1\circ\tau_{K_6}
\circ\sigma$.
\par
By an argument similar to above, we can obtain all involutions $\tau$
of Type V and VI, and determine $\mathfrak{h\cap g}^{\tau}$ and 
$\mathfrak{g}^{\tau}$, which are listed in Table 2.
\vspace{0.5cm}\par

Now we investigate involutions of Type II and III in Table 1. Since
 $\mathfrak{g}^*$ is a normal real form and of type $\mathfrak{e}_7$,
 the pair $(\mathfrak{g}^*,\mathfrak{k})$ is given by
 $(\mathfrak{e}_{7(7)}, \mathfrak{su}(8))$. It is easy to see that
 $\dim(\mathfrak{h\cap k})=13$. 
On the other hand, for an involution $\tau$ of $\mathfrak{g}$, 
we can see that 
$$\dim{\mathfrak{t}}_+=\left\{\begin{array}{ll}
4 & \mbox{if $\tau$ is of Type II}, \\
5 & \mbox{if $\tau$ is of Type III}, 
\end{array}\right.
$$
and if $\tau$ is of Type II (resp. Type III), the number of the subsets
$\{\alpha , \beta\}$ such that $\alpha\in\varDelta^+
(\mathfrak{g}_{\mathbb{C}},\mathfrak{t}_{\mathbb{C}}) , \ \tau(\alpha)
= \beta, \ \alpha\ne\pm\beta$ and $\alpha(K_4)\equiv 0\pmod 4$ is $6$
(resp. $4$). 
Hence we obtain 
$$\left\{
\begin{array}{ll}
\dim(\mathfrak{h\cap g}^{\tau})\geq 16 & \mbox{if $\tau$ is of Type II} \\
\dim(\mathfrak{h\cap g}^{\tau})\geq 13 & \mbox{if $\tau$ is of Type III}.
\end{array}
\right.
$$
Therefore the Cartan involution $\tau^*$ of $\mathfrak{g}^*=\mathfrak{e}_{7(7)}$ 
is conjugate within ${\rm Int}(\mathfrak h)$ to an involution $\tau^{\varPi}_3$ of 
Type III. By an argument similar to Type IV, we can obtain all involutions 
$\tau$ of Type III, which are listed in Table 2.\par

Finally we consider the Type II. Put $\beta_3 :=\alpha_1 + \alpha_2 +
2\alpha_3 + 3\alpha_4 + 2\alpha_5 + \alpha_6$. Let $\tau^{\varPi}_3$ be as above. 
Then, since $\dim(\mathfrak{h\cap g}^{\tau^{\varPi}_3})=\dim(\mathfrak{h\cap k})=13$, 
it follows that 
\begin{equation}
  \begin{array}{lll}
    \tau^{\varPi}_3(E_{\alpha_1}) = -E_{\alpha_1}, & \tau^{\varPi}_3(E_{\alpha_2}) =
     -E_{\alpha_2}, & \tau^{\varPi}_3(E_{\alpha_3}) = E_{\alpha_0}, \\
     \tau^{\varPi}_3(E_{\alpha_4}) = c_3E_{\beta_3}, & \tau^{\varPi}_3(E_{\alpha_5}) =
      E_{\alpha_7}, & \tau^{\varPi}_3(E_{\alpha_6}) = -E_{\alpha_6},
  \end{array}
  \label{eqn:58}
\end{equation}
for some $c_3\in\mathbb{C}$ with $|c_3|=1$. On the other hand, from
Theorem 5.1 of Chapter IX of 
\cite{Helg}, There exists an automorphism $\varphi$ on $\mathfrak
 g$ such that 
$$
\begin{array}{llll}
\varphi(E_{\alpha_1}) = E_{\alpha_6}, & \varphi(E_{\alpha_2}) = E_{\alpha_2},
 & \varphi(E_{\alpha_3}) = E_{\alpha_5}, & \varphi(E_{\alpha_4}) =
 E_{\alpha_4}, \\
\varphi(E_{\alpha_5}) = E_{\alpha_3}, & \varphi(E_{\alpha_6}) = E_{\alpha_1},
 & \varphi(E_{\alpha_7}) = E_{\alpha_0}, & \varphi(E_{\alpha_0})
=\epsilon E_{\alpha_7}, \\
\end{array}$$\\
where $\epsilon = \pm 1$.
If $\epsilon = -1$, then we have 
$$\left\{\begin{array}{ll}
\varphi^2(E_{\alpha_i}) = E_{\alpha_i} & (1\leq i \leq 6), \\
\varphi^2(E_{\alpha_7}) = -E_{\alpha_7} .& \\
\end{array}\right.$$
Thus the inner automorphism $\varphi^2$ has the form
$\tau_{K_7}$. Hence we have
\begin{equation}
   \mathfrak{g}^{\varphi^2} = \mathfrak{t} +
   \sum_{\alpha\in\varDelta^+(\mathfrak{g}_{\mathbb{C}},
   \mathfrak{t}_{\mathbb{C}}) \atop\alpha(K_7)=0} (\mathbb{R} A_{\alpha} + 
   \mathbb{R} B_{\alpha}).
       \label{eqn:59}
\end{equation}
Put $\gamma := \alpha_2 + \alpha_3 + 2\alpha_4 + 2\alpha_5 + 2\alpha_6 +
\alpha_7$. Then we get $\varphi(\gamma) = -\gamma$ and from the proof of 
Theorem 5.1 of Chapter IX of \cite{Helg}, we get
$$\varphi(E_{\gamma}) = \epsilon_{\gamma} E_{-\gamma}, \ \varphi(E_{-\gamma}) = \epsilon_{-\gamma} E_{\gamma} \ \ (\epsilon_{\gamma}
= \epsilon_{-\gamma} = \pm 1).$$
Therefore we obtain
$$\begin{array}{l}
   \varphi(A_{\gamma}) = \epsilon_{\gamma}(E_{-\gamma}-E_{\gamma}) =
   -\epsilon_{\gamma}A_{\gamma}, \\ 
   \varphi(B_{\gamma}) = \sqrt{-1}\epsilon_{\gamma}(E_{\gamma} +
    E_{-\gamma})=\epsilon_{\gamma}B_{\gamma}. 
\end{array}$$
This implies that
$$A_{\gamma} \ {\rm or} \ B_{\gamma}\in\mathfrak
g^{\varphi}\subset\mathfrak{g}^{\varphi^2}.$$
This contradicts (\ref{eqn:59}). Thus $\varphi(E_{\alpha_0}) =
E_{\alpha_7}$, that is,
\begin{equation}
  \begin{array}{llll}
    \varphi(E_{\alpha_1}) = E_{\alpha_6}, & \varphi(E_{\alpha_2}) =
     E_{\alpha_2}, & \varphi(E_{\alpha_3}) = E_{\alpha_5}, &
     \varphi(E_{\alpha_4}) = E_{\alpha_4}, \\ 
    \varphi(E_{\alpha_5}) = E_{\alpha_3}, & \varphi(E_{\alpha_6}) =
     E_{\alpha_1}, & \varphi(E_{\alpha_7}) = E_{\alpha_0}, &
     \varphi(E_{\alpha_0}) = E_{\alpha_7}. \\ 
  \end{array}
  \label{eqn:510}
\end{equation}
Then $\tau^{\varPi}_3\circ\varphi$ maps
$$\alpha_1\mapsto \alpha_6, \ \alpha_2\mapsto \alpha_2, \ \alpha_3\mapsto
 \alpha_7, \ \alpha_5\mapsto \alpha_0, \ \alpha_4\mapsto\beta_3 ,$$
and it is easy to see
\begin{equation}
  \begin{array}{l}
    \tau^{\varPi}_3\circ\varphi(E_{\alpha_1}) = \varphi\circ\tau^{\varPi}_3(E_{\alpha_1})
     = -E_{\alpha_6}, \ \tau^{\varPi}_3\circ\varphi(E_{\alpha_2}) =
     \varphi\circ\tau^{\varPi}_3(E_{\alpha_2}) = -E_{\alpha_2},\\
    \tau^{\varPi}_3\circ\varphi(E_{\alpha_3}) =
    \varphi\circ\tau^{\varPi}_3(E_{\alpha_3}) = E_{\alpha_7}, \
    \tau^{\varPi}_3\circ\varphi(E_{\alpha_5}) =
    \varphi\circ\tau^{\varPi}_3(E_{\alpha_5}) = E_{\alpha_0}, \\
    \tau^{\varPi}_3\circ\varphi(E_{\alpha_4}) = c_3E_{\beta_3}, \
    \varphi\circ\tau^{\varPi}_3(E_{\alpha_4}) = c_3\varphi(E_{\beta_3}).
  \end{array}
      \label{eqn:511}
\end{equation}
Therefore $(\tau^{\varPi}_3\circ\varphi)^2={\rm Id}$ if and only if
$\varphi(E_{\beta_3})=E_{\beta_3}$. Since $\dim\mathfrak{g}_{\beta_3}=1$, we
have $\mathfrak{g}_{\beta_3}=\mathbb{C} X_{\beta_3}$, where
$$
   X_{\beta_3}:=[[[[[[[[[E_{\alpha_2},E_{\alpha_4}],E_{\alpha_3}],
   E_{\alpha_5}], E_{\alpha_4}], E_{\alpha_1}],E_{\alpha_6}],
   E_{\alpha_3}], E_{\alpha_5}],E_{\alpha_4}].
$$
Because
$[E_{\alpha_3},E_{\alpha_5}]=[E_{\alpha_1},E_{\alpha_6}]=0$, we get
$\varphi(X_{\beta_3})=X_{\beta_3}$ and therefore 
$\varphi|_{\mathfrak{g} _{\beta_3}}={\rm Id}$. Thus $\tau^{\varPi}_3\circ\varphi$ is
an involutive automorphism of $\mathfrak{e}_7$. From (\ref{eqn:511}), we obtain
\begin{equation}
 \begin{array}{l}
    \tau^{\varPi}_3\circ\varphi(E_{\alpha_1}) = -E_{\alpha_6},\
     \tau^{\varPi}_3\circ\varphi(E_{\alpha_2}) = -E_{\alpha_2},\ 
    \tau^{\varPi}_3\circ\varphi(E_{\alpha_3}) = E_{\alpha_7},\\
    \tau^{\varPi}_3\circ\varphi(E_{\alpha_5}) = E_{\alpha_0}, 
    \ \tau^{\varPi}_3\circ\varphi(E_{\alpha_4}) = c_3E_{\beta_3}.
 \end{array} 
      \label{eqn:512}
\end{equation}
Hence we can construct an involution of Type II. By an argument similar
to Type IV, we can give all involutions of Type II. \par
 Consequently we obtain the following proposition.
\begin{prop} Suppose that $\dim\mathfrak{z}=0$. Let $\tau$ be an involution of 
$\mathfrak{g}$ such that $\tau\circ\sigma=\sigma^{-1}\circ\tau$. Then $\tau$ 
is conjugate within ${\rm Aut}_{\mathfrak h}(\mathfrak{g})$ to one of 
automorphisms listed in Table 2.
\label{prop:51}
\end{prop}

{\small
\begin{center}
\renewcommand{\arraystretch}{1.1}
\begin{longtable}{l|l|l|l}
\caption{$\dim\mathfrak{z}=0, \ \tau\circ\sigma=\sigma^{-1}\circ\tau$, 
$\sigma=\tau_{(1/2)H}$, $\mathfrak{k}=\mathfrak{g}^{\tau}$.}\\
\hline
\multicolumn{1}{c|}{$(\mathfrak{g}, \mathfrak h, H)$} &
 \multicolumn{1}{|c|}{$\tau$}
 &\multicolumn{1}{|c|}{$\mathfrak k$} & 
 \multicolumn{1}{|c}{$\mathfrak {h\cap k}$}\\
\hline
\hline

$(\mathfrak{e}_8,\mathfrak{su}(8)\oplus\mathfrak{su}(2), K_3)$ & $\tau^{\varPi}_1$
& $D_8$ & $D_4\oplus D_1$ \\

& $\tau^{\varPi}_1\circ\tau_{K_6}$ & $E_7\oplus A_1$ & $C_4\oplus D_1$  \\

& $\tau^{\varPi}_1\circ\tau_{K_6}\circ\tau_{(1/2)K_3}$ & $D_8$ & $C_4\oplus D_1$  \\

\hline

$(\mathfrak{e}_8,\mathfrak{so}(10)\oplus\mathfrak{so}(6), K_6)$ &
 $\tau^{\varPi}_2$ & $D_8$ & $B_2\oplus B_2\oplus B_1\oplus B_1$ \\

& $\tau^{\varPi}_2\circ\tau_{K_1+K_4}$ & $D_8$ & $B_2\oplus B_2\oplus B_1\oplus B_1$ \\

& $\tau^{\varPi}_2\circ\tau_{K_1+K_4}\circ\tau_{(1/2)K_6}$ & $D_8$ & $B_2\oplus B_2\oplus
 B_1\oplus B_1$ \\

& $\tau^{\varPi}_2\circ\tau_{K_1+K_8}$ & $E_7\oplus A_1$ & $B_3\oplus B_2\oplus B_1$ \\

& $\tau^{\varPi}_2\circ\tau_{K_1+K_8}\circ\tau_{(1/2)K_6}$ & $D_8$ & $B_3\oplus B_2\oplus
 B_1$ \\

& $\tau^{\varPi}_2\circ\tau_{K_4+K_8}$ & $E_7\oplus A_1$ & $B_3\oplus B_2\oplus B_1$  \\

& $\tau^{\varPi}_2\circ\tau_{K_4+K_8}\circ\tau_{(1/2)K_6}$ & $D_8$ & $B_3\oplus B_2\oplus
 B_1$ \\

& $\tau^{\varPi}_2\circ\tau_{K_3}$ & $D_8$ & $B_2\oplus B_2\oplus B_1\oplus B_1$ \\

& $\tau^{\varPi}_2\circ\tau_{K_3}\circ\tau_{(1/2)K_6}$ & $D_8$ & $B_2\oplus B_2\oplus
 B_1\oplus B_1$ \\ 

& $\tau^{\varPi}_2\circ\tau_{K_1+K_3+K_4}$ & $E_7\oplus A_1$ & $D_2\oplus
 B_1\oplus B_1$ \\

& $\tau^{\varPi}_2\circ\tau_{K_1+K_3+K_4}\circ\tau_{(1/2)K_6}$ & $E_7\oplus A_1$ &
 $D_2\oplus B_1\oplus B_1$ \\ 

& $\tau^{\varPi}_2\circ\tau_{K_1+K_3+K_8}$ & $D_8$ & $B_3\oplus B_2\oplus B_1$  \\

& $\tau^{\varPi}_2\circ\tau_{K_1+K_3+K_8}\circ\tau_{(1/2)K_6}$ & $E_7\oplus A_1$ &
 $B_3\oplus B_2\oplus B_1$ \\ 

& $\tau^{\varPi}_2\circ\tau_{K_3+K_4+K_8}$ & $E_7\oplus A_1$ & $B_3\oplus
 B_2\oplus B_1$  \\ 

& $\tau^{\varPi}_2\circ\tau_{K_3+K_4+K_8}\circ\tau_{(1/2)K_6}$ & $D_8$ & $B_3\oplus
 B_2\oplus B_1$  \\

\hline

$(\mathfrak{e}_7,\mathfrak{so}(6)\oplus\mathfrak{so}(6)\oplus\mathfrak{su}(2),
 K_4)$
 & $\tau^{\varPi}_3$ & $A_7$ & $B_1\oplus B_1\oplus B_1\oplus B_1\oplus
 B_1\oplus\mathbb{R}$ \\ 

& $\tau^{\varPi}_3\circ\tau_{K_1+K_2}$ & $D_6\oplus A_1$ & $B_2\oplus B_1\oplus
 B_1\oplus A_1$  \\ 

& $\tau^{\varPi}_3\circ\tau_{K_1+K_2}\circ\tau_{(1/2)K_4}$ & $D_6\oplus A_1$ & $B_2\oplus
 B_1\oplus B_1\oplus A_1$ \\ 

& $\tau^{\varPi}_3\circ\tau_{K_1+K_6}$ & $E_6\oplus \mathbb{R}$ & $B_2\oplus
 B_2\oplus \mathbb{R}$  \\

& $\tau^{\varPi}_3\circ\tau_{K_1+K_6}\circ\tau_{(1/2)K_4}$ & $A_7$ & $B_2\oplus
 B_2\oplus \mathbb{R}$ \\

& $\tau^{\varPi}_3\circ\tau_{K_2+K_6}$ & $D_6\oplus A_1$ & $B_2\oplus B_1\oplus
 B_1\oplus A_1$  \\ 

& $\tau^{\varPi}_3\circ\tau_{K_2+K_6}\circ\tau_{(1/2)K_4}$ & $D_6\oplus A_1$ & $B_2\oplus
 B_1\oplus B_1\oplus A_1$ \\ 

\cline{2-4}

& $\tau^{\varPi}_3\circ\varphi$ & $D_6\oplus A_1$ & $D_3\oplus D_1$ \\ 

& $\tau^{\varPi}_3\circ\varphi\circ\tau_{(1/2)K_4}$ & $A_7$ & $D_3\oplus D_1$  \\ 

\hline

$(\mathfrak{f}_4,\mathfrak{so}(6)\oplus\mathfrak{so}(3), K_3)$ & $\tau^{\varPi}_4$
& $C_3\oplus A_1$ & $B_1\oplus B_1\oplus D_1$ \\

& $\tau^{\varPi}_4\circ\tau_{K_1+K_4}$ & $B_4$ & $B_2\oplus B_1$ \\

& $\tau^{\varPi}_4\circ\tau_{K_1+K_4}\circ\tau_{(1/2)K_3}$ & $C_3\oplus A_1$ &
 $B_2\oplus B_1$ \\ 

\hline

\multicolumn{4}{l}{$\tau^{\varPi}_1\ :\ E_{\alpha_1}\mapsto -E_{\alpha_1},\ 
 E_{\alpha_2}\mapsto E_{\alpha_0},\  E_{\alpha_3}\mapsto c_1 E_{\beta_1},\ 
 E_{\alpha_4}\mapsto E_{\alpha_8},\  E_{\alpha_5}\mapsto E_{\alpha_7},\ 
 E_{\alpha_6}\mapsto -E_{\alpha_6},$} \\
\multicolumn{4}{l}{\hspace{25pt}($\beta_1=\alpha_1 + 2\alpha_2 + 3\alpha_3
 + 4\alpha_4 + 3\alpha_5 + 2\alpha_6 + \alpha_7$)}\\
\multicolumn{4}{l}{$\tau^{\varPi}_2\ : \ E_{\alpha_1}\mapsto
 -E_{\alpha_1}, \ E_{\alpha_2}\mapsto E_{\alpha_5}, \
  E_{\alpha_3}\mapsto -E_{\alpha_3}, \ 
 E_{\alpha_4}\mapsto -E_{\alpha_4}, \ E_{\alpha_6}\mapsto c_2
 E_{\beta_2}, \ E_{\alpha_7}\mapsto E_{\alpha_0},$} \\ 
\multicolumn{4}{l}{\hspace{25pt}$E_{\alpha_8}\mapsto
 -E_{\alpha_8},$ ($\beta_2=\alpha_1 + \alpha_2 + 2\alpha_3 + 3\alpha_4 +
 3\alpha_5 + 3\alpha_6 + 2\alpha_7 +\alpha_8$)} \\
\multicolumn{4}{l}{$\tau^{\varPi}_3\ :\ E_{\alpha_1}\mapsto -E_{\alpha_1}, \
 E_{\alpha_2}\mapsto -E_{\alpha_2}, \ E_{\alpha_3}\mapsto E_{\alpha_0}, \
 E_{\alpha_4}\mapsto c_3 E_{\beta_3}, \ E_{\alpha_5}\mapsto E_{\alpha_7}, \
 E_{\alpha_6}\mapsto -E_{\alpha_6},$}\\
\multicolumn{4}{l}{\hspace{25pt}($\beta_3=\alpha_1 + \alpha_2 +
 2\alpha_3 + 3\alpha_4 + 2\alpha_5 + \alpha_6$)} \\
\multicolumn{4}{l}{$\tau^{\varPi}_4\ :\ E_{\alpha_1}\mapsto -E_{\alpha_1},\ 
 E_{\alpha_2}\mapsto E_{\alpha_0},\  E_{\alpha_3}\mapsto c_4 E_{\beta_4},\ 
 E_{\alpha_4}\mapsto -E_{\alpha_4},$ ($\beta_4=\alpha_1 + 2\alpha_2 +
 3\alpha_3 + \alpha_4$)}\\
\multicolumn{4}{l}{where $c_i(i=1,2,3,4)$ is some complex number with
 $|c_i|=1$.} 

\end{longtable}
\end{center}
}


\section{The case where $\dim\mathfrak z =1, \tau\circ\sigma =
 \sigma^{-1}\circ\tau$} 
 In this section we investigate involutions $\tau$ of $\mathfrak{g}$ such 
that $\dim\mathfrak{z}=1$ and $\tau\circ\sigma=\sigma^{-1}\circ\tau$. 
 First, we construct such involutions by using graded Lie algebras. 
Let $\mathfrak{g}^*$ be a noncompact simple Lie algebra over $\mathbb{R}$ 
such that $\mathfrak{g}_{\mathbb{C}}$ is simple. Let 
$\mathfrak{g}^*=\mathfrak{k}+\mathfrak{p}^*$ and $\tau$ be the Cartan
involution as in (\ref{eqn:25}), and $\mathfrak{a}$ be a maximal abelian
subspace of 
$\mathfrak{p}^*$. Let $\mathfrak{t}^*$ be a Cartan subalgebra of 
$\mathfrak{g}^*$ such that $\mathfrak{a}\subset\mathfrak{t}^*$. We take
compatible orderings on $\mathfrak{a}$ and $\mathfrak{t}^*$. \par
Take a gradation $\mathfrak{g}^*
=\sum_{p=-3}^3\mathfrak{g}_p^*$ of the third kind on $\mathfrak
g^*$ corresponding to a partition 
$$\varPi = \varPi_0\cup\varPi_1,\quad\varPi_1=\{\lambda_i\},\quad n_i =3.$$
Put $\sigma := {\rm Ad}(\exp(\pi/2)\sqrt{-1}h_i)$. It is obvious that 
$\sigma\in{\rm Int}(\mathfrak{g})$ ($\mathfrak{g}=\mathfrak{k}
+\sqrt{-1}\mathfrak{p}^*$), $\sigma^4={\rm Id}$ and $\tau\circ\sigma=
\sigma^{-1}\circ\tau$. Considering the classification of the Satake diagram, 
there exists a unique $\alpha_j\in\varDelta(\mathfrak{g}^*_{\mathbb{C}}, 
\mathfrak{t}^*_{\mathbb{C}})$ such that $\alpha_j|_{\mathfrak{a}}=\lambda_i$ 
with $m_j=3$, and it follows from Lemma \ref{lem:25} that $h_i=K_j$. 
Therefore, by Remark 2.2 we have $\dim\mathfrak{z}=1$. 
\par 
Generally, let $(G/H,\langle,\rangle,\sigma)$ be a compact Riemannian 
$4$-symmetric space such that $G$ is simple and $\sigma$ is inner.  
As before, let $\mathfrak{t}$ be a maximal abelian subalgebra of 
$\mathfrak{g}$ contained in $\mathfrak{h}$. 
We suppose that $\sigma=\tau_{(1/2)K_i}$ for some $\alpha_i\in
\varPi(\mathfrak{g}_{\mathbb{C}},\mathfrak{t}_{\mathbb{C}})$ with $m_i=3$. 
Then $\mathfrak{z}=\mathbb{R}\sqrt{-1}K_i$. From \cite{J}, a pair
 $(\mathfrak{g},\mathfrak{h})$ is one of the following: 
\begin{eqnarray*}
& & (\mathfrak{e}_6,\mathfrak{su}(3)\oplus\mathfrak{su}(3)\oplus
\mathfrak{su}(2)\oplus\mathbb{R}),\ \ (\mathfrak{e}_7,\mathfrak{su}(2)
\oplus\mathfrak{su}(6)\oplus\mathbb{R}), \ \ 
(\mathfrak{e}_7,\mathfrak{su}(5)\oplus\mathfrak{su}(3)\oplus\mathbb{R}), \\
& & (\mathfrak{e}_8,\mathfrak{su}(8)\oplus\mathbb{R}),\ \ 
(\mathfrak{e}_8,\mathfrak{su}(2)\oplus\mathfrak{e}_6\oplus\mathbb{R}),\ \ 
(\mathfrak{f}_4,\mathfrak{su}(2)\oplus\mathfrak{su}(3)\oplus\mathbb{R}), \ \
(\mathfrak{g}_2,\mathfrak{su}(2)\oplus\mathbb{R}).
\end{eqnarray*}

{\sc Remark 6.1.} 
Each $4$-symmetric pair described in the above is neither symmetric nor 
$3$-symmetric. Indeed, except for $(\mathfrak{g}_2,
\mathfrak{su}(2)\oplus\mathbb{R})$, it follows from the classifications of 
compact $k$-symmetric spaces ($k=2,3$) that each $4$-symmetric pair described 
in the above is not $k$-symmetric ($k=2,3$). 
Now, for $(\mathfrak{g}_2,\mathfrak{su}(2)\oplus\mathbb{R})$, we prove that 
it is not isomorphic to a $k$-symmetric pair ($k=2,3$). First, we note 
that $\sigma=\tau_{(1/2)K_1}$ with $m_1=3$. From the classification of 
compact symmetric spaces, it is obvious that the pair 
$(\mathfrak{g}_2,\mathfrak{su}(2)\oplus\mathbb{R})$ is not symmetric. 
Let $(\mathfrak{g}_2,\theta)$ be a $3$-symmetric pair. Then $\theta$ is 
conjugate to $\tau_{(2/3)K_2}$ and 
$${\mathfrak{g}_2}^{\tau_{(2/3)K_2}}=\mathfrak{su}_{\alpha_1}(2)\oplus
\mathbb{R}\sqrt{-1}K_2\cong\mathfrak{su}(2)\oplus\mathbb{R}.$$
If there exists $\mu\in{\rm Aut}(\mathfrak{g})$ such that 
$\mu({\mathfrak{g}_2}^{\sigma})=\mu({\mathfrak{g}_2}^{\tau_{(1/2)K_1}})
={\mathfrak{g}_2}^{\tau_{(2/3)K_2}}$, 
then we have $\mu(\mathfrak{su}_{\alpha_2}(2))=\mathfrak{su}_{\alpha_1}(2)$. 
Therefore it follows that there exists $k\in\mathbb{C}$ with $|c|=1$ 
such that 
$$\mu(E_{\alpha_2})=cE_{\pm\alpha_1},\ \ \mu(E_{-\alpha_2})
=c^{-1}E_{\mp\alpha_1},$$
which implies that $\mu(H_{\alpha_2})=\pm H_{\alpha_1}$. 
However, this is a contradiction because $|\alpha_1|\not=|\alpha_2|$. 
Consequently, the $4$-symmetric pair $(\mathfrak{g}_2,\sigma)$ is not 
$k$-symmetric ($k=2,3$). 

\medskip\par

Now we assume that $\sigma=\tau_{(1/2)K_i}$ for some
$\alpha_i\in\varPi(\mathfrak{g}_{\mathbb{C}},\mathfrak{t}_{\mathbb{C}})$
with $m_i=3$. Let $\tau$ be an involution 
of $\mathfrak{g}$ such that $\tau\circ\sigma=\sigma^{-1}\circ\tau$. 
Then it is easy to see that $\tau(\mathfrak h)=\mathfrak h$ and
$\tau(\mathfrak z)=\mathfrak z$. Thus we have
$\tau(\sqrt{-1}K_i)=\pm\sqrt{-1}K_i$. If $\tau(\sqrt{-1}K_i)=\sqrt{-1}K_i$, 
then $\tau\circ\sigma=\sigma\circ\tau$. Hence we get
$\tau(\sqrt{-1}K_i)=-\sqrt{-1}K_i$. Let 
$\mathfrak{g}=\mathfrak{k}+\mathfrak{p}$ be the canonical decomposition 
of $\mathfrak{g}$ corresponding to $\tau$. Then we have 
$\sqrt{-1}K_i\in\mathfrak{p}$. Put
$$
\varDelta_{\mathfrak m}^+:=\{\alpha\in\varDelta^+(\mathfrak{g}_{\mathbb{C}},
\mathfrak{t}_{\mathbb{C}}) \ ; \ A_{\alpha},  B_{\alpha}\in\mathfrak m\},\ \ 
   \varDelta_{\mathfrak h}^+:=\{\alpha\in\varDelta^+(\mathfrak{g}_{\mathbb{C}},
\mathfrak{t}_{\mathbb{C}}) \ ; \ A_{\alpha}, B_{\alpha}\in\mathfrak h\}.
$$
\begin{lem}
    $$\varDelta_{\mathfrak h}^+=\{\alpha=\sum_{j=1}^nk_j\alpha_j\in
      \varDelta^+(\mathfrak{g}_{\mathbb{C}},\mathfrak{t}_{\mathbb{C}}) \ ;
     \ k_i=0\}.$$
\end{lem}

 Proof. Since $m_i=3$ and $E_{\alpha}=\sigma(E_{\alpha})
=e^{(\pi\sqrt{-1}/2)\alpha(K_i)}E_{\alpha}$ for any 
$\alpha\in\varDelta_{\mathfrak{h}}^+$, we have $\alpha(K_i)=\sum_{j=1}^n
k_j\alpha_j(K_i)=k_i=0$.\hfill$\Box$
\smallskip\par
We define a subset $\varDelta^+_s \ (s=1,2,3)$ of $\varDelta_{\mathfrak
 m}^+$ as follows:
$$
  \varDelta_s^+:=\{\alpha=\sum_{j=1}^m\in\varDelta^+ \ ; \ k_i=s\}.
$$
Then we have an orthogonal decomposition $\mathfrak m 
=\mathfrak m_1\oplus\mathfrak m_2\oplus\mathfrak m_3$, where 
$$\mathfrak m_s=\sum_{\alpha\in\varDelta_s^+}(\mathbb
RA_{\alpha}+\mathbb{R}B_{\alpha}).$$
\begin{lem}
   $$\tau(\mathfrak m_s)=\mathfrak{m}_s\quad s=1,2,3.$$
\end{lem}

 Proof. Since 
$$
  \sigma(E_{\alpha})=e^{(\pi\sqrt{-1}/2)\alpha(K_i)}E_{\alpha}=\left\{
    \begin{array}{ll}
      \sqrt{-1}E_{\alpha}  & \alpha\in\varDelta_1^+, \\
      -E_{\alpha}          & \alpha\in\varDelta_2^+,\\
      -\sqrt{-1}E_{\alpha} & \alpha\in\varDelta_3^+,
    \end{array}\right.
$$
it follows that 
$$\sigma(X)=-X \Longleftrightarrow X\in\mathfrak m_2.$$ 
Hence if $X\in\mathfrak m_2$, then
$$\sigma(\tau(X))=\tau\circ\sigma^{-1}(X)=-\tau(X).$$
Thus we obtain $\tau(\mathfrak m_2)=\mathfrak m_2$.\par
 Next for $\alpha\in\varDelta_1^+$ (resp. $\varDelta_3^+$), we get
 $\tau(\alpha)\in -\varDelta_1^+$ (resp. $-\varDelta_3^+$). Indeed, since
$$[K_i,\tau(E_{\alpha})]=\tau[\tau(K_i),E_{\alpha}]=-\tau[K_i,E_{\alpha}]
=-\alpha(K_i)\tau(E_{\alpha}),$$ 
and $\tau(E_{\alpha})\in\mathfrak{g}_{\tau(\alpha)}$,
we get $\tau(\alpha)(K_i)=-\alpha(K_i)=-1$ (resp. $-3$).  This completes
the proof of the lemma. \hfill$\Box$\\ 
Put 
$$
\mathfrak{h}^{\pm}:=\{X\in\mathfrak{h} \ ; \ \tau(X)=\pm X\},\ \ 
\mathfrak m^{\pm}_s:=\{X\in\mathfrak m_s  \ ; \ \tau(X)=\pm X\}.
$$
Since $\tau(\mathfrak{h})=\mathfrak{h}$ and $\tau(\mathfrak{m})=\mathfrak{m}$, 
we can write
$$\displaystyle
\begin{array}{l}
   \displaystyle\mathfrak{g} =(\mathfrak{h^++h^-})\oplus\sum_{s=1}^3
    (\mathfrak m_s^++\mathfrak m_s^-),\\
   \displaystyle\mathfrak k=\mathfrak h^+\oplus\mathfrak m_1^+
    \oplus\mathfrak m_2^+\oplus\mathfrak m_3^+,\ \ 
   \displaystyle\mathfrak p=\mathfrak h^-\oplus\mathfrak m_1^-
    \oplus\mathfrak m_2^-\oplus\mathfrak m_3^-.  
\end{array}
$$
Put
$\mathfrak{g^*:=k}+\sqrt{-1}\mathfrak p$. Then we have 
$Z:=K_i\in\sqrt{-1}\mathfrak{p}$. 
We shall prove the following lemma.
\begin{lem}
   The eigenvalues of ${\rm ad}(Z):\mathfrak{g}^*\to\mathfrak{g}^*$ are $0,\pm
   1,\pm 2$ and $\pm 3$. 
     \label{lem:63}
\end{lem}

 Proof. First we note that $\mathfrak{h}^{+}+\sqrt{-1}\mathfrak{h}^-$ is 
the $0$-eigenspace of ${\rm ad}(Z)$. \par
It is easy to see that 
\begin{equation}
  \left\{
   \begin{array}{l}
     \alpha\in\varDelta_1^+\Longrightarrow\sigma(A_{\alpha})=B_{\alpha}
      ,\ \sigma(B_{\alpha})=-A_{\alpha}, \\
     \alpha\in\varDelta_2^+\Longrightarrow\sigma(A_{\alpha})=-A_{\alpha}, 
       \ \sigma(B_{\alpha})=-B_{\alpha},\\ 
     \alpha\in\varDelta_3^+\Longrightarrow\sigma(A_{\alpha})=-B_{\alpha},
      \ \sigma(B_{\alpha})=A_{\alpha},
   \end{array}\right.
     \label{eqn:62}
\end{equation}
and 
\begin{equation}
  \sigma^2|_{\mathfrak m_1}=-{\rm Id}_{\mathfrak
   m_1},\quad\sigma^2|_{\mathfrak m_3}=-{\rm Id}_{\mathfrak m_3}.
     \label{eqn:63}
\end{equation}
If $X\in\mathfrak m_1^+$, then by (\ref{eqn:63}), we have 
$$\tau(\sigma(X))=\sigma^{-1}(\tau(X))=\sigma^3(X)=-\sigma(X).$$
Thus we have $\sigma(\mathfrak m_1^+)\subset\mathfrak m_1^-$. Similarly, we
get $\sigma(\mathfrak m_1^-)\subset\mathfrak m_1^+$. Therefore it follows that
\begin{equation}
  \sigma(\mathfrak m_1^+)=\mathfrak m_1^-,\quad\sigma(\mathfrak
  m_1^-)=\mathfrak m_1^+.
    \label{eqn:64}
\end{equation}
 Similarly, we obtain 
\begin{equation}
  \sigma(\mathfrak m_3^+)=\mathfrak m_3^-, \quad\sigma(\mathfrak
   m_3^-)=\mathfrak m_3^+.  
     \label{eqn:65}
\end{equation}
By a straightforward computation we have
\begin{equation}
   [\sqrt{-1}H,A_{\alpha}]=\alpha(H)B_{\alpha},\quad
   [\sqrt{-1}H,B_{\alpha}]=-\alpha(H)A_{\alpha}.
     \label{eqn:66}
\end{equation}
Put
$X_1=\sum_{\alpha\in\varDelta_1^+}(a_{\alpha}A_{\alpha}+b_{\alpha}B_{\alpha})
\in\mathfrak m_1$. Then by (\ref{eqn:62}), we have 
$\sigma(X_1)=\sum_{\alpha\in\varDelta_1^+}(a_{\alpha}B_{\alpha}-b_{\alpha}
A_{\alpha})$. Using (\ref{eqn:66}), it is easy to see that
$$
    [\sqrt{-1}Z,X_1]=\sigma(X_1).
$$
Similarly, we get
$$
    [\sqrt{-1}Z,X_3]=-3\sigma(X_3),\ \ 
[\sqrt{-1}Z,[\sqrt{-1}Z,X_2]]=-4X_2,
$$
for $X_j\in\mathfrak{m}_j$ ($j=2,3$). Therefore it follows from (\ref{eqn:63}) 
that 
\begin{equation}
  \begin{array}{ll}
    [Z,X_1\pm\sqrt{-1}\sigma(X_1)]=\mp(X_1\pm\sqrt{-1}\sigma(X_1)),\\ 
  {}[Z,X_3\pm\sqrt{-1}\sigma(X_3)]=\pm 3(X_3\pm\sqrt{-1}\sigma(X_3)).
  \end{array}
      \label{eqn:67}
\end{equation}
Note that $X_s\pm\sqrt{-1}\sigma(X_s)\in\mathfrak{g}^*$ for 
$X_s\in\mathfrak{m}_s^+$ ($s=1,3$) from (\ref{eqn:64}) and (\ref{eqn:65}). 
Moreover, $Y_2:=[\sqrt{-1}Z,X_2]\ne 0$ and $Y_2\in\mathfrak{m}_2^-$ for 
$X_2\in\mathfrak{m}_2^-$, and 
\begin{equation}
    [Z,X_2\pm \frac{1}{2}\sqrt{-1}Y_2]=\mp 2(X_2\pm \frac{1}{2}\sqrt{-1}Y_2). 
      \label{eqn:68}
\end{equation}
Consequentry, from (\ref{eqn:67}) and (\ref{eqn:68}) the lemma is
proved. \hfill$\Box$ 
\par

Now, we are in a position to prove the following proposition which classifies 
involutions preserving $\mathfrak h$ for this case. 
\begin{prop} $(1)$ Let $\mathfrak{g}^*=\sum_{p=-3}^3\mathfrak{g}^*_p$ be a 
graded simple Lie algebra of the third kind with a grade-reversing Cartan 
involution $\tau$, which is corresponding to a partition 
$\{\varPi_0,\varPi_1\}$ of $\varPi=\{\lambda_1,\ldots,\lambda_l\}$ such 
that $\varPi_1=\{\lambda_i\}$ with $n_i=3$. Put 
$\sigma={\rm Ad}(\exp(\pi/2)\sqrt{-1}h_i)$. Then $\sigma$ is an automorphism 
of order $4$ on the compact dual $\mathfrak{g}$ of $\mathfrak{g}^*$ such that 
$\dim\mathfrak{z}=1$ and $\tau\circ\sigma=\sigma^{-1}\circ\tau$. 
\par
$(2)$ Let $\sigma={\rm Ad}(\exp(\pi/2)\sqrt{-1}K_i)$ for some 
$\alpha_i\in\varPi(\mathfrak{g}_{\mathbb{C}},\mathfrak{t}_{\mathbb{C}})$ 
with $m_i=3$. Then for each involution $\tau$ of $\mathfrak{g}$ satisfying 
$\tau\circ\sigma=\sigma^{-1}\circ\tau$, there exists $\theta\in
{\rm Aut}(\mathfrak{g})$ such that $\theta\circ\sigma\circ\theta^{-1}$ and 
$\theta\circ\tau\circ\theta^{-1}$ are obtained from a graded Lie algebra 
by the method described in $(1)$.
\label{prop:61}
\end{prop}

 Proof. We have proved (1) in the above. 
\par
Now we prove (2). For each $\sigma={\rm Ad}(\exp(\pi/2)\sqrt{-1}K_i)$ and 
$\tau$ with $\tau\circ\sigma=\sigma^{-1}\circ\tau$, it follows from 
Lemma \ref{lem:63} that there exists a graded Lie algebra 
$\mathfrak{g}^*=\sum_{p=-3}^{3}\mathfrak{g}^*_p$ with the characteristic 
element $Z:=K_i$ such that $\tau$ is the Cartan involution. As above, 
let $\mathfrak{g}^*=\mathfrak{k+p}^*$ be the Cartan decomposition of 
$\mathfrak{g}^*$ corresponding to $\tau$ and let $\mathfrak{a}$ be a maximal 
abelian subspace of $\mathfrak{p}^*$ such that $Z\in\mathfrak{a}$. Moreover, 
let $\mathfrak{t}^*$ be a Cartan subalgebra of $\mathfrak{g}^*$ containing 
$\mathfrak{a}$ equipped with a compatible ordering. 
By Lemma \ref{lem:63}, we have $\lambda(Z)=0,\pm 1,\pm2$ or $\pm3$ for 
any $\lambda\in\varDelta:=\varDelta(\mathfrak{g}^*,\mathfrak{a})$.
If $\varDelta$ is a reduced root system, then from Lemma \ref{lem:24} together 
with Lemma 2.4 of \cite{WG} there exists 
$w\in W(\mathfrak{g}^*,\mathfrak{a})$ such that 
$$\frac{1}{4}w(Z)=\frac14 h+T.$$
Here $T$ is an element in $\mathfrak{a}$ satisfying 
$\lambda(T)\in{\mathbb Z}$ for any $\lambda\in\varDelta$, and $h$ is one of the following: 
$$h_p,\ h_{q_1}+h_{q_2}, \ 2h_{r_1}+h_{r_2},\ h_{s_1}+h_{s_2}+h_{s_3}, 
$$
with $n_p=1, 2, 3$ or $4$, $(n_{q_1},n_{q_2})=(1,1)$, $(1,2)$ or $(2,2)$, 
$n_{r_1}=n_{r_2}=1$ and $n_{s_1}=n_{s_2}=n_{s_3}=1$. If $\varDelta$ is 
a nonreduced root system, then 
$\varDelta':=\{\lambda\in\varDelta \ ; \ 2\lambda\not\in\varDelta\}$ is a
reduced root system of type $B_l$ with the fundamental root system $\varPi$. 
Applying Lemma \ref{lem:24} together with Lemma 2.4 of \cite{WG} to 
$\varDelta'$, we can see that there exists $w\in W(\mathfrak{g}^*,
\mathfrak{a})$ such that $(1/4) w(Z)=(1/4) h+T$ with $\lambda(T)\in
\mathbb Z$ for any $\lambda\in\varDelta'$ and $h$ is one of 
$$h_a, \ h_b+h_c,\ \ n_a=n_b=n_c=2.$$
Hence we may assume that there exists $\theta\in{\rm Int}(\mathfrak{k})$ 
such that 
$$
\theta\circ\sigma\circ\theta^{-1}={\rm Ad}(\exp\frac{\pi}{2}\sqrt{-1}h).
$$
Note that $\theta\circ\tau\circ\theta^{-1}=\tau$ because 
$\theta\in{\rm Int}(\mathfrak{k})$. 
\par

Next, we shall prove that $h=h_p$ for some $\lambda_p\in\varPi$ with $n_p=3$.
In the case where $h=h_p$ with $n_p=1$, there exists a unique 
$\alpha_{i_p}\in\varPi(\mathfrak{g}_{\mathbb{C}},\mathfrak{t}_{\mathbb{C}})$ 
such that $m_{i_p}=1$ and $\alpha_{i_p}|_{\mathfrak{a}}=\lambda_p$. 
Therefore by Lemma \ref{lem:25} together with Remark 2.2 we have 
$h_p=K_{i_p}$ and $(\mathfrak{g},\mathfrak{h})$ is a symmetric pair, which 
also contradicts Remark 6.1. 
Similarly, if $h=h_p$, $h_{q_1}+h_{q_2}$, $2h_{r_1}+h_{r_2}$ or $h_a$ with 
$n_p=n_a=2$ and $n_{q_1}=n_{q_2}=n_{r_1}=n_{r_2}=1$, then a pair 
$(\mathfrak{g},\mathfrak{h})$ is $3$-symmetric, which contradicts Remark 6.1. 
\par

In the case where $h=h_{s_1}+h_{s_2}+h_{s_3}$ with 
$n_{s_1}=n_{s_2}=n_{s_3}=1$, there exist unique $\alpha_{i_1}$, 
$\alpha_{i_2}$, $\alpha_{i_3}\in\varPi(\mathfrak{g}_{\mathbb{C}},
\mathfrak{t}_{\mathbb{C}})$ such that $\alpha_{i_k}|_{\mathfrak{a}}
=\lambda_{s_k}$ ($k=1,2,3$). Then we have $h=K_{i_1}+K_{i_2}+K_{i_3}$ 
and hence $\dim\mathfrak{z}=3$, which is a contradiction. \par

In the case where $h=h_{q_1}+h_{q_2}$, then we obtain 
$$h=K_{i_1}+K_{i_2} \ \mbox{or}\ K_{i_1}+K_{j_1}+K_{j_2}.$$
Here $\alpha_{i_k}|_{\mathfrak{a}}=\lambda_{q_k}$, $m_{i_1}=1$, $m_{i_2}=2$, 
or $\alpha_{i_1}|_{\mathfrak{a}}=\lambda_{q_1}$, $\alpha_{j_k}|_{\mathfrak{a}}
=\lambda_{q_k}$, $m_{i_k}=1$ ($k=1,2$). Therefore by Remark 2.2 we 
have $\dim\mathfrak{z}\not=1$. 
\par

In the case where $h=h_p$ with $n_p=4$, then we have
\begin{eqnarray*}
&{\rm (i)}& h=K_{i_p} \ \mbox{with}\ m_{i_p}=4,\ \mbox{or}\\
&{\rm (ii)}& h=K_{i_1}+K_{i_2}\ \mbox{with}\ \alpha_{i_k}|_{\mathfrak{a}}
=\lambda_{p}, \ m_{i_k}=2 \ (k=1,2).
\end{eqnarray*} 
For the case (i), it follow from Remark 2.2 that 
$\dim\mathfrak{z}=0$. For the case (ii), it is easy to see that the center 
$\mathfrak{z}(\mathfrak{g}^{\tau_{(1/2)h}})$ of $\mathfrak{g}^{\tau_{(1/2)h}}$ 
coincides with 
\begin{equation}
\mathfrak{z}(\mathfrak{g}^{\tau_{(1/2)h}})
={\mathbb{R}}\sqrt{-1}(K_{i_1}-K_{i_2}).
\label{eqn:68}
\end{equation}
Note that if $\sigma={\rm Ad}(\exp(\pi/2)\sqrt{-1}K_i)$ with 
$m_i=3$, then the center $\mathfrak{z}$ of $\mathfrak{h}$ coincides with 
$\mathbb{R}\sqrt{-1}K_i$ as mentioned before. It is easy to see that 
$\mathfrak{h}$ is the centralizer of $\mathfrak{z}$ in $\mathfrak{g}$. 
However, $\mathfrak{g}^{\tau_{(1/2)h}}$ is not the 
centralizer of $\mathfrak{z}(\mathfrak{g}^{\tau_{(1/2)h}})$. 
Indeed, let $\alpha=\sum_{j}k_j\alpha_j$ be a root satisfying 
$k_{i_1}=k_{i_2}=1$. Since $\alpha(h)=2$ and $\alpha(K_{i_1}-K_{i_2})=0$, 
we obtain 
$$[\sqrt{-1}(K_{i_1}-K_{i_2}),A_{\alpha}]=0,\ \ 
\tau_{(1/2)h}(A_{\alpha})=-A_{\alpha},$$
which implies that $A_{\alpha}$ belongs to the centralizer of 
$\mathfrak{z}(\mathfrak{g}^{\tau_{(1/2)h}})$ and 
$A_{\alpha}\not\in\mathfrak{g}^{\tau_{(1/2)h}}$. 
Hence $\sigma$ is not conjugate to $\tau_{(1/2)h}$. \par

Finally, consider the case where $h=h_j+h_k$ with $n_j=n_k=2$. In this case, 
we have $h=K_{i_j}+K_{i_k}$ with $m_{i_j}=m_{i_k}=2$, or 
$h=K_{j_1}+K_{j_2}+K_{i_k}$ with $m_{j_1}=m_{j_2}=1$, $m_{i_k}=2$. 
By the same argument as (i) above, the first case is impossible. Moreover, 
if $h=K_{j_1}+K_{j_2}+K_{i_k}$, then the center of 
$\mathfrak{g}^{\tau_{(1/2)h}}$ coincides with 
$${\mathbb{R}}\sqrt{-1}(2K_{j_1}-K_{i_k})
+{\mathbb{R}}\sqrt{-1}(K_{j_2}-K_{i_k}),$$
since $\mathfrak{g}^{\tau_{(1/2)h}}$ is  generated by $\mathfrak{t}$ and 
$\{A_{\alpha},\ B_{\alpha} \ ; \ \alpha(h)\equiv0\pmod{4}\}$.
However, this is a contradiction. \par
Consequently we obtain $h=h_p$ with $n_p=3$ which completes the proof 
of $(2)$ of the proposition.  \hfill$\Box$


\section{The case where $\dim\mathfrak z = 0$ and $\tau\circ\sigma =
 \sigma\circ\tau$}
 In this section we consider the case where $\dim\mathfrak z = 0$ and
$\tau\circ\sigma = \sigma\circ\tau$. In this case, it follows from 
Proposition \ref{prop:41} that $\tau |_{\mathfrak{t}} 
={\rm Id}_{\mathfrak{t}}$ or $\tau$ is of Type I in Table 1.
\par
 First we consider the Type I in Table 1. From Section 5 there exists an
 automorphim $\varphi$ satisfying (\ref{eqn:510}). We note that $\varphi$ 
is an involution of Type I in Table 1. 
Let $\mathfrak{t}_{\pm}$ be the $(\pm 1)$-eigenspaces of $\varphi$. 
Then we have
$$\begin{array}{l}
   \mathfrak{t}_+ = {\rm span}\{K_1 + K_6 - 2K_7, \ K_2 - K_7, \ K_3+K_5-
   3K_7, \  K_4 -2K_7 \}, \\
   \mathfrak{t}_- = {\rm span}\{-K_1 + K_6 , \ -K_3 + K_5, \ K_7 \}.
\end{array}$$ 
For any involution $\tau$ of Type I, it follows from Proposiion 5.3 of
Chapter IX of \cite{Helg} that there exists $\sqrt{-1} h 
\in\mathfrak{t}$ such that $\tau =\varphi\circ \tau_h$. We put $h
=h_+ + h_-, \ h_{\pm}\in\sqrt{-1}\mathfrak{t}_{\pm}$. Then since $\tau^2 ={\rm Id}$, we can write
\begin{equation}
   h_+  =  k_1(K_1 + K_6 - 2K_7) + k_2(K_2 + K_7) + k_3(K_3 + K_5 - 3K_7) +
   k_4(K_4 - 2K_7), \ \	k_i\in\mathbb Z.
     \label{eqn:71}
\end{equation}
As in the case of Type IV in Section 5, we may assume 
$\tau(E_{\alpha_1})=E_{\alpha_6}$, $\tau(E_{\alpha_3})=E_{\alpha_5}$ and 
$\tau(E_{\alpha_7})=E_{\alpha_0}$. Indeed, for example, if $\tau(E_{\alpha_1})
=b_1E_{\alpha_6}$ for some $b_1\in\mathbb{C}$ with $|b_1|=1$, then 
using $h_-=k(-K_1+K_6)$ with $e^{2\pi\sqrt{-1}k}=b_1$, we have 
$(\tau_{h_-})^{-1}\circ\tau\circ\tau_{h_-}(E_{\alpha_1})=E_{\alpha_6}$. 
\par  
Using (\ref{eqn:71}), by an argument similar to
the case of Type IV in Section 5 we can prove that $\tau$ is conjugate within 
${\rm Int}(\mathfrak{h})$ to one of the following involutions:
$$
\varphi, \quad  \varphi\circ\tau_{K_2}, \quad \varphi\circ\tau_{K_4}, \quad
\varphi\circ\tau_{K_2}\circ\tau_{K_4}.  
$$
Note that $\mathfrak{su}_{\alpha_2}(2)\subset\mathfrak{h\cap g}^{\varphi}$, 
and hence $t_{\alpha_2}\in{\rm Int}(\mathfrak{h\cap g}^{\varphi})$. Therefore 
we have $\varphi\circ t_{\alpha_2}=t_{\alpha_2}\circ\varphi$. Moreover, since 
$t_{\alpha_2}(K_2) = -K_2 + K_4$, it follows that 
$\varphi\circ\tau_{K_2}$ is conjugate within ${\rm Int}(\mathfrak h)$ to 
$\varphi\circ\tau_{K_2}\circ\tau_{K_4}$. 
\par
Put $\nu:=\tau|_{\mathfrak{t}}$. It is easy to see that the set 
$\varDelta_{\nu}^+$ of positive
roots $\alpha$ satisfying $\nu(\alpha)=\alpha$ coincides with 
$$\varDelta_{\nu}^+ =\left\{
\begin{tabular}{p{300pt}}
  $\alpha_2, \ \alpha_4, \ \alpha_2 + \alpha_4, \ \alpha_3 +
  \alpha_4 + \alpha_5 , \ \alpha_2 + \alpha_3 +\alpha_4 + \alpha_5, \
  \alpha_2 + \alpha_3 + 2\alpha_4 + \alpha_5, \ \alpha_1 + \alpha_3
  +\alpha_4 + \alpha_5 + \alpha_6, \ \alpha_1 + \alpha_2 + \alpha_3
  + \alpha_4 + \alpha_5 + \alpha_6, \ \alpha_1 + \alpha_2 + \alpha_3
  + 2\alpha_4 + \alpha_5 + \alpha_6, \ \alpha_1 + \alpha_2 +
  2\alpha_3 + 2\alpha_4 + 2\alpha_5 + \alpha_6, \ \alpha_1 + \alpha_2
  + 2\alpha_3 + 3\alpha_4 + 2\alpha_5 + \alpha_6, \ \delta$ 
\end{tabular}
\right\}.$$
Using (\ref{eqn:510}), we can check that $\varphi(E_{\alpha})=E_{\alpha}$ 
for any $\alpha\in\varDelta_{\nu}^+$. For example 
$$\varphi ([E_{\alpha_5},[E_{\alpha_3},E_{\alpha_4}]]) =
[E_{\alpha_3},[E_{\alpha_5},E_{\alpha_4}]] =
[E_{\alpha_5},[E_{\alpha_3},E_{\alpha_4}]],$$
and thus $\varphi(E_{\alpha_3+\alpha_4+\alpha_5})=
E_{\alpha_3+\alpha_4+\alpha_5}$. 
Hence it follows from Lemma \ref{lem:41} that 
$\dim\mathfrak{g}^{\varphi}=79$. By using the classification of
symmetric spaces,  we get $\mathfrak{g}^{\varphi}\cong
E_6\oplus\mathbb{R}$.\par

 The number of subsets $\{\alpha , \beta\}$ such that
$\alpha\in\varDelta^+(\mathfrak{g}_{\mathbb{C}},\mathfrak{t}_{\mathbb{C}}), \ 
\tau(\alpha) = \beta, \ \alpha\ne\pm\beta$ and
$\alpha(K_4)\equiv 0\pmod 4$ is $6$. Furthermore
$\alpha\in\varDelta^+(\mathfrak{g}_{\mathbb{C}},\mathfrak{t}_{\mathbb
C})$ such that $\tau(\alpha)=\alpha$ and $\alpha(K_4)\equiv 0\pmod 4$ is
only $\alpha_2$. 
Since $\dim\mathfrak{t}_+ = 4$, we get 
$$\dim(\mathfrak{h\cap g}^{\varphi})=4 + ((6+1)\times 2) = 18.$$ 
Therefore we get $\mathfrak{h\cap g}^{\varphi}\cong D_8\oplus C_1$.\par

 Similarly as above we can compute $\dim(\mathfrak{h\cap g}^{\tau})$ and
$\dim\mathfrak{g}^{\tau}$ for the other types.\\
\smallskip\par 

Next we consider the case $\tau|_{\mathfrak{t}} = {\rm Id}$. First we suppose 
that $\mathfrak{g}$ is of type $\mathfrak{e}_8$ and 
$\sigma={\rm Ad}(\exp(\pi/2)\sqrt{-1}K_3)$. Then by (\ref{eqn:42}) (i), 
we have 
$$\mathfrak{h} \cong A_7\oplus A _1.$$
Furthermore a maximal abelian subalgebra $\mathfrak{t}$ is decomposed into 
$\mathfrak{t} =(A_7\cap\mathfrak{t})\oplus(A_1\cap\mathfrak{t}).$ 
 Hence we can write 
$$\tau = \tau_{T_1}\circ\tau_{T_2}, \ \
\sqrt{-1}T_1\in A_7\cap\mathfrak{t}, \ \sqrt{-1}T_2\in 
A _1\cap\mathfrak{t}.$$
We define $v_i\in\sqrt{-1}(A_7\cap\mathfrak{t}), \ i\in\Lambda
 :=\{0,2,4,5,6,7,8\}$ and $v_1\in\sqrt{-1}(A_1\cap\mathfrak{t})$ 
by $\alpha_i (v_j) = \delta_{ij}$. 
Since $(\tau_{T_1}|_{A_7})^2 ={\rm Id}_{A_7}$ and 
$(\tau_{T_2}|_{A_1})^2 ={\rm Id}_{A_1}$, it follows from Lemma
\ref{lem:22} and Remark 2.1 that there exist $\mu_1\in{\rm Int} (A_7)$ and
 $\mu_2\in{\rm Int} (A_1)$ such that
\begin{equation}
  \mu_1(T_1)\equiv\left\{
   \begin{array}{ll}
    0  \ &{\rm mod} \ 2\varPi_{A_7},\vspace{0.1cm}\\ 
    v_i \ &{\rm mod} \ 2\varPi_{A_7} \ (i\in\Lambda),
   \end{array}
  \right.\quad
  \mu_2(T_2)\equiv\left\{
   \begin{array}{ll}
    0  \ &{\rm mod} \ 2\varPi_{A_1},\vspace{0.1cm}\\ 
    v_1 \ &{\rm mod} \ 2\varPi_{A_1},
   \end{array}
  \right.
\end{equation}
where $\varPi_{A_l}$ denotes the fundamental root system of Type $A_l$. 
Therefore considering Lemma \ref{lem:23} we may assume 
\begin{equation}
  T_1=\left\{
   \begin{array}{ll}
    2m_0v_0+2m_2v_2 +2m_3v_3 +\cdots +2m_8v_8,\vspace{0.1cm}\\ 
    v_i + 2m_0v_0+2m_2v_2 +2m_3v_3 +\cdots +2m_8v_8,
   \end{array}
  \right.\quad
  T_2=\left\{
   \begin{array}{ll}
    2m_1v_1,\vspace{0.1cm}\\ 
    v_1 + 2m_1v_1,
   \end{array}
  \right.
\end{equation}
where $i=2,4,5,6$ and $m_0,m_1,m_2,m_4,\ldots
,m_8\in\mathbb Z$. 
Consequently $\tau$ is
conjugate within ${\rm Int}(\mathfrak h)$ to one of the following 
automorphisms:
\begin{equation}
  \left\{
   \begin{array}{l}
     {\rm Ad}(\exp\pi\sqrt{-1}(2m_0v_0 + 2m_1v_1 + 2m_2v_2 + 2m_4v_4
     +\cdots + 2m_8v_8)),
     \vspace{0.1cm}\\  
     {\rm Ad}(\exp\pi\sqrt{-1}(v_i + 2m_0v_0 + 2m_1v_1 + 2m_2v_2 + 2m_4v_4
     +\cdots + 2m_8v_8)), 
     \vspace{0.1cm}\\  
     {\rm Ad}(\exp\pi\sqrt{-1}(v_1 + v_j + 2m_0v_0 + 2m_1v_1 + 2m_2v_2 +
     2m_4v_4 +\cdots + 2m_8v_8)),  
   \end{array}\right.
     \label{eqn:74}
\end{equation}
where $i= 1,2,4,5,6, \ j = 2,4,5,6$ and $m_0,m_1,m_2,m_4,\ldots
,m_8\in\mathbb Z$. \par
Now we compute $v_j$. Put $v_1 = \sum_{i=1}^8 a_i K_i, \ a_i\in\mathbb
 R$. Since $A_1\cap\mathfrak{t}={\mathbb{R}}\sqrt{-1}H_{\alpha_1}$ 
and 
$$A_1\cap\mathfrak{t} = \{\sqrt{-1} H \in\mathfrak{t} \ ; \
\alpha_j(H)=0,\ j=0,2,4,5,6,7,8\},$$
we have $a_1=1$, $a_2 = a_4= \cdots = a_8 =0$ and $a_1+2a_3=0$. 
Hence we obtain $v_1 = K_1 -(1/2)K_3$. \par 
Moreover, since $A_7\cap\mathfrak{t} = \{\sqrt{-1} H \in\mathfrak{t} \ ;
\ \alpha_1(H)= 0 \}$,
we can put $v_i = \sum_{k=2}^8 b_k^i K_k, \ b^i_k\in\mathbb{R}, \
i\in\Lambda$. Then computing simultaneus equations
$\alpha_i(v_j)=\delta_{ij}, \ i,j\in\Lambda$, we obtain
\begin{equation}
  \begin{array}{ll}
    &\displaystyle v_0 = - \frac{1}{4}K_3, \ v_1 =  K_1 -
     \frac{1}{2}K_3, \ v_2 = K_2 -\frac{3}{4} K_3, \ 
     v_4 =   K_4 - \frac{3}{2}K_3, \vspace{0.2cm}\\ 
    &\displaystyle v_5 = K_5 -\frac{5}{4}K_3, \  v_6 = K_6 - K_3, \ 
    v_7 = K_7 - \frac{3}{4}K_3 , \ v_8 = K_8 - \frac{1}{2}K_3.
  \end{array}
\label{eqn:75}
\end{equation}
Thus (\ref{eqn:74}) implies that $\tau$ is conjugate within 
${\rm Int}(\mathfrak h)$ to one of the following: 
\begin{equation}
 \tau_{mK_3}, \ \ \tau_{v_i + mK_3}, 
\ \ \tau_{v_1 + v_j + mK_3},  
\label{eqn:76}
\end{equation}
where $m = -((1/2)m_0 + m_1 + (3/2)m_2 + 3m_4 +
(5/2)m_5 + (3/2)m_7 + m_8)$. From (\ref{eqn:75}) if $i,j=2,5$, then
$\tau^2\ne {\rm Id}$. Therefore 
$i=1,4,6$ and $j=4,6$. Hence $\tau$ is conjugate within ${\rm
Int}(\mathfrak h)$ to some $\tau_h$ where $h$ is one of the following: 
$$K_i,\ \ K_3+K_j,\ \ K_1+K_k,\ \ K_1+K_3+K_k,$$
where $i=1,3,4,6$, $j=1,4,6$ and $k=4,6$. 
If $h = K_1$, then $\mathfrak{g}^{\tau_{K_1}}\cong D_8$ (cf. 
Theorem 5.15 of Chapter X of \cite{Helg}). Furthermore 
$$\mathfrak {h\cap g}^{\tau_{K_1}} = \mathfrak{t}
\oplus\sum_{\substack{\alpha\in\varDelta^+(\mathfrak{g}_{\mathbb{C}},
\mathfrak{t}_{\mathbb{C}}) \\ \alpha(K_3)\equiv 0 \bmod 4 
\\ \alpha(K_1)\equiv 0 \bmod 2}} (\mathbb{R} A_{\alpha} + \mathbb{R} 
B_{\alpha})\subset\mathfrak h \left(\cong \mathfrak
 {su}(8)  \oplus \mathfrak{su}(2)\right).$$
In this case, $\tau_{K_1}|_{A_7}={\rm Id}$ and 
$A_1^{\tau_{K_1}}\cong\mathbb{R}$, and hence 
$$\mathfrak {h\cap k}\cong A_7 \oplus \mathbb {R} .$$
Similarly as above, we can get $(\mathfrak{g}^{\tau}, 
\mathfrak{h\cap g}^{\tau})$ for each $\tau=\tau_h$.

Now we consider the reflection $t_{\alpha_1}\in
{\rm Int}(\mathfrak{su}_{\alpha_1}(2))\subset{\rm Int}(\mathfrak h)$. 
It is easy to check that $t_{\alpha_1}$ maps
 $K_1\mapsto -K_1 + K_3, \ K_1 + K_4\mapsto -K_1
 + K_3 + K_4, \ K_1 + K_6\mapsto -K_1 + K_3 + K_6$ and $K_3\mapsto
 K_3$. Therefore we have $\tau_{K_1}\approx\tau_{K_1 + K_3}, \
 \tau_{K_1 + K_4}\approx\tau_{K_1 + K_3 + K_4}$ and $\tau_{K_1 +
 K_6}\approx\tau_{K_1 + K_3 + K_6}$, where we write
 $\tau_H\approx\tau_{H'}$ if $\tau_H$ is conjugate to $\tau_{H'}$ within
 ${\rm Aut}_{\mathfrak h}(\mathfrak{g})$.\vspace{0.5cm}\par
Next, we consider the case where $\mathfrak{g}$ is of type $\mathfrak{e}_8$ 
and $\sigma={\rm Ad}(\exp(\pi/2)\sqrt{-1}K_6)$. Then by (\ref{eqn:42}) (ii), 
we have $\mathfrak{h}\cong A_3\oplus D_5$. By a computation similar to the 
above case, we obtain 
\begin{equation}
  \begin{array}{l}
    \displaystyle v_0=- \frac{1}{4}K_6, \ v_1=K_1-\frac{1}{2}K_6, \
     v_2=K_2 -\frac{3}{4} K_6, \ v_3=K_3-K_6, \vspace{0.2cm}\\  
    \displaystyle  v_4=K_4-\frac{3}{2}K_6, \
     v_5=K_5-\frac{5}{4}K_3, \ v_7=K_7-\frac{3}{4}K_6 , 
  \ v_8=K_8-\frac{1}{2}K_6.
  \end{array}
   \label{eqn:77}
\end{equation}
Then, considering Lemma \ref{lem:23}, $\tau$ is conjugate within 
${\rm Int}(\mathfrak h)$ to one of the following automorphisms:  
\begin{equation}
\tau_{mK_6},\ \ \tau_{v_a + mK_6},\ \ \tau_{v_b + v_c + mK_6},
\label{eqn:78}
\end{equation}
where $a =1,2,3,7,8, \  b = 1,2,3, \ c = 7,8$, and $m$ is equal to that
of the above case. Since $\tau^2={\rm Id}$, it follows from (\ref{eqn:77}) 
that $\tau$ is conjugate within ${\rm Int}(\mathfrak h)$ to some $\tau_h$ 
where $h$ is one of the following:  
\begin{eqnarray*}
& & K_i,\ \  K_j + K_6, \ \ K_2 + K_7,\ \ K_k + K_8,\ \ K_k + K_6 + K_8, \ \
K_2 + K_6 + K_7,
\end{eqnarray*}
where $i=1,3,6,8$, $j=1,3,8$ and $k=1,3$. 
By a computation similar to the above case, we obtain
$\mathfrak{g}^{\tau}$ and $\mathfrak{h\cap g}^{\tau}$.
We put
\begin{equation}
  \begin{array}{ll}
    \beta_1 :=\left(
      \begin{array}{lllllll}
         &   &   &   & 3 &   &   \\
       1 & 2 & 4 & 5 & 6 & 4 & 2 \\
      \end{array}\right),& \ 
    \beta_2 := \left(
      \begin{array}{lllllll}
         &   &   &   & 1 &   &   \\
       0 & 0 & 0 & 1 & 2 & 2 & 1 \\
      \end{array}\right), \\ 
    \beta_3:=\left(
      \begin{array}{lllllll}
         &   &   &   & 1 &   &   \\
       0 & 0 & 0 & 1 & 2 & 1 & 0 \\
      \end{array}\right). \\ 
  \end{array}
\end{equation}
 It is easy to check that
$$t_{\beta_1}\circ t_{\alpha_8}(K_6)=-3K_6 + 4K_7, \ \  t_{\beta_2}\circ
t_{\alpha_1}(K_6)=K_6, \ \ t_{\beta_3}\circ t_{\alpha_3}(K_6)=K_6,$$
and so $t_{\beta_1}\circ t_{\alpha_8}$, $t_{\beta_2}\circ t_{\alpha_1}$, 
$t_{\beta_3}\circ t_{\alpha_3}\in{\rm Aut}_{\mathfrak h}(\mathfrak{g})$. 
Moreover we have 
$$\begin{array}{ll}
 t_{\beta_1}\circ t_{\alpha_8}(K_8) = -K_6 + 2K_7 - K_8, &
 t_{\beta_2}\circ t_{\alpha_1}(K_1) = -K_1 + K_6, \\
 t_{\beta_2}\circ t_{\alpha_1}(K_1 + K_8) = -K_1 + K_6 + K_8, &
 t_{\beta_3}\circ t_{\alpha_3}(K_3) = 2K_1 - K_3 + K_6, \\
 t_{\beta_3}\circ t_{\alpha_3}(K_3 + K_8) = 2K_1 - K_3 + K_6 + K_8. &
\end{array}$$
Therefore we have 
$$
\tau_{K_6 + K_8}\approx\tau_{K_8}, \ \ \tau_{K_1 +
K_6}\approx\tau_{K_1}, \ \ \tau_{K_3 + K_6}\approx\tau_{K_3}, \ \
\tau_{K_1 + K_6 + K_8}\approx\tau_{K_1 + K_8},  \ \ \tau_{K_3 + K_6 +
K_8}\approx\tau_{K_3 + K_8}.$$
\vspace{0.2cm}\par 
For the case where $\mathfrak{g}=\mathfrak{e}_7$ and 
$\sigma={\rm Ad}(\exp(\pi/2)\sqrt{-1}K_4)$, we can check that 
$\tau$ is conjugate within ${\rm Int}(\mathfrak h)$ to some $\tau_h$
where $h$ is one of the following:  
\begin{eqnarray*}
& & K_i, \ \ K_j + K_k, \ \ K_3 + K_7, \ \ K_l + K_m + K_n, \\
& & K_p + K_3 + K_7,\ \ K_1 + K_2 + K_4 + K_6, \ \ K_2 + K_3 + K_4 + K_7,
\end{eqnarray*}
where $i=1,2,4,6$, $j,k=1,2,4,6$ $(j<k)$, $l,m,n=1,2,4,6$ ($l<m<n$) and 
$p=2,4$. Using the reflection $t_{\alpha_2}\in{\rm Int}(\mathfrak{h})$ we have 
$$
  \begin{array}{l}
    \tau_{K_2 + K_4}\approx\tau_{K_2}, \ \ \tau_{K_1 + K_2 +
     K_4}\approx\tau_{K_1 + K_2}, \ \ \tau_{K_2 + K_4 +
     K_6}\approx\tau_{K_2 + K_6}, \\ 
    \tau_{K_1 + K_2 + K_4 + K_6}\approx\tau_{K_1 + K_2 + K_6}, \ \
     \tau_{K_2 + K_3 + K_4 + K_7}\approx\tau_{K_2 + K_3 + K_7}, 
  \end{array} 
$$
and since $t_{\alpha_5 + \alpha_6 + \alpha_7}\circ t_{\alpha_6}\in{\rm 
Int}(\mathfrak h)$, we have
$$\tau_{K_6}  \approx  \tau_{K_4 + K_6}.$$ 
Furthermore put $\gamma_1 := \alpha_1 + 2\alpha_2 + 2\alpha_3 + 4\alpha_4 +
3\alpha_5 + 2\alpha_6 + \alpha_7$. Then $t_{\gamma_1}\circ t_{\alpha_1}\in{\rm 
Int}(\mathfrak h)$ gives the following conjugations:
$$
\tau_{K_1 + K_4}\approx\tau_{K_1}, \ \ \tau_{K_1 +
K_4 + K_6}\approx\tau_{K_1 + K_6}.
$$\par
Finally we consider an involution 
$\varphi\in{\rm Aut}_{\mathfrak h}(\mathfrak{g})$ (see (\ref{eqn:510})). 
Then it is easy to see that
\begin{eqnarray*}
& & \varphi(K_1)= K_6-2K_7, \ \varphi(K_2)= K_2 - 2K_7, \ 
\varphi(K_3)= K_5-3K_7, \ \varphi(K_4) = K_4 + 4K_7,\\
& & \varphi(K_5)=K_3-3K_7, \ \varphi(K_6)=K_1-2K_7, \ \varphi(K_7)=-K_7,
\end{eqnarray*}
and therefore $\varphi$ gives the following conjugations: 
$$
\tau_{K_1}\approx\tau_{K_6}, \ \ \tau_{K_1+K_2}\approx\tau_{K_2 + K_6}.
$$\vspace{0.5cm}\par
 For the case where $\mathfrak{g}=\mathfrak{f}_4$ and 
$\sigma={\rm Ad}(\exp(\pi/2)\sqrt{-1}K_3)$, we can check that $\tau$ is 
conjugate within ${\rm Int}(\mathfrak h)$ to some $\tau_h$ where $h$ is one 
of the following: 
$$
K_i, \ \ K_j + K_k,\ \ K_1 + K_3 + K_4.
$$
Here $i=1,3,4$ and $j,k=1,3,4$ ($j<k$).
Using the reflection $t_{\alpha_4}\in{\rm Int}(\mathfrak{h})$ and 
$t_{\alpha_1 +\alpha_2 + 2\alpha_3 + \alpha_4}$ we have
$$
\tau_{K_3 + K_4}\approx\tau_{K_4}, \ \ \tau_{K_1 + K_3 +
K_4}\approx\tau_{K_1 + K_4},  \ \ \tau_{K_1 + K_3}\approx\tau_{K_1}. 
$$
Consequently we obtain the following proposition.
\begin{prop} Suppose that $\dim\mathfrak{z}=0$. Let $\tau$ be an involution 
of $\mathfrak{g}$ such that $\tau\circ\sigma=\sigma\circ\tau$. Then $\tau$ is 
conjugate within ${\rm Aut}_{\mathfrak h}(\mathfrak{g})$ to one of 
automorphisms listed in Table 3.
\label{prop:71}
\end{prop}
{\small
\begin{center}
\renewcommand{\arraystretch}{1.2}
\setlength{\unitlength}{0.7pt}
\begin{longtable}{l|l|l|l}
\caption{$\dim\mathfrak z=0, \ \tau\circ\sigma=\sigma\circ\tau$, 
$\sigma=\tau_{(1/2)H}$ and $\mathfrak{k}=\mathfrak{g}^{\tau}$.} \\
\hline
\multicolumn{1}{c|}{$(\mathfrak{g}, \mathfrak h, H)$} &
 \multicolumn{1}{|c|}{$h \ (\tau=\tau_h)$}
 &\multicolumn{1}{|c|}{$\mathfrak k$} & 
 \multicolumn{1}{|c}{$\mathfrak {h\cap k}$}\\
\hline
\hline

$(\mathfrak{e}_8,
 \mathfrak{su}(8)\oplus\mathfrak{su}(2), K_3)$ & $K_1$ & $D_8$  &
 $A_7\oplus\mathbb{R}$ \\   

& $K_3$ & $E_7\oplus A_1$ & $A_7\oplus A_1$ \\

& $K_4$ & $E_7\oplus A_1$ & $A_5\oplus A_1\oplus A_1\oplus\mathbb{R}$ \\

& $K_6$ & $D_8$ & $A_3\oplus A_3\oplus A_1\oplus\mathbb{R}$ \\

& $K_3 + K_4$ & $D_8$ & $A_5\oplus A_1\oplus A_1\oplus\mathbb{R}$ \\

& $K_3 + K_6$ & $E_7\oplus A_1$ & $A_3\oplus A_3\oplus A_1\oplus\mathbb{R}$\\

& $K_1 + K_4$ & $E_7\oplus A_1$ & $A_5\oplus A_1\oplus\mathbb{R}$ \\

& $K_1 + K_6$ & $D_8$  & $A_5\oplus A_1\oplus\mathbb{R}$ \\
\hline
$(\mathfrak{e}_8,
 \mathfrak{so}(10)\oplus\mathfrak{so}(6), K_6)$ & $K_1$ & $D_8$ & 
 $D_4\oplus D_3\oplus\mathbb{R}$ \\   

& $K_3$ & $E_7\oplus A_1$ & $D_3\oplus D_3\oplus D_2$ \\

& $K_6$ & $D_8$ & $D_5\oplus D_3$ \\

& $K_8$ & $E_7\oplus A_1$ & $D_5\oplus D_2\oplus\mathbb{R}$ \\

& $K_1 + K_8$ & $E_7\oplus A_1$ & $D_4\oplus D_2\oplus\mathbb{R}^2$ \\

& $K_2 + K_7$ & $E_7\oplus A_1$ & $A_4\oplus A_2\oplus\mathbb{R}^2$ \\

& $K_3 + K_8$ & $D_8$ & $D_3\oplus D_2\oplus D_2\oplus\mathbb{R}$ \\

& $K_2 + K_6 + K_7$ & $D_8$ & $A_4\oplus A_2\oplus\mathbb{R}^2$ \\

\hline
$(\mathfrak{e}_7,
 \mathfrak{so}(6)\oplus\mathfrak{so}(6)\oplus\mathfrak{su}(2), K_4)$ & $K_1$ & 
 $D_6\oplus A_1$ & $D_3\oplus D_2\oplus A_1\oplus\mathbb{R}$\\

& $K_2$ & $A_7$ & $D_3\oplus D_3\oplus\mathbb{R}$ \\

& $K_4$ & $D_6\oplus A_1$ & $D_3\oplus D_3\oplus A_1$ \\

& $K_1+K_2$ & $E_6\oplus\mathbb{R}$ & $D_3\oplus D_2\oplus\mathbb{R}^2$ \\

& $K_1+K_3$ & $D_6\oplus A_1$ & $D_3\oplus D_2\oplus A_1\oplus\mathbb{R}$\\

& $K_1+K_6$ & $D_6\oplus A_1$ & $D_2\oplus D_2\oplus A_1\oplus\mathbb{R}^2$\\

& $K_3+K_4$ & $D_6\oplus A_1$ & $D_3\oplus D_3\oplus A_1$ \\

& $K_3 + K_7$ & $A_7$ & $A_2\oplus A_2\oplus A_1\oplus \mathbb{R}$ \\

& $K_1 + K_2 + K_6$ & $A_7$ & $D_2\oplus D_2\oplus\mathbb{R}^3$ \\

& $K_1 + K_3 + K_4$ & $D_6\oplus A_1$ & $D_3\oplus D_2\oplus
	 A_1\oplus\mathbb{R}$\\

& $K_2 + K_3 + K_7$ & $D_6\oplus A_1$ & $A_2\oplus A_2\oplus
	 \mathbb{R}^3$\\

& $K_3 + K_4 + K_7$ & $E_6\oplus\mathbb{R}$ & $A_2\oplus A_2\oplus
	 A_1\oplus\mathbb{R}^2$\\

\hline
$(\mathfrak{f}_4,
 \mathfrak{so}(6)\oplus\mathfrak{so}(3), K_3)$ & $K_1$ & $B_3\oplus A_1$ &
 $D_2\oplus C_2\oplus\mathbb{R}$ \\ 

& $K_3$ & $C_4$ & $D_3 \oplus C_2$ \\  

& $K_4$ & $C_4$ & $D_3\oplus\mathbb{R}$ \\

& $K_1 + K_4$ & $B_3\oplus A_1$ & $D_3\oplus\mathbb{R}$ \\
\hline
\hline
\multicolumn{1}{c|}{$(\mathfrak{g}, \mathfrak{h}, H)$} &
 \multicolumn{1}{|c|}{$\tau$}
 &\multicolumn{1}{|c|}{$\mathfrak k$} & 
 \multicolumn{1}{|c}{$\mathfrak {h\cap k}$}\\
\hline
\hline
$(\mathfrak{e}_7,
 \mathfrak{so}(6)\oplus\mathfrak{so}(6)\oplus\mathfrak{su}(2), K_4)$
& $\varphi$ & $E_6\oplus \mathbb{R}$ & $D_8\oplus C_1$ \\

& $\varphi\circ\tau_{K_2}$ & $A_7$ & $D_8\oplus D_1$ \\

& $\varphi\circ\tau_{K_4}$ & $A_7$ & $D_8\oplus C_1$ \\

& $\varphi\circ\tau_{K_2}\circ\tau_{K_4}$ & $A_7$ & $D_8\oplus D_1$ \\
\hline

\multicolumn{4}{l}{$\varphi\ : \ E_{\alpha_1}\mapsto E_{\alpha_6}, \
 E_{\alpha_2}\mapsto E_{\alpha_2}, \ E_{\alpha_3}\mapsto
 E_{\alpha_5}, \ E_{\alpha_4}\mapsto E_{\alpha_4}, \ E_{\alpha_7}\mapsto
 E_{\alpha_0}$}

\end{longtable}
\end{center}
}


\section{The case where $\dim\mathfrak z = 1, \tau\circ\sigma 
= \sigma\circ\tau$}
 Let $(G/H,\langle,\rangle,\sigma)$ be a compact Riemannian 4-symmetric space
such that $G$ is simple and $\sigma={\rm Ad}(\exp(\pi/2)\sqrt{-1}K_i)$ for 
some $\alpha_i\in\varPi(\mathfrak{g}_{\mathbb{C}},\mathfrak{t}_{\mathbb{C}})$ 
with $m_i=3$. By Remark 2.2, we have $\dim\mathfrak{z}=1$. 
We shall classify the equivalence classes 
of involutions $\tau$ such that $\tau\circ\sigma = \sigma\circ\tau$. 
According to Section 3 and Jim\'{e}nez \cite {J},
4-symmetric pairs $(\mathfrak {g,h})$ satisfying the condition $\dim
\mathfrak z =1$ are given by
\begin{equation}
  \begin{array}{l}
     (\mathfrak{e} _6, \mathfrak{su}(3)\oplus\mathfrak{su}(3)\oplus
      \mathfrak{su}(2)\oplus\mathbb {R}), \ (\mathfrak{e} _7,
      \mathfrak{su}(5)\oplus\mathfrak{su}(3)\oplus\mathbb {R}), \\ 
     (\mathfrak{e} _7, \mathfrak{su}(6) \oplus \mathfrak{su}(2)
      \oplus\mathbb{R}), (\mathfrak{e} _8, \mathfrak{su}(8)  \oplus
      \mathbb{R}), \  (\mathfrak{e} _8, \mathfrak{su}(2)  \oplus
      \mathfrak{e} _6  \oplus \mathbb{R}), \\ 
     (\mathfrak{f} _4, \mathfrak{su}(3)  \oplus \mathfrak{su}(2)  \oplus
      \mathbb{R}), \ (\mathfrak{g}_2, \mathfrak{su}(2) \oplus
      \mathbb{R}).   
  \end{array}
    \label{eqn:81}
\end{equation}\par
 Suppose that $\mathfrak{g}$ is of type $\mathfrak{e} _8$. From Section
 3, the Dynkin diagram of $\mathfrak{h}$ is one of the 
following:\vspace{0.3cm}
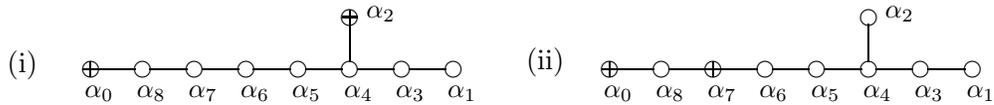
\begin{center}
\setlength{\unitlength}{0.7pt}
\begin{picture}(270,50)
 \put(10,13){(i)}
 \put(55,15){\circle{8}}
 \put(51,15){\line(1,0){8}}
 \put(55,11){\line(0,1){8}}
 \put(52,0){\small $\alpha_0$}
 \put(59,15){\line(1,0){20}}
 \put(83,15){\circle{8}}
 \put(80,0){\small $\alpha_8$}
 \put(87,15){\line(1,0){20}}
 \put(111,15){\circle{8}}
 \put(108,0){\small $\alpha_7$}
 \put(115,15){\line(1,0){20}}
 \put(139,15){\circle{8}}
 \put(136,0){\small $\alpha_6$}
 \put(143,15){\line(1,0){20}}
 \put(167,15){\circle{8}}
 \put(164,0){\small $\alpha_5$}
 \put(171,15){\line(1,0){20}}
 \put(195,15){\circle{8}}
 \put(192,0){\small $\alpha_4$}
 \put(195,19){\line(0,1){20}}
 \put(195,43){\circle{8}}
 \put(191,43){\line(1,0){8}}
 \put(195,39){\line(0,1){8}}
 \put(204,43){\small $\alpha_2$}
 \put(199,15){\line(1,0){20}}
 \put(223,15){\circle{8}}
 \put(220,0){\small $\alpha_3$}
 \put(227,15){\line(1,0){20}}
 \put(251,15){\circle{8}}
 \put(248,0){\small $\alpha_1$}
\end{picture} 
\setlength{\unitlength}{0.7pt}
\begin{picture}(300,50)
 \put(10,15){(ii)}
 \put(55,15){\circle{8}}
 \put(51,15){\line(1,0){8}}
 \put(55,11){\line(0,1){8}}
 \put(52,0){\small $\alpha_0$}
 \put(59,15){\line(1,0){20}}
 \put(83,15){\circle{8}}
 \put(80,0){\small $\alpha_8$}
 \put(87,15){\line(1,0){20}}
 \put(111,15){\circle{8}}
 \put(107,15){\line(1,0){8}}
 \put(111,11){\line(0,1){8}}
 \put(108,0){\small $\alpha_7$}
 \put(115,15){\line(1,0){20}}
 \put(139,15){\circle{8}}
 \put(136,0){\small $\alpha_6$}
 \put(143,15){\line(1,0){20}}
 \put(167,15){\circle{8}}
 \put(164,0){\small $\alpha_5$}
 \put(171,15){\line(1,0){20}}
 \put(195,15){\circle{8}}
 \put(192,0){\small $\alpha_4$}
 \put(195,19){\line(0,1){20}}
 \put(195,43){\circle{8}}
 \put(204,43){\small $\alpha_2$}
 \put(199,15){\line(1,0){20}}
 \put(223,15){\circle{8}}
 \put(220,0){\small $\alpha_3$}
 \put(227,15){\line(1,0){20}}
 \put(251,15){\circle{8}}
 \put(248,0){\small $\alpha_1$}
\end{picture}
\end{center}
\vspace{0.1cm}\par
{\it The case} (i) : In this case, $\sigma={\rm
Ad}(\exp(\pi/2)\sqrt{-1}K_2)$. From Lemma \ref{lem:32}, the
possibilities of positive roots whose 
coefficients of $\alpha_2$ are $3$ are as follows:
\begin{equation}
 \begin{array}{lll}
  \left(
   \begin{array}{lllllll}
       &   &   &   & 3 &   &   \\
     1 & 2 & 3 & 4 & 5 & 3 & 1 \\
   \end{array}\right),& \left(
   \begin{array}{lllllll}
       &   &   &   & 3 &   &   \\
     1 & 2 & 3 & 4 & 5 & 3 & 2 \\
   \end{array}\right),& \left(
   \begin{array}{lllllll}
       &   &   &   & 3 &   &   \\
     1 & 2 & 3 & 4 & 5 & 4 & 2 \\
   \end{array}\right), \\ 
   \left(
   \begin{array}{lllllll}
       &   &   &   & 3 &   &   \\
     1 & 2 & 3 & 4 & 6 & 4 & 2 \\
   \end{array}\right),& \left(
   \begin{array}{lllllll}
       &   &   &   & 3 &   &   \\
     1 & 2 & 3 & 5 & 6 & 4 & 2 \\
   \end{array}\right), & \left(
   \begin{array}{lllllll}
       &   &   &   & 3 &   &   \\
     1 & 2 & 4 & 5 & 6 & 4 & 2 \\
   \end{array}\right),\\
   \left(
   \begin{array}{lllllll}
       &   &   &   & 3 &   &   \\
     1 & 3 & 4 & 5 & 6 & 4 & 2 \\
   \end{array}\right),&  \left(
   \begin{array}{lllllll}
       &   &   &   & 3 &   &   \\
     2 & 3 & 4 & 5 & 6 & 4 & 2 \\
   \end{array}\right).
 \end{array}
   \label{eqn:82}
\end{equation}
Since $\tau(\varPi(\mathfrak{h}))=\varPi(\mathfrak{h})$ and 
$\delta + \alpha_j\not\in\varDelta(\mathfrak{g}_{\mathbb{C}},
\mathfrak{t}_{\mathbb{C}}) \ (j\ne 2)$, we have 
$\tau(\delta) + \alpha_k\notin \varDelta(\mathfrak{g}_{\mathbb{C}},
\mathfrak{t}_{\mathbb{C}})\ (k\ne 2)$. Thus it follows from (\ref{eqn:82}) 
that $\tau(\delta) =\delta$. If $\tau$ satisfies 
$$\tau(\alpha_1) = \alpha_8, \ \tau(\alpha_3) = \alpha_7, \
\tau(\alpha_4) = \alpha_6, \ \tau(\alpha_5) = \alpha_5,$$
we get
$$\left(\begin{array}{lllllll}
  &   &   &   & 3 &   &   \\
2 & 3 & 4 & 5 & 6 & 4 & 2 \\
\end{array}\right)
 = \tau(\delta) = 3\tau(\alpha_2)+
\left(\begin{array}{lllllll}
  &   &   &   & 0 &   &   \\
2 & 4 & 6 & 5 & 4 & 3 & 2 \\
\end{array}\right).$$
Hence we have
$$3\tau(\alpha_2)=
\left(\begin{array}{lllllll}
  &    &    &   & 3 &   &   \\
0 & -1 & -2 & 0 & 2 & 1 & 0 \\
\end{array}\right),$$ 
which is a contradiction. Therefore $\tau$ satisfies 
$\tau|_{\mathfrak{t}} ={\rm Id}_{\mathfrak{t}}$. 
Hence from Proposition 5.3 of Chapter IX of \cite{Helg}, $\tau$ has a 
form $\tau_H$ for a suitable element $H\in\sqrt{-1}\mathfrak{t}$. From
(\ref{eqn:81}), we have 
\begin{equation}
  \mathfrak h \cong A_7\oplus\mathbb{R} \sqrt{-1} K_2\cong \mathfrak
   {su}(8)  \oplus \mathbb {R} \quad \mbox{ and } \quad \mathfrak{t}
   =(A_7\cap\mathfrak{t})\oplus\mathbb{R} \sqrt{-1} K_2,
\end{equation}
and we can write
$$
\tau = \tau_H = \tau_{T+kK_2} = \tau_T\circ\tau_{kK_2}, \ \ 
\sqrt{-1}T\in A_7\cap\mathfrak{t}, \ k\in\mathbb{R}.  
$$
Note that $\tau_T = \tau|_{A_7} : A_7 \to A_7$ and $(\tau_T)^2 ={\rm Id}$ 
on $A_7$.

 We define $v_i \in \sqrt{-1}(A_7\cap\mathfrak{t}), \ i\in
\Lambda:=\{1,3,4,5,6,7,8\}$ by $\alpha_i (v_j) = \delta_{ij}$.  From
Lemma \ref{lem:23}, we may suppose that $\tau$ is
conjugate within ${\rm Int}(\mathfrak h)$ to one of the following 
automorphisms:
\begin{equation}
 \left\{
  \begin{array}{l}
    {\rm Ad}(\exp\pi\sqrt{-1}(2m_1v_1 +2m_3v_3 +\cdots +2m_8v_8+kK_2)),
     \vspace{0.1cm}\\  
    {\rm Ad}(\exp\pi\sqrt{-1}(v_i + 2m_1v_1 +2m_3v_3 +\cdots +2m_8v_8+kK_2)),
  \end{array}\right.
   \label{eqn:84}
\end{equation}
where $i= 1,3,4,5$ and $m_1,m_3,m_4,\ldots , m_8\in\mathbb Z$. 
Put $K_2 = \sum _{i=1}^8 b_i H_{\alpha_i}$. Then we have
$$\delta_{j2} = \alpha_j (K_2) = \sum_{i =1}^8 b_i \alpha_j
(H_{\alpha_i}), \ \mbox{for} \ j=1,2,\ldots ,8,$$
and therefore
$$
\begin{array}{l}
  \displaystyle{b_1 - \frac{b_3}{2} = 0, \ \  b_2 - \frac{b_4}{2}\alpha_2
  (H_{\alpha_2})=1, \ \ -\frac{b_1}{2}+ b_3 - \frac{b_4}{2}=0, \ \
  -\frac{b_2}{2}-\frac{b_3}{2}+b_4-\frac{b_5}{2}=0,}\vspace{0.2cm}\\  
  \displaystyle{-\frac{b_4}{2}+\frac{b_5}{2}-\frac{b_6}{2}=0, \ \
   -\frac{b_5}{2}+b_6-\frac{b_7}{2}=0,\ \ -\frac{b_6}{2}+
   b_7-\frac{b_8}{2}= 0, \ \ -\frac{b_7}{2}+b_8=0.} 
\end{array}$$
Indeed if $j=1$, considering the $\alpha_1$ series containing
$\alpha_i$, we have $\alpha_i(H_{\alpha_1}) = 0$ for $i\ne 1,3$ and
$2\alpha_3(H_{\alpha_1})/\alpha_1 (H_{\alpha_1})=-1$. Thus we
get 
\begin{eqnarray*}
  0 &\!\!\! = &\!\!\! \alpha_1 (K_2) = \sum_{i=1}^8 b_i
  \alpha_1(H_{\alpha_i}) = b_1 \alpha_1(H_{\alpha_1}) + b_3
  \alpha_1(H_{\alpha_3}) \\
  &\!\!\! = &\!\!\! b_1 \alpha_1(H_{\alpha_1}) + b_3 (-\frac{1}{2}\alpha_1(H_{\alpha_1}))
  = (b_1-\frac{b_3}{2})\alpha_1(H_{\alpha_1}).
\end{eqnarray*}
We can obtain the other equations by a similar computation as above. 
\vspace{0.2cm}

Computing these simultaneous equations we have 
\begin{equation}
  K_2 = \frac{c_8}{3}(5H_{\alpha_1} + 8H_{\alpha_2} + 10H_{\alpha_3} +
  15H_{\alpha_4} + 12H_{\alpha_5} + 9H_{\alpha_6} + 6H_{\alpha_7} +
  3H_{\alpha_8}).
\label{eqn:85}
\end{equation}
Now put $v_1 = a_1 K_1 + \cdots +a_8 K_8, \ a_1, \ldots , a_8 \in \mathbb
R$. Then we get $a_3 = \cdots =a_8 = 0$ and $a_1 =1$, since $\alpha _i
(v_1) = \delta_{i1}$ for $i\in\Lambda$. Thus we have $v_1
= K_1 + a_2 K_2$. Since $v_1\bot K_2$, it follows from (\ref{eqn:85}) that 
$0 = (5/3)c_8 + (8/3)a_2 c_8$
and therefore $a_2=-5/8$. Hence we have $v_1 = K_1 -
(5/8)K_2$. By a similar computation, we obtain 
\begin{equation}
 \begin{array}{l}
  \displaystyle v_1 = K_1 - \frac{5}{8}K_2, \ v_3 =  - \frac{5}{4}K_2 +
   K_3, \ v_4 = - \frac{15}{8}K_2 + K_4, \ v_5 =  - \frac{3}{2}K_2 + K_5, 
\\  
  \displaystyle v_6 =  - \frac{9}{8}K_2 + K_6, \ v_7 =  - \frac{3}{4}K_2 + K_7,
 \ v_8 =  - \frac{3}{8}K_2 + K_8.
 \end{array}
   \label{eqn:86}
\end{equation}
Therefore by (\ref{eqn:84}) and (\ref{eqn:86}) it follows that $\tau$ is 
conjugate within ${\rm Int}(\mathfrak h)$ to one of the following:
\begin{equation}
 \tau_{mK_2},\ \ \tau_{v_i + mK_2},
\end{equation}
where $m= - (1/4)(5m_1 + 10m_3 + 15m_4 + 12m_5 + 9m_6 + 6m_7 +
3m_8 - 2k)$. 
Moreover, since $\tau^2={\rm Id}$, it ifollows from (\ref{eqn:86}) that 
$\tau$ is conjugate within ${\rm Int}(\mathfrak h)$ to some $\tau_h$ where 
$h$ is one of the following:
$$
K_i, \ K_j + K_2,\ \ i=1,2,3,4,5, \ j=1,3,4,5.
$$
By a computation similar to Section 7 we can obtain $(\mathfrak{g}^{\tau_h},
\mathfrak{h}\cap\mathfrak{g}^{\tau_h})$ for each $h$. 
\vspace{0.5cm}\par

{\sc Remark 8.1.} From Lemma \ref{lem:23}, we can see that
$\tau_{v_8}|_{A_7}$ is conjugate 
within ${\rm Int}(A_7)(\subset{\rm Int}(\mathfrak h))$ to $\tau_{v_1}|_{A_7}$. 
Therefore by the above argument, $\tau_{K_8}$ is conjugate within 
${\rm Int}(\mathfrak h)$ to $\tau_{K_1}$ or $\tau_{K_1+K_2}$. However, 
$\mathfrak{g}^{\tau_{K_8}}\not\cong\mathfrak{g}^{\tau_{K_1}}$, and hence 
$\tau_{K_8}\approx\tau_{K_1+K_2}$. 
\vspace{0.5cm}\par

{\it The case} (ii): In this case, $\sigma={\rm Ad}(\exp(\pi/2)\sqrt{-1}K_7)$ 
and 
\begin{eqnarray*}
& & \mathfrak h \cong A_1\oplus E _6\oplus\mathbb{R} \sqrt{-1} K_7
\cong\mathfrak{su}(2)\oplus\mathfrak{e}_6\oplus{\mathbb{R}}, \\
& & \mathfrak{t} =(A_1\cap\mathfrak{t})\oplus(E_6\cap\mathfrak
t)\oplus\mathbb{R} \sqrt{-1}K_7.
\end{eqnarray*}
By a computation similar to the case (i), we have
$\tau|_{\mathfrak{t}} = {\rm Id}_{\mathfrak{t}}$. 
Hence we can write
$$\tau = \tau_{T_1}\circ\tau_{T_2}\circ\tau_{kK_7},$$    
where $\sqrt{-1}T_1\in A_1\cap\mathfrak{t}, \ \sqrt{-1}T_2\in 
E _6\cap\mathfrak{t}, \  k\in\mathbb{R}$. 
We define $v_8\in\sqrt{-1}(A_1\cap\mathfrak{t})$ and
$v_a\in\sqrt{-1}(E_6\cap\mathfrak{t}), \ a\in\Lambda :=\{1,2,3,4,5,6\}$
by $\alpha_i (v_j) = \delta_{ij}$. Then from Lemma \ref{lem:23}, we may 
assume $\tau$ is conjugate within ${\rm Int}(\mathfrak h)$ to one of following 
automorphisms:
\begin{equation}
 \left\{
  \begin{array}{l}
   {\rm Ad}(\exp\pi\sqrt{-1}(2m_1v_1 +\cdots +2m_6v_6 + 2m_8v_8+kK_7)),
    \vspace{0.1cm}\\  
   {\rm Ad}(\exp\pi\sqrt{-1}(v_a + 2m_1v_1  +\cdots +2m_6v_6
    +2m_8v_8+kK_7)),
    \vspace{0.1cm}\\  
   {\rm Ad}(\exp\pi\sqrt{-1}(v_8 + v_b + 2m_1v_1  +\cdots +2m_6v_6
    +2m_8v_8+kK_7)),
  \end{array}\right.
\end{equation}
where $a= 1,2,8, \ b = 1,2$ and $m_1,\ldots ,m_6,m_8\in\mathbb Z$.
Furthermore we obtain 
$$
\begin{array}{l}
  \displaystyle v_1 = K_1 - \frac{2}{3}K_7, \ v_2 =  K_2 - K_7, \ v_3 =
  K_3 -\frac{4}{3} K_7, \ v_4 =   K_4 - 2K_7, \vspace{0.2cm}\\ 
  \displaystyle v_5 = K_6 - \frac{5}{3}K_7, \ 
  v_6 = K_6 - \frac{4}{3}K_7, \ v_8=K_8-\frac{1}{2}K_7.
\end{array}$$
Similarly as in the case (i), we can see that $\tau$ is conjugate within 
${\rm Int}(\mathfrak{h})$ to some $\tau_h$ where $h$ is one of the following:
$$
K_i, \ K_j + K_7, \ K_k + K_8, \ K_k + K_7 + K_8, 
$$
where $i=1,2,7,8$, $j=1,2,8$ and $k=1,2$. 
It is easy to check that the reflection 
$t_{\alpha_8}\in{\rm Int}(\mathfrak h)$ maps
 $K_7 + K_8\mapsto 2K_7-K_8, \ K_1 + K_7 + K_8\mapsto K_1+2K_7
 - K_8, \ K_2 + K_7 + K_8\mapsto K_2 + 2K_7 -K_8$ and $K_7\mapsto
 K_7$. Therfore we have $\tau_{K_7 + K_8}\approx\tau_{K_8}, \ \tau_{K_1
 + K_7 + K_8}\approx\tau_{K_1 + K_8}$ and
 $\tau_{K_2 + K_7 + K_8}\approx\tau_{K_2 + K_8}$.\vspace{0.3cm}\par

In the case where $\mathfrak{g}=\mathfrak{e}_7$, the Dynkin diagram of 
$\mathfrak h$ is one of the following:

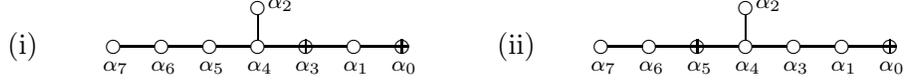
\begin{center}
(i) \ \ \ \
\setlength{\unitlength}{1.3pt}
\begin{picture}(92,20)
 \put(5,2){\circle{4}}
 \put(2,-5){{\scriptsize $\alpha_7$}}
 \put(7,2){\line(1,0){10}}
 \put(19,2){\circle{4}}
 \put(16,-5){{\scriptsize $\alpha_6$}}
 \put(21,2){\line(1,0){10}}
 \put(33,2){\circle{4}}
 \put(30,-5){{\scriptsize $\alpha_5$}}
 \put(35,2){\line(1,0){10}}
 \put(47,2){\circle{4}}
 \put(44,-5){{\scriptsize $\alpha_4$}}
 \put(47,13){\circle{4}}
 \put(50,13){{\scriptsize $\alpha_2$}}
 \put(47,4){\line(0,1){7}}
 \put(49,2){\line(1,0){10}}
 \put(61,2){\circle{4}}
 \put(58,-5){{\scriptsize $\alpha_3$}}
 \put(59,2){\line(1,0){4}}
 \put(61,0){\line(0,1){4}}
 \put(63,2){\line(1,0){10}}
 \put(75,2){\circle{4}}
 \put(72,-5){{\scriptsize $\alpha_1$}}
 \put(77,2){\line(1,0){10}}
 \put(89,2){\circle{4}}
 \put(86,-5){{\scriptsize $\alpha_0$}}
 \put(87,2){\line(1,0){4}}
 \put(89,0){\line(0,1){4}}
\end{picture}\hspace{1cm}
(ii) \ \ \ \
\begin{picture}(92,20)
 \put(5,2){\circle{4}}
 \put(2,-5){{\scriptsize $\alpha_7$}}
 \put(7,2){\line(1,0){10}}
 \put(19,2){\circle{4}}
 \put(16,-5){{\scriptsize $\alpha_6$}}
 \put(21,2){\line(1,0){10}}
 \put(33,2){\circle{4}}
 \put(30,-5){{\scriptsize $\alpha_5$}}
 \put(31,2){\line(1,0){4}}
 \put(33,0){\line(0,1){4}}
 \put(35,2){\line(1,0){10}}
 \put(47,2){\circle{4}}
 \put(44,-5){{\scriptsize $\alpha_4$}}
 \put(47,13){\circle{4}}
 \put(50,13){{\scriptsize $\alpha_2$}}
 \put(47,4){\line(0,1){7}}
 \put(49,2){\line(1,0){10}}
 \put(61,2){\circle{4}}
 \put(58,-5){{\scriptsize $\alpha_3$}}
 \put(63,2){\line(1,0){10}}
 \put(75,2){\circle{4}}
 \put(72,-5){{\scriptsize $\alpha_1$}}
 \put(77,2){\line(1,0){10}}
 \put(89,2){\circle{4}}
 \put(86,-5){{\scriptsize $\alpha_0$}}
 \put(87,2){\line(1,0){4}}
 \put(89,0){\line(0,1){4}}
\end{picture}
\end{center}
${}$\par
 {\it The case} (i): In this case, $\sigma={\rm Ad}(\exp(\pi/2)\sqrt{-1}K_3)$ 
and $\mathfrak{h}\cong A_1\oplus A_5\oplus{\mathbb{R}}\sqrt{-1}K_3$. By an 
argument similar to the above, we can see that $\tau$ is conjugate within 
${\rm Int}(\mathfrak h)$ to some $\tau_h$ where $h$ is one of the following:  
$$K_i, \ K_j + K_3, \ K_k + K_1, \ K_k + K_1 + K_3,$$
where $i=1,2,3,4,5$, $j=1,2,4,5$ and $k=2,4,5$. 
Using the reflection $t_{\alpha_1}\in{\rm Int}(\mathfrak h)$ we obtain
$$
\tau_{K_1 + K_3}\approx\tau_{K_1}, \ \ \tau_{K_1 + K_2 +
K_3}\approx\tau_{K_1 + K_2},  \ \ \tau_{K_1 + K_3 + K_4}\approx\tau_{K_2 +
K_4}, \ \ \tau_{K_1 + K_3 + K_5}\approx\tau_{K_1 + K_5}. 
$$
Furthermore, similarly as in Remark 8.1 we get $\tau_{K_2 + K_3}
 \approx\tau_{K_7}$.\vspace{0.5cm}\par 
{\it The case} (ii): In this case, $\sigma={\rm Ad}(\exp(\pi/2)\sqrt{-1}K_5)$ 
and $\mathfrak{h}\cong A_2\oplus A_4\oplus{\mathbb{R}}\sqrt{-1}K_5$. Moreover, 
$\tau$ is conjugate within ${\rm Int}(\mathfrak
h)$ to some $\tau_h$ where $h$ is one of the following: 
$$
K_i, \ K_j + K_5, \  K_k + K_6, \ K_k + K_5 + K_6,
$$ 
where $i=1,3,5,6$, $j=1,3,6$ and $k=1,3$. 
Similarly as in Remark 8.1 we get $\tau_{K_5 + K_6}
\approx\tau_{K_7}$, $\tau_{K_1 + K_5 + K_6} \approx\tau_{K_1 + K_7}$ and
$\tau _{K_3 + K_5 + K_6} \approx \tau_{K_3 + K_7}$.\vspace{0.2cm}\par 
If $\mathfrak{g}=\mathfrak{f}_4$, then 
$\sigma={\rm Ad}(\exp(\pi/2)\sqrt{-1}K_2)$ and $\mathfrak{h}\cong 
A_1\oplus A_2\oplus{\mathbb{R}}\sqrt{-1}K_2$. In this case, 
$\tau$ is conjugate within ${\rm Int}(\mathfrak{h})$ to some 
$\tau_h$ where $h$ is one of the following: 
$$
K_i,\  K_j + K_k, \ K_1 + K_2 + K_3,$$
where $i=1,2,3$ and $j,k=1,2,3$ ($j\not=k$).
Using the reflection $t_{\alpha_1}\in{\rm Int}(\mathfrak{h})$ we have
$\tau_{K_1 + K_2}\approx\tau_{K_1}$ and
$\tau_{K_1 + K_2 + K_3}\approx\tau_{K_1 + K_3}$.\vspace{0.2cm}\par

If $\mathfrak{g}=\mathfrak{g}_2$, then 
$\sigma={\rm Ad}(\exp(\pi/2)\sqrt{-1}K_1)$ and $\mathfrak{h}\cong 
A_1\oplus{\mathbb{R}}\sqrt{-1}K_1$. In this case, 
$\tau$ is conjugate within ${\rm Int}(\mathfrak
 h)$ to some $\tau_h$ where $h$ is one of 
$$
K_i,\ K_1 + K_2,\ \ i=1,2.
$$
Using the reflection $t_{\alpha_2}\in{\rm Int}(\mathfrak{h})$ we have
$\tau_{K_1 + K_2}\approx \tau_{K_2}$. \vspace{0.2cm}\par

If $\mathfrak{g} =\mathfrak{e} _6$, then 
$\sigma={\rm Ad}(\exp(\pi/2)\sqrt{-1}K_4)$ and $\mathfrak{h}\cong 
A_1\oplus A_2\oplus A_2\oplus{\mathbb{R}}\sqrt{-1}K_4$. By an argument 
similar to the case where $\mathfrak{g}=\mathfrak{e}_8$, we obtain
 $\tau|_{\mathfrak{t}}={\rm Id}_{\mathfrak{t}}$ or 
$$\tau(\alpha_1)=\alpha_6, \ \tau(\alpha_3)=\alpha_5, \
\tau(\alpha_4)=\alpha_4, \ \tau(\alpha_2)=\alpha_2.$$
If $\tau|_{\mathfrak{t}} = {\rm Id}_{\mathfrak{t}}$, then $\tau$ is conjugate 
within ${\rm Int}(\mathfrak h)$ to some $\tau_h$ 
where $h$ is one of the following: 
$$
K_i,\ K_j + K_k,\ K_l + K_m + K_n, \ K_1 + K_2 + K_4 + K_5,
$$
where $i=1,2,4,5$, $j,k=1,2,4,5$ ($j<k$) and $l,m,n=1,2,4,5$ 
($l<m<n$).  
Using the reflection $t_{\alpha_2}\in{\rm Int}(\mathfrak h)$ we have 
$$
\tau_{K_2 + K_4}\approx\tau_{K_2}, \ \ \tau_{K_1 + K_2 +
K_4}\approx\tau_{K_1 + K_2},  \ \ \tau_{K_2 + K_4 + K_5}\approx\tau_{K_2 +
K_5}, \ \ \tau_{K_1 + K_2 + K_4 + K_5}\approx\tau_{K_1 + K_2 + K_5}.
$$\vspace{0.2cm}\par 
Next suppose that $\mathfrak{g}=\mathfrak{e}_6$ and $\tau$ satisfies  
$$\tau(\alpha_1)=\alpha_6, \ \tau(\alpha_3)=\alpha_5, \
\tau(\alpha_4)=\alpha_4, \ \tau(\alpha_2)=\alpha_2.$$
Let $\mathfrak{t}_{\pm}$ be the $(\pm 1)$-eigenspaces of $\tau|_{\mathfrak{t}}$, 
respectively. Then we have
$$\mathfrak{t}_+ = {\rm span}\{K_1 + K_6, \ K_2, \  K_3 + K_5, \ K_4\}, \ \
\mathfrak{t}_- = {\rm span}\{K_1 - K_6 ,  \ K_3 - K_5 \}$$
It is known that there exists an involutive automorphism
$\psi$ of outer type satisfying 
\begin{equation}
   \psi(E_{\alpha_1})=E_{\alpha_6}, \ \psi(E_{\alpha_2})=E_{\alpha_2},
   \ \psi(E_{\alpha_3})=E_{\alpha_5}, \ \psi(E_{\alpha_4})=E_{\alpha_4}.
     \label{eqn:88}
\end{equation}
Therefore there exists $\sqrt{-1}h_+\in\mathfrak{t}_+$ such that 
$\tau_{h_+}^2={\rm Id}$ and $\tau\approx\psi\circ\tau_{h_+}$. Then 
by an argument similar to that in Section 7, we can see that $\tau$ is 
conjugate within ${\rm Aut}_{\mathfrak h}(\mathfrak{g})$ to one of the 
following involutions: 
$$
\psi, \ \ \psi\circ\tau_{K_2},\ \ \psi\circ\tau_{K_4},\ \ 
\psi\circ\tau_{K_2+K_4}.  
$$

Since $\mathfrak{su}_{\alpha_2}(2)\subset\mathfrak{g}^{\psi}$ and 
$t_{\alpha_2}(K_2)=-K_2+K_4$, we obtain $\psi\circ t_{\alpha_2}=
t_{\alpha_2}\circ\psi$ and 
$$
\psi\circ\tau_{K_2 + K_4}=\psi\circ\tau_{t_{\alpha_2}(K_2)}
=\psi\circ t_{\alpha_2}\circ\tau_{K_2}\circ t_{\alpha_2}^{-1} =
 t_{\alpha_2}(\psi\circ\tau_{K_2}) t_{\alpha_2}^{-1}. 
$$
Thus we obtain
$\psi\circ\tau_{K_2}\approx\psi\circ\tau_{K_2+K_4}$.
Furthermore by an argument similar to that in Section 7, we can compute
$\dim(\mathfrak{h\cap g}^{\tau_h})$ and $\dim\mathfrak{g}^{\tau_h}$ 
for each $h$. Consequently we have the following proposition. 
\begin{prop} 
Suppose that $\dim\mathfrak{z}=1$ and 
$\sigma={\rm Ad}(\exp(\pi/2)\sqrt{-1}K_i)$ for some $\alpha_i\in
\varPi(\mathfrak{g}_{\mathbb{C}},\mathfrak{t}_{\mathbb{C}})$ with $m_i=3$. 
Let $\tau$ be an involution of $\mathfrak{g}$ such that 
$\tau\circ\sigma=\sigma\circ\tau$. Then $\tau$ is conjugate within 
${\rm Aut}_{\mathfrak h}(\mathfrak{g})$ to one of automorphisms listed in 
Table 4. 
\label{prop:81}
\end{prop}

{\small
\begin{center}
\renewcommand{\arraystretch}{1.2}
\begin{longtable}{l|l|l|l}
\caption{$\dim\mathfrak z=1, \ \tau\circ\sigma=\sigma\circ\tau$, 
$\sigma=\tau_{(1/2)H}$ and $\mathfrak{k}=\mathfrak{g}^{\tau}$}\\
\hline
\multicolumn{1}{c|}{$(\mathfrak{g}, \mathfrak{h}, H)$} &
 \multicolumn{1}{|c|}{$h \ (\tau=\tau_h)$}
 &\multicolumn{1}{|c|}{$\mathfrak k$} & 
 \multicolumn{1}{|c}{$\mathfrak {h\cap k}$}\\
\hline
\hline

$(\mathfrak{e}_8,
 \mathfrak{su}(8)\oplus\mathbb{R}, K_2)$ & $K_1$ & $D_8$ &
 $A_6\oplus\mathbb{R}^2 $ \\  

& $K_2$ & $D_8$ & $A_7\oplus\mathbb{R} $ \\ 

& $K_3$ & $E_7\oplus A_1$ &
	     $A_5\oplus A_1\oplus\mathbb{R}^2$ \\   

& $K_4$ & $E_7\oplus A_1$ & $A_4\oplus A_2\oplus\mathbb{R}^2$ \\ 

& $K_5$ & $D_8$ & $A_3\oplus A_3\oplus\mathbb{R}^2$ \\

& $K_8$ & $E_7\oplus A_1$ & $A_6\oplus\mathbb{R}^2$ \\

& $K_2 + K_3$ & $E_7\oplus A_1$ & $A_5\oplus A_1\oplus\mathbb{R}^2$ \\

& $K_2 + K_4$ & $D_8$ & $A_4\oplus A_2\oplus\mathbb{R}^2$ \\ 

& $K_2 + K_5$ & $D_8$ & $A_3\oplus A_3\oplus\mathbb{R}^2$ \\ 
\hline
$(\mathfrak{e}_8,
 \mathfrak{e}_6\oplus\mathfrak{su}(2)\oplus\mathbb{R}, K_7)$ & $K_1$ & $D_8$ &
 $D_5\oplus A_1\oplus\mathbb{R}^2 $ \\

& $K_2$ & $D_8$ & $A_5\oplus A_1\oplus A_1\oplus\mathbb{R}$ \\  

& $K_7$ & $E_7\oplus A_1$ & $E_6\oplus A_1\oplus\mathbb{R} $ \\

& $K_8$ & $E_7\oplus A_1$ & $E_6\oplus\mathbb{R}^2$ \\ 

& $K_1 + K_7$ & $E_7\oplus A_1$ & $D_5\oplus A_1\oplus\mathbb{R}^2$ \\

& $K_2 + K_7$ & $E_7\oplus A_1$ & $A_5\oplus A_1\oplus A_1\oplus\mathbb{R}$ \\ 

& $K_1 + K_8$ & $E_7\oplus A_1$ & $D_5\oplus\mathbb{R}^3$ \\

& $K_2 + K_8$ & $D_8$ & $A_5\oplus A_1\oplus\mathbb{R}^2$ \\

\hline
$(\mathfrak{e} _7,
\mathfrak{su}(6)\oplus\mathfrak{su}(2)\oplus\mathbb{R}, K_3)$ & $K_1$ &
 $D_6\oplus A_1$ & $A_5\oplus\mathbb{R}^2$ \\   

& $K_2$ & $A_7$ & $A_4\oplus A_1\oplus\mathbb{R}^2$ \\ 

& $K_3$ & $D_6\oplus A_1$ & $A_5\oplus A_1\oplus\mathbb{R} $ \\

& $K_4$ & $D_6\oplus A_1$ & $A_3\oplus A_1\oplus A_1\oplus\mathbb{R}^2$ \\ 

& $K_5$ & $A_7$ & $A_2\oplus A_2\oplus\mathbb{R}^2$ \\

& $K_7$ & $E_6\oplus\mathbb{R}$ & $A_4\oplus A_1\oplus\mathbb{R}^2$ \\

& $K_1 + K_2$ & $E_6\oplus\mathbb{R}$ & $A_4\oplus\mathbb{R}^3$ \\  

& $K_1 + K_4$ & $D_6\oplus A_1$ & $A_3\oplus A_1\oplus\mathbb{R}^3$ \\

& $K_1 + K_5$ & $A_7$ & $A_2\oplus A_2\oplus\mathbb{R}^3$ \\ 

& $K_3 + K_4$ & $D_6\oplus A_1$ & $A_3\oplus A_1\oplus A_1\oplus\mathbb{R}^2$ \\ 

& $K_3 + K_5$ & $E_6\oplus\mathbb{R}$ & $A_2\oplus A_2\oplus
	 A_1\oplus\mathbb{R}^2$ \\
\hline
$(\mathfrak{e}_7,
 \mathfrak{su}(5)\oplus\mathfrak{su}(3)\oplus\mathbb{R}, K_5)$ & $K_1$ & 
$D_6\oplus A_1$ & $A_3\oplus A_2\oplus\mathbb{R}^2$ \\

& $K_3$ & $D_6\oplus A_1$ & $A_2\oplus A_2\oplus A_1\oplus\mathbb{R}^2$ \\

& $K_5$ & $A_7$ & $A_4\oplus A_2\oplus\mathbb{R}$ \\

& $K_6$ & $D_6\oplus A_1$ & $A_4\oplus A_1\oplus\mathbb{R}^2$ \\

& $K_7$ & $E_6\oplus\mathbb{R}$ & $A_4\oplus A_1\oplus\mathbb{R}^2$ \\

& $K_1 + K_5$ & $A_7$ & $A_3\oplus A_2\oplus\mathbb{R}^2$ \\

& $K_1 + K_6$ & $D_6\oplus A_1$ & $A_3\oplus A_1\oplus\mathbb{R}^3$ \\

& $K_3 + K_5$ & $E_6 \oplus\mathbb{R}$ & 
$A_2\oplus A_2\oplus A_1\oplus\mathbb{R}^2$ \\ 

& $K_3 + K_6$ & $D_6\oplus A_1$ & $A_2\oplus A_1\oplus A_1\oplus\mathbb{R}^3$ \\

& $K_1 + K_7$ & $E_6\oplus\mathbb{R}$ & $A_3\oplus A_1\oplus\mathbb{R}^3$ \\

& $K_3 + K_7$ & $A_7$ & $A_2\oplus A_1\oplus A_1\oplus\mathbb{R}^3$ \\

\hline

$(\mathfrak{e}_6,
\mathfrak{su}(3)\oplus\mathfrak{su}(3)\oplus\mathfrak{su}(2)\oplus\mathbb{R},
 K_4)$ & 
$K_1$ & $D_5\oplus\mathbb{R}$ & $A_2\oplus A_1\oplus A_1\oplus\mathbb{R}^2$ \\ 

& $K_4$ & $A_5\oplus A_1$ & $A_2\oplus A_2\oplus A_1\oplus\mathbb{R}$ \\

& $K_5$ & $A_5\oplus A_1$ & $A_2\oplus A_1\oplus A_1\oplus\mathbb{R}^2$ \\

& $K_1 + K_2$ & $D_5\oplus\mathbb{R}$ & $A_2\oplus
 A_1\oplus\mathbb{R}^3$ \\

& $K_1 + K_4$ & $A_5\oplus A_1$ & $A_2\oplus A_1\oplus A_1\oplus\mathbb{R}^2$ \\

& $K_1 + K_5$ & $A_5\oplus A_1$ & $A_1\oplus A_1\oplus A_1\oplus\mathbb{R}^3$ \\

& $K_2 + K_4$ & $A_5\oplus A_1$ & $A_2\oplus A_2\oplus\mathbb{R}^2$ \\

& $K_2 + K_5$ & $D_5\oplus\mathbb{R}$ & $A_2\oplus A_1\oplus\mathbb{R}^3$ \\

& $K_4 + K_5$ & $D_5\oplus\mathbb{R}$ & $A_2\oplus A_1\oplus
 A_1\oplus\mathbb{R}^2$ \\ 

& $K_1 + K_2 + K_5$ & $A_5\oplus A_1$ & $A_1\oplus A_1\oplus\mathbb{R}^4$ \\

& $K_1 + K_4 + K_5$ & $D_5\oplus\mathbb{R}$  & $A_1\oplus A_1\oplus
	 A_1\oplus\mathbb{R}^3$ \\ 
\hline

$(\mathfrak{f}_4,
 \mathfrak{su}(3)\oplus\mathfrak{su}(2)\oplus\mathbb{R}, K_2)$ & $K_1$ & 
 $C_3\oplus A_1$ & $A_2\oplus\mathbb{R}^2 $ \\ 
	    
& $K_2$ & $C_3\oplus A_1$ & $A_2\oplus A_1\oplus\mathbb{R}$ \\

& $K_4$ & $B_4$ & $A_1\oplus A_1\oplus\mathbb{R}^2$ \\  
	    
& $K_1 + K_3$ & $C_3\oplus A_1$ & $A_1\oplus\mathbb{R}^3$ \\  

& $K_2 + K_4$ & $C_3\oplus A_1$ & $A_1\oplus A_1\oplus\mathbb{R}^2$ \\ 
\hline
$(\mathfrak{g}_2,
 \mathfrak{su}(2)\oplus\mathbb{R}, K_1)$ & $K_1$ & $A_1\oplus A_1$ & 
 $A_1\oplus\mathbb{R}$ \\  

& $K_2$ & $A_1\oplus A_1$ & $\mathbb{R}^2$ \\
\hline 
\hline
\multicolumn{1}{c|}{$(\mathfrak{g}, \mathfrak{h}, H)$} &
 \multicolumn{1}{|c|}{$\tau$}
 &\multicolumn{1}{|c|}{$\mathfrak k$} & 
 \multicolumn{1}{|c}{$\mathfrak {h\cap k}$}\\
\hline
\hline
$(\mathfrak{e} _6,\mathfrak{su}(3) \oplus\mathfrak{su}(3)
 \oplus\mathfrak{su}(2)\oplus\mathbb{R}, K_4)$ & $\psi$ & $F _4$ &
 $A_2\oplus A_1\oplus C_1\oplus\mathbb{R}$ \\ 

& $\psi\circ\tau_{K_2}$ & $C_4$ & $A_2\oplus A_1\oplus D_1\oplus\mathbb{R}$ \\

& $\psi\circ\tau_{K_4}$ & $C_4$ & $A_2\oplus A_1\oplus C_1\oplus\mathbb{R}$ \\ 

\hline
\multicolumn{4}{l}{$\psi\ : \ E_{\alpha_1}\mapsto E_{\alpha_6}, \
 E_{\alpha_2}\mapsto E_{\alpha_2}, \ E_{\alpha_3}\mapsto E_{\alpha_5}, \
 E_{\alpha_4}\mapsto E_{\alpha_4}$}
\end{longtable}
\end{center}
}


\section{Remarks on conjugations}

\subsection{$\dim\mathfrak z = 1, \ \tau\circ\sigma = \sigma\circ\tau$}
 {\it The case where $\mathfrak{g}=\mathfrak{e}_6$ and 
$\sigma=\tau_{(1/2)K_4}$.} 
We shall show that $\tau_{K_2+K_5}\approx\tau_{K_1+K_2}, \tau_{K_1+K_4}
\approx\tau_{K_5}$ and $\tau_{K_4+K_5}\approx\tau_{K_1}$. We consider $\psi\in
  {\rm Aut}_{\mathfrak h}(\mathfrak{g})$ (see (\ref{eqn:88})).
For $\mu_1 :=t_{\alpha_5 + \alpha_6}\circ t_{\alpha_2}\circ\psi\in
{\rm Aut}(\mathfrak{g})$, we have 
$$
\begin{array}{l}
\mu_1(\alpha_1)=-\alpha_5, \ \mu_1(\alpha_2)=-\alpha_2, \
\mu_1(\alpha_3)=-\alpha_6, \\
\mu_1(\alpha_4)=\alpha_2 + \alpha_4 + \alpha_5 + \alpha_6, \
\mu_1(\alpha_5)=\alpha_3, \ \mu_1(\alpha_6)=\alpha_1.
\end{array}$$
Hence we get $\mu_1^{-1}(K_4)=K_4$ and $\mu_1^{-1}(K_2 + K_5)=K_1 - K_2 +
2K_4$. Thus $\mu_1^{-1}$ is in ${\rm Aut}_{\mathfrak{h}}(\mathfrak{g})$
and gives a conjugation between $\tau_{K_2 + K_5}$ and $\tau_{K_1 +
K_2}$.\par
 Similarly as above, by using $\mu_2 :=t_{\alpha_3} \circ t_{\alpha_1 
+\alpha_3}\circ\psi$, $\mu_3 :=t_{\alpha_1} \circ t_{\alpha_3} \circ\psi 
\in{\rm Aut}_{\mathfrak h}(\mathfrak{g})$, we obtain 
$\tau_{K_1 + K_4}\approx\tau_{K_5}$ and $\tau_{K_4 + K_5}\approx\tau_{K_1}$. 
\vspace{0.3cm}
\par
{\it The case where $\mathfrak{g=e}_8$ and $\sigma=\tau_{(1/2)K_2}$.} From Proposition \ref{prop:81}, we see that $\mathfrak{g}^{\tau_{K_3}}
\cong\mathfrak{g}^{\tau_{K_2+K_3}}$ and $\mathfrak{h\cap g}^{\tau_{K_3}}
\cong\mathfrak{h\cap g}^{\tau_{K_2+K_3}}$. 
Now, we shall show that $\tau_{K_3}$ is  not conjugate within 
${\rm Aut}_{\mathfrak h} (\mathfrak{g})$ to
 $\tau_{K_2+K_3}$. Put $\mathfrak{k}_1 := \mathfrak{g}^{\tau_{K_3}}$ and
 $\mathfrak{k}_2 :=\mathfrak{g}^{\tau_{K_2 + K_3}}$,
  then we have $\mathfrak{k}_1\cong\mathfrak{k}_2\cong 
 A_1\oplus E_7$ and $\mathfrak{k}_1 \cap\mathfrak{h} \cong\mathfrak{k}_2
 \cap\mathfrak{h}$. We denote $\alpha\in\varDelta (\mathfrak{g}_{\mathbb
 C}, \mathfrak{t}_{\mathbb{C}})$ by $\alpha = \sum _{i=1}^{8}n_i
 \alpha_i$ and put
\begin{eqnarray*}
& & \varDelta_{\mathfrak{k}_1}:=\{\alpha\in
\varDelta^+(\mathfrak{g}_{\mathbb{C}},\mathfrak{t}_{\mathbb{C}}) \ ; \ 
\alpha(K_3)=0,2,4 \}, \\ 
& & \varDelta_{\mathfrak{k}_2}:=\{\alpha\in
\varDelta^+(\mathfrak{g}_{\mathbb{C}},\mathfrak{t}_{\mathbb{C}}) \ ; \  
\alpha(K_3)=0,2,4,6 \}. 
\end{eqnarray*}
Then 
$$\mathfrak{k}_i =\mathfrak{t} \oplus
 \sum_{\alpha\in\varDelta_{\mathfrak k_i}} (\mathbb{R}
 A_{\alpha} + \mathbb{R} B_{\alpha}),\ \ i=1,2.
$$
Put $\gamma:=2\alpha_1+2\alpha_2+4\alpha_3+5\alpha_4+4\alpha_5+3\alpha_6+
2\alpha_7+\alpha_8$. Then for any $\alpha\in\varDelta_{\mathfrak k _1}$ and 
$\alpha'\in\varDelta_{\mathfrak k _2}$, we can see that $\alpha_1\pm\alpha
\not\in\varDelta(\mathfrak{g}_{\mathbb{C}},\mathfrak{t}_{\mathbb{C}})$ and 
$\gamma\pm\alpha'\not\in\varDelta(\mathfrak{g}_{\mathbb{C}},
\mathfrak{t}_{\mathbb{C}})$ (cf. \cite{FV}). Therefore, we get
\begin{equation}
\mathfrak{k}_1=\mathfrak{su}_{\alpha_1}(2)\oplus
\mathfrak{su}_{\alpha_1}(2)^{\perp}, \ \ 
\mathfrak{k}_2=\mathfrak{su}_{\gamma}(2)\oplus
\mathfrak{su}_{\gamma}(2)^{\perp},
\label{eqn9-1}
\end{equation}
where $\mathfrak{su}_{\alpha_1}(2)^{\perp}\cong 
\mathfrak{su}_{\beta}(2)^{\perp}\cong\mathfrak{e}_7$.  
For any $\nu\in{\rm Aut}(\mathfrak{g})$ satisfying $\nu(\mathfrak{k}_1)=
\mathfrak{k}_2$, it follows from (\ref{eqn9-1}) that 
$\nu(\mathfrak{su}_{\alpha_1}(2))=\mathfrak{su}_{\beta}(2)$. Hence there is 
no automorphism in ${\rm Aut}_{\mathfrak h}(\mathfrak{g})$ such that it maps 
$\mathfrak{k}_1$ to $\mathfrak{k}_2$ because 
$\mathfrak{su}_{\alpha_1}(2)\subset\mathfrak{h}$ and 
$\mathfrak{su}_{\beta}(2)\not\subset\mathfrak{h}$. 
\vspace{0.2cm}
\par

Next we shall show that $\tau_{K_5}$ is conjugate within 
${\rm Aut}_{\mathfrak h} ({\mathfrak{g}})$ to $\tau_{K_2 + K_5}$. 
Set 
$$
\begin{array}{l}
\gamma_1 :=\alpha_8, \ \gamma_2 :=-\alpha_1 -\alpha_2
-2\alpha_3 -3\alpha_4 -3\alpha_5 -3\alpha_6-2\alpha_7 -\alpha_8, \\
\gamma_3:=\alpha_7, \ \gamma_4 :=\alpha_8, \ \gamma_5 :=\alpha_5, \
\gamma_6 :=\alpha_4, \ \gamma_7 :=\alpha_3, \ \gamma_8 :=\alpha_1.
\end{array}
$$ 
It is easy to see that $\varPi' :=\{\gamma_1,\ldots ,\gamma_8\}$ is a
fundamental root system of $\mathfrak{e}_8$ (cf.\cite{FV}). Therefore
there exists a unique $\nu$ in $W(\mathfrak{g, h})$ such that $\nu(\varPi)
=\varPi'$. Hence we have $\nu(\alpha_i)=\gamma_i \ (i=1,\ldots ,8)$.
Then it is easy to see that $\nu^{-1}(K_2) =
-K_2$ and $\nu^{-1}(K_5) = -3K_2 + K_5$. Hence $\nu^{-1}\in
{\rm Aut}_{\mathfrak h}(\mathfrak{g})$ and 
 $\tau_{K_5}\approx\tau_{K_2+K_5}$. 
 
\vspace{0.3cm}\par

{\it The case where $\mathfrak{g=e}_7$ and $\sigma=\tau_{(1/2)K_3}$.} From Proposition \ref{prop:81}, we can see that $\mathfrak{g}^{\tau_{K_4}}
\cong\mathfrak{g}^{\tau_{K_3+K_4}}$ and $\mathfrak{h\cap g}^{\tau_{K_4}}
\cong\mathfrak{h\cap g}^{\tau_{K_3+K_4}}$. However $\tau_{K_4}$ is not 
conjugate within ${\rm Aut}_{\mathfrak h} (\mathfrak{g})$ to 
$\tau_{K_3+K_4}$. Indeed, note that
 $\mathfrak{g}^{\tau_{K_4}}\cong\mathfrak{g}^{\tau_{K_3+K_4}}\cong
  A_1\oplus D_6$, 
where $A_1\subset\mathfrak{g}^{\tau_{K_4}}$ and 
$A_1\subset\mathfrak{g}^{\tau_{K_3+K_4}}$ coincide with 
$\mathfrak{su}_{\alpha_2}(2)(\subset\mathfrak h)$ and 
$\mathfrak{su}_{\beta}(2)(\not\subset\mathfrak h)$, respectively. Here 
$\beta=\alpha_1+2\alpha_2+2\alpha_3+4\alpha_4+3\alpha_5+2\alpha_6+\alpha_7$. 
Therefore, there is no automorphism in ${\rm Aut}_{\mathfrak h}(\mathfrak{g})$ 
which maps $\mathfrak{g}^{\tau_{K_4}}$ to $\mathfrak{g}^{\tau_{K_3+K_4}}$, 
since $\nu(\mathfrak{su}_{\alpha_2}(2))=\mathfrak{su}_{\beta}(2)$ for any 
$\nu\in{\rm Aut}(\mathfrak{g})$ satisfying
$\nu(\mathfrak{g}^{\tau_{K_4}})=\mathfrak{g}^{\tau_{K_3+K_4}}$.

\subsection{$\dim\mathfrak{z}=0, \ \tau\circ\sigma = \sigma^{-1}\circ\tau$}
 {\it The case where $\mathfrak{g=e}_8$ and $\sigma=\tau_{(1/2)K_6}$.} 
For the reflection $t_{\alpha_1}\in{\rm Int}(\mathfrak{h})$, 
it follows from Proposition \ref{prop:51} that 
$t_{\alpha_1}^{-1}\circ\tau^{\varPi}_2\circ t_{\alpha_1}(E_{\alpha_3})=
-\tau^{\varPi}_2(E_{\alpha_3})$ and $t_{\alpha_1}^{-1}\circ\tau^{\varPi}_2\circ
 t_{\alpha_1}(E_{\alpha_i})=\tau^{\varPi}_2(E_{\alpha_i}) \ (i\ne 3)$. Hence we get 
\begin{equation}
t_{\alpha_1}^{-1}\circ\tau^{\varPi}_2\circ t_{\alpha_1}=\tau^{\varPi}_2\circ\tau_{K_3},
\label{eqn:92}
\end{equation}
and therefore $\tau^{\varPi}_2\circ\tau_{K_3}\approx\tau^{\varPi}_2\approx\tau^{\varPi}_2\circ\sigma
\approx\tau^{\varPi}_2\circ\tau_{K_3}\circ\sigma$.
  By an argument similar to the case $\mathfrak{g=e}_7$, for the
 reflection $t_{\alpha_3}\in{\rm Int}(\mathfrak{su}_{\alpha_3}(2))
\subset{\rm Int}(\mathfrak h)$, it follows that  
$$
\tau^{\varPi}_2\circ t_{\alpha_3}(E_{\alpha_i})=-t_{\alpha_3}\circ\tau^{\varPi}_2(E_{\alpha_i}),\ \ 
\tau^{\varPi}_2\circ t_{\alpha_3}(E_{\alpha_j})=t_{\alpha_3}\circ\tau^{\varPi}_2(E_{\alpha_j}),
$$
where $i=1,4$ and $j\ne1,4$. Therefore 
\begin{equation}
t_{\alpha_3}^{-1}\circ\tau^{\varPi}_2\circ t_{\alpha_3}=\tau^{\varPi}_2\circ\tau_{K_1+K_4},
\label{eqn:93}
\end{equation}
which implies $\tau^{\varPi}_2\circ\tau_{K_1+K_4}\approx\tau^{\varPi}_2\approx\tau^{\varPi}_2\circ\sigma
\approx\tau^{\varPi}_2\circ\tau_{K_1+K_4}\circ\sigma$.
\par
 Next, we shall prove 
\begin{equation}
\begin{array}{ll}
 \tau^{\varPi}_2\circ\tau_{K_1+K_8}\approx\tau^{\varPi}_2\circ\tau_{K_4+K_8}\approx
 \tau^{\varPi}_2\circ\tau_{K_3+K_4+K_8}\approx
 \tau^{\varPi}_2\circ\tau_{K_1+K_3+K_8}\circ\sigma, \\
 \tau^{\varPi}_2\circ\tau_{K_1+K_8}\circ\sigma\approx
 \tau^{\varPi}_2\circ\tau_{K_4+K_8}\circ\sigma\approx
 \tau^{\varPi}_2\circ\tau_{K_3+K_4+K_8}\circ\sigma\approx\tau^{\varPi}_2\circ\tau_{K_1+K_3+K_8}.
\end{array}
\label{eqn:94}
\end{equation}
It is easy to see that $t_{\alpha_i}(K_4+K_8)=K_4+K_8$ ($i=1,3$), and 
it follows from (\ref{eqn:92}) and (\ref{eqn:93}) that 
\begin{eqnarray*}
\tau^{\varPi}_2 \circ\tau_{K_1+K_8}&\!\!\! =&\!\!\!\tau^{\varPi}_2 \circ\tau_{K_1+K_4}
 \circ\tau_{K_4+K_8} =t_{\alpha_3}^{-1} \circ\tau^{\varPi}_2 \circ t_{\alpha_3}
 \circ\tau_{K_4+K_8} \\ 
&\!\!\! =&\!\!\! t_{\alpha_3}^{-1} \circ\tau^{\varPi}_2 \circ\tau_{K_4+K_8} \circ
 t_{\alpha_3},\\ 
\tau^{\varPi}_2 \circ\tau_{K_3+K_4+K_8}&\!\!\! =&\!\!\! t_{\alpha_1}^{-1}
 \circ\tau^{\varPi}_2 \circ t_{\alpha_1} \circ\tau_{K_4+K_8} =t_{\alpha_1}^{-1}
 \circ\tau^{\varPi}_2 \circ\tau_{K_4+K_8} \circ t_{\alpha_1}.
\end{eqnarray*}
For $t_{\alpha_4}\in{\rm Int}(\mathfrak h)$, we get 
$$\tau^{\varPi}_2\circ t_{\alpha_4}(E_{\alpha_3})
=-t_{\alpha_4}\circ\tau^{\varPi}_2(E_{\alpha_3}), \ \
\tau^{\varPi}_2\circ t_{\alpha_4}(E_{\alpha_i})=t_{\alpha_4}\circ\tau^{\varPi}_2(E_{\alpha_i}),
$$
where $i=1,4,6,7,8$. Moreover, we obtain 
\begin{eqnarray*}
\tau^{\varPi}_2\circ
 t_{\alpha_4}(E_{\alpha_2})&\!\!\! =&\!\!\! b_2\tau^{\varPi}_2(E_{\alpha_2+\alpha_4})
=b_2kE_{\alpha_4+\alpha_5},\\
t_{\alpha_4}\circ\tau^{\varPi}_2(E_{\alpha_2})&\!\!\! =&\!\!\! t_{\alpha_4}(E_{\alpha_5})
=b_5E_{\alpha_4+\alpha_5},\\
\tau^{\varPi}_2\circ t_{\alpha_4}(E_{\alpha_5})&\!\!\! =&\!\!\! b_5\tau^{\varPi}_2(E_{\alpha_4+\alpha_5})
=b_5k^{-1}E_{\alpha_2+\alpha_4},\\
t_{\alpha_4}\circ\tau^{\varPi}_2(E_{\alpha_5})&\!\!\! =&\!\!\! t_{\alpha_4}(E_{\alpha_2})
=b_2E_{\alpha_2+\alpha_4},
\end{eqnarray*}
for some $b_2, \ b_5, \ k\in{\mathbb{C}}$ with $|b_2|=|b_5|=|k|=1$.  Put
$a:=b_2k/b_5$. Then we have
$$\tau^{\varPi}_2\circ  t_{\alpha_4}(E_{\alpha_2})
=at_{\alpha_4}\circ\tau^{\varPi}_2(E_{\alpha_2}),\ \ 
\tau^{\varPi}_2\circ t_{\alpha_4}(E_{\alpha_5})
=a^{-1}t_{\alpha_4}\circ\tau^{\varPi}_2(E_{\alpha_5}).
$$
Hence we have $t_{\alpha_4}^{-1}\circ\tau^{\varPi}_2\circ t_{\alpha_4} =
\tau^{\varPi}_2\circ\tau_{K_3}\circ\tau_{s(K_2-K_5)}$, where 
$a=e^{s\pi\sqrt{-1}}$. Note that $s\in\mathbb Z$ since
$(t_{\alpha_4}^{-1}\circ\tau^{\varPi}_2\circ t_{\alpha_4})^2 = {\rm Id}$, and thus 
we may assume $s=0$ or $1$. If $s=0$, then we have 
$t_{\alpha_4}^{-1}\circ\tau^{\varPi}_2\circ t_{\alpha_4}\circ\tau_{K_1+K_8}
=t_{\alpha_4}^{-1}\circ\tau^{\varPi}_2\circ\tau_{K_1+K_8}\circ t_{\alpha_4}
=\tau^{\varPi}_2\circ\tau_{K_1+K_3+K_8}$, which contradicts Proposition \ref{prop:51}. 
Thus $t_{\alpha_4}^{-1}\circ\tau^{\varPi}_2\circ t_{\alpha_4}
=\tau^{\varPi}_2\circ\tau_{K_3}\circ\tau_{K_2-K_5}$ and 
$$\begin{array}{lll}
t_{\alpha_4}^{-1} \circ\tau^{\varPi}_2 \circ\tau_{K_1+K_8} \circ
 t_{\alpha_4}&\!\!\! =&\!\!\!t_{\alpha_4}^{-1} \circ\tau^{\varPi}_2 \circ
 t_{\alpha_4} \circ\tau_{K_1+K_8} =\tau^{\varPi}_2 \circ\tau_{K_3+K_2-K_5+K_1+K_8}\\
&\!\!\! =&\!\!\!\tau^{\varPi}_2 \circ\tau_{K_1+K_3+K_8} \circ\tau_{K_7}
 \circ\tau_{K_2-K_5+K_7}\\ 
&\!\!\! =&\!\!\!\tau^{\varPi}_2 \circ\tau_{K_1+K_3+K_8} \circ\sigma
 \circ\tau_{(-(1/2)K_6+K_7)+(K_2-K_5+K_7))}. 
\end{array}$$
It is easy to see that $\sqrt{-1}h:=\sqrt{-1}(-(1/2)K_6+K_7+K_2-K_5+K_7)$ is 
a $(-1)$-eigenvector of $\tau^{\varPi}_2$. Therefore $t_{\alpha_4}^{-1}\circ\tau^{\varPi}_2\circ
\tau_{K_1+K_8}\circ\tau_{\alpha_4}=\tau_{-(1/2)h}\circ\tau^{\varPi}_2\circ
\tau_{K_1+K_3+K_8}\circ\sigma\circ\tau_{(1/2)h}$, which implies 
$\tau^{\varPi}_2\circ\tau_{K_1+K_8}\approx\tau^{\varPi}_2\circ\tau_{K_1+K_3+K_8}\circ\sigma$.  
Moreover since 
$$\tau^{\varPi}_2\circ\tau_{K_1+K_8}\circ\sigma\approx\tau^{\varPi}_2\circ\tau_{K_1+K_3+K_8}
\circ\sigma^2=\sigma^{-1}\circ\tau^{\varPi}_2\circ\tau_{K_1+K_3+K_8}\circ\sigma,$$
we have $\tau^{\varPi}_2\circ\tau_{K_1+K_8}\circ\sigma\approx
\tau^{\varPi}_2\circ\tau_{K_1+K_3+K_8}$. We have thus proved (\ref{eqn:94}). \par
 Finally, by using $t_{\alpha_8}\in{\rm Int}(\mathfrak h)$, we shall show that 
$\tau^{\varPi}_2\circ\tau_{K_1+K_3+K_4}\approx
\tau^{\varPi}_2\circ\tau_{K_1+K_3+K_4}\circ\sigma$. It is easy to see that 
$\tau^{\varPi}_2\circ t_{\alpha_8}(E_{\alpha_i})=t_{\alpha_8}\circ
 \tau^{\varPi}_2(E_{\alpha_i})$, $i=1,\ldots ,6$, and 
$$
  \begin{array}{l}
     \tau^{\varPi}_2\circ t_{\alpha_8}(E_{\alpha_7})
      =b_7\tau^{\varPi}_2(E_{\alpha_7+\alpha_8}) =b_7kE_{\alpha_0+\alpha_8},\\
      t_{\alpha_8} \circ\tau^{\varPi}_2(E_{\alpha_7}) =t_{\alpha_8}(E_{\alpha_0})
       =b_0E_{\alpha_0+\alpha_5},
   \end{array}
$$
for some $b_0$, $b_7$, $k\in{\mathbb{C}}$ with $|b_0|=|b_7|=|k|=1$. 
By an argument similar to the above, we obtain 
$t_{\alpha_8}^{-1}\circ\tau^{\varPi}_2\circ t_{\alpha_8}=\tau^{\varPi}_2$
or $\tau^{\varPi}_2\circ\tau_{K_7}$. If $t_{\alpha_8}^{-1}\circ\tau^{\varPi}_2\circ
t_{\alpha_8}=\tau^{\varPi}_2$, then since $t_{\alpha_8}(K_1)=K_1$, 
$t_{\alpha_8}(K_8)=K_7-K_8$, it follows that 
\begin{eqnarray*}
   \tau^{\varPi}_2\circ\tau_{K_1} \circ\tau_{K_8}
   &\!\!\! =&\!\!\! t_{\alpha_8}^{-1} \circ\tau^{\varPi}_2 \circ t_{\alpha_8}
   \circ\tau_{K_1+K_8} =t_{\alpha_8}^{-1} \circ\tau^{\varPi}_2
   \circ\tau_{K_1+K_7+K_8} \circ t_{\alpha_8}\\ 
   &\!\!\! =&\!\!\! t_{\alpha_8}^{-1} \circ\tau^{\varPi}_2 \circ\tau_{K_1+K_8}
    \circ\sigma \circ\tau_{K_7-(1/2)K_6} \circ t_{\alpha_8}.
\end{eqnarray*}
Therefore $\tau^{\varPi}_2\circ\tau_{K_1+K_8}\approx \tau^{\varPi}_2\circ\tau_{K_1+K_8}\circ
\sigma$ since $\sqrt{-1}(K_7-(1/2)K_6)$ is a $(-1)$-eigenvector of $\tau^{\varPi}_2$. 
This contradicts Proposition \ref{prop:51}, and hence 
$t_{\alpha_8}^{-1}\circ\tau^{\varPi}_2\circ t_{\alpha_8}=\tau^{\varPi}_2\circ\tau_{K_7}$. Thus 
\begin{eqnarray*}
t_{\alpha_8}^{-1} \circ\tau^{\varPi}_2 \circ\tau_{K_1+K_3+K_4} \circ
 t_{\alpha_8}&\!\!\! =&\!\!\! t_{\alpha_8}^{-1} \circ\tau^{\varPi}_2 \circ
 t_{\alpha_8} \circ\tau_{K_1+K_3+K_4} =\tau^{\varPi}_2 \circ\tau_{K_1+K_3+K_4}
 \circ\tau_{K_7}\\ 
 &\!\!\! =&\!\!\! \tau^{\varPi}_2 \circ\tau_{K_1+K_3+K_4} \circ\sigma
  \circ\tau_{K_7-(1/2)K_6} \approx \tau^{\varPi}_2 \circ\tau_{K_1+K_3+K_4} \circ\sigma,
\end{eqnarray*}
which implies $\tau^{\varPi}_2\circ \tau_{K_1+K_3+K_4}\approx
\tau^{\varPi}_2\circ\tau_{K_1+K_3+K_4}\circ\sigma$.\vspace{0.3cm}
\par

{\it The case where $\mathfrak{g=e}_7$ and $\sigma=
\tau_{(1/2)K_4}$.} We consider
$\tau^{\varPi}_3\circ\tau_{K_1+K_2}, \ \tau^{\varPi}_3\circ\tau_{K_2+K_6}, \
\tau^{\varPi}_3\circ\tau_{K_1+K_2}\circ\sigma$ and
$\tau^{\varPi}_3\circ\tau_{K_2+K_6}\circ\sigma$ (see Proposition \ref{prop:51}). 
Let $\varphi$ be the automorphism of $\mathfrak{g}$ given by (\ref{eqn:510}).
It follows from (\ref{eqn:58}) and (\ref{eqn:510}) that $\varphi\circ\tau^{\varPi}_3=
\tau^{\varPi}_3\circ\varphi$. Since $\varphi(K_1+K_2)=K_2+K_6-4K_7$ and 
$\varphi(K_4)=K_4-4K_7$, we have $\varphi\circ\sigma=\sigma\circ\varphi$ and 
$$\varphi\circ\tau^{\varPi}_3\circ\tau_{K_1+K_2}\circ\varphi^{-1}
=\tau^{\varPi}_3\circ\tau_{K_2+K_6},$$
and therefore 
$$\varphi\circ\tau^{\varPi}_3\circ\tau_{K_1+K_2}\circ\sigma\circ\varphi^{-1}
=\varphi\circ\tau^{\varPi}_3\circ\tau_{K_1+K_2}\circ\varphi^{-1}\circ\sigma
=\tau^{\varPi}_3\circ\tau_{K_2+K_6}\circ\sigma.$$
Thus we obtain $\tau^{\varPi}_3\circ\tau_{K_1+K_2}\approx\tau^{\varPi}_3\circ\tau_{K_2+K_6}$
and
$\tau^{\varPi}_3\circ\tau_{K_1+K_2}\circ\sigma\approx
\tau^{\varPi}_3\circ\tau_{K_2+K_6}\circ\sigma$.\par

Next considering reflection $t_{\alpha_1}\in{\rm
Int}(\mathfrak{su}_{\alpha_1}(2))\subset{\rm Int}(\mathfrak h)$, we get
$t_{\alpha_1}(E_{\alpha_i}) = E_{\alpha_i}$ for $i=2,4,5,6,7$ because 
$\alpha_1\pm\alpha_i$ are not roots. Put
$\beta :=\alpha_1 + \alpha_2 + 2\alpha_3 + 3\alpha_4 + 2\alpha_5 +
\alpha_6$. Then $\tau^{\varPi}_3(\alpha_4) = \beta$ and 
$\beta\pm\alpha_1\notin\varDelta(\mathfrak{g}_{\mathbb{C}},
\mathfrak{t}_{\mathbb{C}})$. Hence we have $t_{\alpha_1}(E_{\beta}) =
E_{\beta}$. Since 
$$
\left\{\begin{array}{l} \tau^{\varPi}_3\circ
t_{\alpha_1}(E_{\alpha_1})=b_1\tau^{\varPi}_3(E_{-\alpha_1})=-b_1E_{-\alpha_1}, \\
t_{\alpha_1}\circ\tau^{\varPi}_3(E_{\alpha_1})=-t_{\alpha_1}(E_{\alpha_1}))
=-b_1E_{-\alpha_1},
\end{array}\right.
\left\{\begin{array}{l}
\tau^{\varPi}_3\circ t_{\alpha_1}(E_{\alpha_3})=b_3\tau^{\varPi}_3(E_{\alpha_1+\alpha_3})
=b_3kE_{\alpha_0 + \alpha_1}, \\
t_{\alpha_1}\circ\tau^{\varPi}_3(E_{\alpha_3})=t_{\alpha_1}(E_{\alpha_0}))
=b_0E_{\alpha_0+\alpha_1},
\end{array}\right.
$$
for some $b_0,b_1,b_3,k\in{\mathbb{C}}$ with $|b_0|=|b_1|=|b_3|=|k|=1$, 
there exists $s\in\mathbb{R}$ such that 
$t_{\alpha_1}^{-1}\circ\tau^{\varPi}_3\circ t_{\alpha_1}=\tau^{\varPi}_3\circ\tau_{sK_3}$.
Moreover, since $(t_{\alpha_1}^{-1}\circ\tau^{\varPi}_3\circ t_{\alpha_1})^2={\rm Id}$ 
and $\tau^{\varPi}_3(K_3)=-3K_3 + 2K_4$, we get $s\in\mathbb Z$, and thus
$t_{\alpha_1}^{-1}\circ\tau^{\varPi}_3\circ t_{\alpha_1}=\tau^{\varPi}_3$ or
$\tau^{\varPi}_3\circ\tau_{K_3}$. If $t_{\alpha_1}^{-1}\circ\tau^{\varPi}_3\circ
t_{\alpha_1}=\tau^{\varPi}_3$, then
$$\begin{array}{lll}
\tau^{\varPi}_3\circ\tau_{K_1+K_6}&\!\!\! =&\!\!\! t_{\alpha_1}^{-1}\circ\tau^{\varPi}_3\circ
 t_{\alpha_1}\circ\tau_{K_1+K_6}=t_{\alpha_1}^{-1}\circ\tau^{\varPi}_3\circ
\tau_{t_{\alpha_1}(K_1+K_6)}\circ t_{\alpha_1}\\
&\!\!\! =&\!\!\!
 t_{\alpha_1}^{-1} \circ\tau^{\varPi}_3 \circ\tau_{K_1+K_3} \circ\tau_{K_6} \circ
 t_{\alpha_1} =t_{\alpha_1}^{-1} \circ\tau^{\varPi}_3 \circ\tau_{K_1+K_6}
 \circ\tau_{K_3} \circ t_{\alpha_1}\\ 
&\!\!\! =&\!\!\! t_{\alpha_1}^{-1} \circ\tau^{\varPi}_3 \circ\tau_{K_1+K_6}
 \circ\sigma \circ \tau_{K_3-(1/2)K_4} \circ t_{\alpha_1}.
\end{array}
$$
Put $\mu:=\tau_{(1/2)(K_3-(1/2)K_4)}$. Then since
$K_3-(1/2)K_4$ is a $(-1)$-eigenvector of $\tau^{\varPi}_3$, it is easy to see that 
$$\tau^{\varPi}_3\circ\tau_{K_1+K_6}=
t_{\alpha_1}^{-1}\circ\mu^{-1}\circ\tau^{\varPi}_3\circ\tau_{K_1+K_6}\circ\sigma\circ
\mu\circ t_{\alpha_1}.$$
Because $\mu\circ t_{\alpha_1}\in{\rm Int}(\mathfrak h)$, it follows that 
$\tau^{\varPi}_3\circ\tau_{K_1+K_6}\approx\tau^{\varPi}_3\circ\tau_{K_1+K_6}\circ\sigma$,
which contradics Proposition \ref{prop:51}. Hence 
$t_{\alpha_1}^{-1}\circ\tau^{\varPi}_3\circ t_{\alpha_1}=\tau^{\varPi}_3\circ\tau_{K_3}$, 
and 
\begin{eqnarray*}
t_{\alpha_1}^{-1} \circ\tau^{\varPi}_3 \circ\tau_{K_2+K_6} \circ
 t_{\alpha_1}&\!\!\! =&\!\!\! t_{\alpha_1}^{-1} \circ\tau^{\varPi}_3 \circ
 t_{\alpha_1} \circ\tau_{K_2+K_6} =\tau^{\varPi}_3 \circ\tau_{K_3}
 \circ\tau_{K_2+K_6} \\ 
&\!\!\! =&\!\!\! \tau^{\varPi}_3 \circ\sigma \circ\mu^2 \circ\tau_{K_2+K_6}
 =\tau^{\varPi}_3 \circ\tau_{K_2+K_6} \circ\sigma \circ\mu^2\\
&\!\!\! =&\!\!\! \mu^{-1} \circ\tau^{\varPi}_3 \circ\tau_{K_2+K_6} \circ\sigma \circ\mu,
\end{eqnarray*}
since $t_{\alpha_1}(K_2+K_6)=K_2+K_6$ and $\tau^{\varPi}_3(K_3-(1/2)K_4)=
-(K_3-(1/2)K_4)$. Consequently we obtain $\tau^{\varPi}_3\approx\tau^{\varPi}_3\circ\tau_{K_3}$ 
and $\tau^{\varPi}_3\circ\tau_{K_2+K_6}\approx\tau^{\varPi}_3\circ\tau_{K_2+K_6}\circ\sigma$.


\section{Classifications}
 From Propositions \ref{prop:51}, \ref{prop:61}, \ref{prop:71} 
and \ref{prop:81} together with the results in Section 9, we obtain 
the following theorem which gives the complete classification of involutions 
preserving $\mathfrak h$. 
\begin{thm}
  Let $(G/H,\langle,\rangle,\sigma)$ be a Riemannian $4$-symmetric space
 such that $G$ is compact and simple. Suppose that 
$\sigma={\rm Ad}(\exp(\pi/2)\sqrt{-1}K_i)$ for some $\alpha_i\in
\varPi(\mathfrak{g}_{\mathbb{C}},\mathfrak{t}_{\mathbb{C}})$ with 
$m_i=3$ or $4$. Then the following Tables 5, 6, 7 and 8 give the 
complete lists of the equivalence classes within 
${\rm Aut}_{\mathfrak h}(\mathfrak{g})$ of involutions $\tau$ satisfying 
$\tau(\mathfrak{h})=\mathfrak{h}$. 
\end{thm}
{\small
\begin{center}
\renewcommand{\arraystretch}{1.2}
\begin{longtable}{l|l|l|l}
\caption{$\dim\mathfrak{z} =0, \ \tau\circ\sigma = \sigma^{-1}\circ\tau$, 
$\sigma=\tau_{(1/2)H}$ and $\mathfrak{k=g}^{\tau}$.}\\
\hline
\multicolumn{1}{c|}{$(\mathfrak{g}, \mathfrak{h}, H)$} &
 \multicolumn{1}{|c|}{$\tau$}
 &\multicolumn{1}{|c|}{$\mathfrak k$} & 
 \multicolumn{1}{|c}{$\mathfrak {h\cap k}$}\\
\hline
\hline

$(\mathfrak{e}_8,\mathfrak{su}(8)\oplus\mathfrak{su}(2), K_3)$ & $\tau^{\varPi}_1$
 & $\mathfrak{so}(16)$ & $\mathfrak{so}(8)\oplus\mathfrak{so}(2)$ \\

& $\tau^{\varPi}_1\circ\tau_{K_6}$ & $\mathfrak{e}_7\oplus\mathfrak{su}(2)$ &
 $\mathfrak{sp}(4)\oplus\mathfrak{so}(2)$ \\

& $\tau^{\varPi}_1\circ\tau_{K_6+(1/2)K_3}$ & $\mathfrak{so}(16)$ &
 $\mathfrak{sp}(4)\oplus\mathfrak{so}(2)$ \\
\hline
$(\mathfrak{e}_8,\mathfrak{so}(10)\oplus\mathfrak{so}(6), K_6)$ &
 $\tau^{\varPi}_2$ & $\mathfrak{so}(16)$ &
 $\begin{array}{l}\!\!\!(\mathfrak{so}(5)+\mathfrak{so}(5))\\
\hspace{12pt}\oplus(\mathfrak{so}(3)+\mathfrak{so}(3))\end{array}$
 \\ 

& $\tau^{\varPi}_2\circ\tau_{K_1+K_8}$ & $\mathfrak{e}_7\oplus\mathfrak{su}(2)$ &
 $(\mathfrak{so}(7)+\mathfrak{so}(3))\oplus\mathfrak{so}(5)$ \\

& $\tau^{\varPi}_2\circ\tau_{K_1+K_3+K_4}$ &  $\mathfrak{e}_7\oplus\mathfrak{su}(2)$ &
 $\mathfrak{so}(9)\oplus(\mathfrak{so}(3)+\mathfrak{so}(3))$ \\

& $\tau^{\varPi}_2\circ\tau_{K_1+K_3+K_8}$ & $\mathfrak{so}(16)$  &
 $(\mathfrak{so}(7)+\mathfrak{so}(3))\oplus\mathfrak{so}(5)$ \\ 
\hline
$(\mathfrak{e}_7,\mathfrak{so}(6)\oplus\mathfrak{so}(6)\oplus\mathfrak{su}(2),
 K_4)$
 & $\tau^{\varPi}_3$ & $\mathfrak{su}(8)$ &
 $\begin{array}{l}\!\!\!(\mathfrak{so}(3)+\mathfrak{so}(3))\\
\hspace{4pt}\oplus(\mathfrak{so}(3)+\mathfrak{so}(3))\oplus\mathfrak{so}(2)
\end{array}$ \\

& $\tau^{\varPi}_3\circ\tau_{K_1+K_2}$ &
 $\mathfrak{so}(12)\oplus\mathfrak{su}(2)$ &
 $\begin{array}{l}\!\!\!\mathfrak{so}(5)\\
\hspace{4pt}\oplus(\mathfrak{so}(3)+\mathfrak{so}(3)) \oplus
 \mathfrak{su}(2)\end{array}$ 
 \\

& $\tau^{\varPi}_3\circ\tau_{K_1+K_6}$ & $\mathfrak{e}_6\oplus \mathbb{R}$ &
 $\mathfrak{so}(5)\oplus\mathfrak{so}(5)\oplus\mathfrak{su}(2)$ \\ 

& $\tau^{\varPi}_3\circ\tau_{K_1+K_6+(1/2)K_4}$ & $\mathfrak{su}(8)$ &
 $\mathfrak{so}(5)\oplus\mathfrak{so}(5)\oplus\mathfrak{su}(2)$ \\ 

\cline{2-4}

& $\tau^{\varPi}_3\circ\varphi$ &
 $\mathfrak{so}(12)\oplus\mathfrak{su}(2)$ &
 $\mathfrak{so}(6)\oplus\mathfrak{so}(2)$ \\

& $\tau^{\varPi}_3\circ\varphi\circ\tau_{(1/2)K_4}$ &
 $\mathfrak{su}(8)$ &
 $\mathfrak{so}(6)\oplus\mathfrak{so}(2)$ \\
\hline
$(\mathfrak{f}_4,\mathfrak{so}(6)\oplus\mathfrak{so}(3), K_3)$ & $\tau^{\varPi}_4$
& $\mathfrak{sp}(3)\oplus\mathfrak{su}(2)$ &
 $(\mathfrak{so}(3)+\mathfrak{so}(3))\oplus\mathfrak{so}(2)$ \\

& $\tau^{\varPi}_4\circ\tau_{K_1+K_4}$ & $\mathfrak{so}(9)$ &
 $\mathfrak{so}(5)\oplus\mathfrak{so}(3)$ \\ 

& $\tau^{\varPi}_4\circ\tau_{K_1+K_4+(1/2)K_3}$ &
 $\mathfrak{sp}(3)\oplus\mathfrak{su}(2)$ &
 $\mathfrak{so}(5)\oplus\mathfrak{so}(3)$ \\ 
\hline

\multicolumn{4}{l}{$\tau^{\varPi}_1\ :\ E_{\alpha_1}\mapsto -E_{\alpha_1},\ 
 E_{\alpha_2}\mapsto E_{\alpha_0},\  E_{\alpha_3}\mapsto c_1 E_{\beta_1},\ 
 E_{\alpha_4}\mapsto E_{\alpha_8},\  E_{\alpha_5}\mapsto E_{\alpha_7},\ 
 E_{\alpha_6}\mapsto -E_{\alpha_6},$} \\
\multicolumn{4}{l}{\hspace{25pt}($\beta_1=\alpha_1 + 2\alpha_2 + 3\alpha_3
 + 4\alpha_4 + 3\alpha_5 + 2\alpha_6 + \alpha_7$)}\\
\multicolumn{4}{l}{$\tau^{\varPi}_2\ : \ E_{\alpha_1}\mapsto
 -E_{\alpha_1}, \ E_{\alpha_2}\mapsto E_{\alpha_5}, \
  E_{\alpha_3}\mapsto -E_{\alpha_3}, \ 
 E_{\alpha_4}\mapsto -E_{\alpha_4}, \ E_{\alpha_6}\mapsto c_2
 E_{\beta_2}, \ E_{\alpha_7}\mapsto E_{\alpha_0},$} \\ 
\multicolumn{4}{l}{\hspace{25pt}$E_{\alpha_8}\mapsto
 -E_{\alpha_8},$ ($\beta_2=\alpha_1 + \alpha_2 + 2\alpha_3 + 3\alpha_4 +
 3\alpha_5 + 3\alpha_6 + 2\alpha_7 +\alpha_8$)} \\
\multicolumn{4}{l}{$\tau^{\varPi}_3\ :\ E_{\alpha_1}\mapsto -E_{\alpha_1}, \
 E_{\alpha_2}\mapsto -E_{\alpha_2}, \ E_{\alpha_3}\mapsto E_{\alpha_0}, \
 E_{\alpha_4}\mapsto c_3 E_{\beta_3}, \ E_{\alpha_5}\mapsto E_{\alpha_7}, \
 E_{\alpha_6}\mapsto -E_{\alpha_6},$}\\
\multicolumn{4}{l}{\hspace{25pt}($\beta_3=\alpha_1 + \alpha_2 +
 2\alpha_3 + 3\alpha_4 + 2\alpha_5 + \alpha_6$)} \\
\multicolumn{4}{l}{$\tau^{\varPi}_4\ :\ E_{\alpha_1}\mapsto -E_{\alpha_1},\ 
 E_{\alpha_2}\mapsto E_{\alpha_0},\  E_{\alpha_3}\mapsto c_4 E_{\beta_4},\ 
 E_{\alpha_4}\mapsto -E_{\alpha_4},$ ($\beta_4=\alpha_1 + 2\alpha_2 +
 3\alpha_3 + \alpha_4$)}\\
\multicolumn{4}{l}{where $c_i(i=1,2,3,4)$ is some complex number with
 $|c_i|=1$.} 

\end{longtable}
\end{center}
}

{\small
\begin{center}
\renewcommand{\arraystretch}{1.2}
\begin{longtable}{l|l|l|l|l}
\caption{$\dim\mathfrak{z}=1, \ \tau\circ\sigma = \sigma^{-1}\circ\tau$, 
$\sigma=\tau_{(1/2)H}$ and $\mathfrak{k=g}^{\tau}$.}\\
\hline
\multicolumn{1}{c|}{$(\mathfrak{g}, \mathfrak{h}, H)$} &
 \multicolumn{1}{|c|}{$\varPi$}
 &\multicolumn{1}{|c|}{$\varPi_1$} &
 \multicolumn{1}{|c|}{$(\mathfrak{g}^*,
 \mathfrak k$)} &
 \multicolumn{1}{|c}{$\mathfrak{h\cap k}$}  
 \\
\hline
\hline
$(\mathfrak{e}_8,
 \mathfrak{su}(8)\oplus\mathbb{R}, K_2)$ & $E_8$ & $\alpha_2$ & $(\mathfrak{e}
 _{8(8)}, \mathfrak{so}(16))$ &
 $\mathfrak{so}(8)$ \\
\hline
\lw {$(\mathfrak{e}_8,
 \mathfrak{e}_6\oplus\mathfrak{su}(2)\oplus\mathbb{R}, K_7)$} & $E_8$ &
 $\alpha_7$ & $(\mathfrak{e} _{8(8)}, \mathfrak{so}(16))$ &
 $\mathfrak{sp}(4)$ \\
\cline{2-5}
 & $F_4$ &
 $\alpha_7$ & $(\mathfrak{e}_{8(-24)}, \mathfrak{e}
 _7\oplus\mathfrak{su}(2))$ & $\mathfrak{f} _4$ \\
\hline
\lw {$(\mathfrak{e}_7,
 \mathfrak{su}(6)\oplus\mathfrak{su}(2)\oplus\mathbb{R}, K_3)$} &  $E_7$ & 
 $\alpha_3$ & $(\mathfrak{e}_{7(7)}, \mathfrak{su}(8))$ &
 $\mathfrak{so}(6)\oplus\mathfrak{so}(2)$ \\
\cline{2-5}
 & $F_4$ &
 $\alpha_3$ &
 $(\mathfrak{e} _{7(5)}, \mathfrak{so}(12)\oplus\mathfrak{su}(2))$ &
 $\mathfrak{sp}(3)\oplus\mathfrak{so}(2)$ \\
\hline
$(\mathfrak{e}_7,
 \mathfrak{su}(5)\oplus\mathfrak{su}(3)\oplus\mathbb{R}, K_5)$ & $E_7$ &
 $\alpha_5$ & $(\mathfrak{e}_{7(7)}, \mathfrak{su}(8))$ &
 $\mathfrak{so}(5)\oplus\mathfrak{so}(3)$ \\
\hline
\lw {$(\mathfrak{e}_6,
 \mathfrak{su}(3)\oplus\mathfrak{su}(3)\oplus\mathfrak{su}(2)\oplus\mathbb{R}, K_4)$} & $E_6$ & $\alpha_4$ & $(\mathfrak{e}_{6(6)}, \mathfrak{sp}(4))$ & 
 $\mathfrak{so}(3)\oplus\mathfrak{so}(3)\oplus\mathfrak{so}(2)$ \\
\cline{2-5}
 & $F_4$ & $\alpha_4$ & $(\mathfrak{e}_{6(2)},\mathfrak{su}(6)\oplus\mathfrak{su}(2))$ & $
 \mathfrak{su}(3)\oplus\mathfrak{su}(2) $ \\
\hline
$(\mathfrak{f}_4, \mathfrak{su}(3)\oplus
 \mathfrak{su}(2)\oplus\mathbb{R}, K_2)
 $ & $F_4$ & $\alpha_2$ & $(\mathfrak{f}_{4(4)}, \mathfrak{sp}(3)
 \oplus \mathfrak{su}(2))$ & $\mathfrak{so}(3)\oplus\mathfrak{so}(2)$ \\ 
\hline
$(\mathfrak{g}_2, \mathfrak{su}(2) \oplus\mathbb{R}, K_1)$ & $G_2$
 &$\alpha_1$ & $(\mathfrak{g}_{2(2)}, \mathfrak{su}(2)
 \oplus\mathfrak{su}(2))$ & $\mathfrak{so}(2)$ \\
\hline
\end{longtable}
\end{center}}

{\small
\begin{center}
\renewcommand{\arraystretch}{1.2}
\begin{longtable}{l|l|l|l}
\caption{$\dim\mathfrak{z}=0, \ \tau\circ\sigma=\sigma\circ\tau$,
$\sigma=\tau_{(1/2)H}$ and $\mathfrak{k=g}^{\tau}$.}\\
\hline
\multicolumn{1}{c|}{$(\mathfrak{g}, \mathfrak{h}, H)$} &
 \multicolumn{1}{|c|}{$h \ (\tau=\tau_h)$}
 &\multicolumn{1}{|c|}{$\mathfrak{k}$} & 
 \multicolumn{1}{|c}{$\mathfrak{h\cap k}$}  
 \\ 
\hline
\hline

$(\mathfrak{e}_8,
 \mathfrak{su}(8)\oplus\mathfrak{su}(2), K_3)$ & $K_1$ & $\mathfrak{so}(16)$
 & $\mathfrak{su}(8)\oplus\mathfrak{so}(2) $ \\

& $K_3$ & $\mathfrak{e} _7\oplus\mathfrak{su}(2)$  &
	     $\mathfrak{su}(8)\oplus\mathfrak{su}(2)$ \\  

& $K_4$ & $\mathfrak{e}_7\oplus\mathfrak{su}(2)$ &
	     $\mathfrak{sp}(4)\oplus \mathfrak{su}(2)$ \\

& $K_6$ & $\mathfrak{so}(16)$  & $\mathfrak{s(u}(4) + \mathfrak
	     u(4))\oplus\mathfrak{su}(2)$ \\ 

& $K_3 + K_4$ & $\mathfrak{so}(16)$ &
	     $\mathfrak{sp}(4)\oplus\mathfrak{su}(2)$ \\

& $K_3 + K_6$ & $\mathfrak{e}_7\oplus\mathfrak{su}(2)$ &
 $\mathfrak{s(u}(4) + \mathfrak u(4))\oplus\mathfrak{su}(2)$ \\  

& $K_1 + K_4$ &  $\mathfrak{e}_7\oplus\mathfrak{su}(2)$ &
 $\mathfrak{s(u}(6) + \mathfrak u(2))\oplus\mathfrak{so}(2)$ \\

& $K_1 + K_6$ &  $\mathfrak{so}(16)$  & $\mathfrak{s(u}(6) +
 \mathfrak u(2))\oplus\mathfrak{so}(2)$ \\

\hline
$(\mathfrak{e}_8,
 \mathfrak{so}(10)\oplus\mathfrak{so}(6), K_6)$ &  $K_1$ & $\mathfrak{so}(16)$  &
 $(\mathfrak{so}(8) + \mathfrak{so}(2))\oplus\mathfrak{so}(6)$ \\   

& $K_3$ & $\mathfrak{e}_7\oplus\mathfrak{su}(2)$  &
	     $(\mathfrak{so}(6) + \mathfrak{so}(4))\oplus\mathfrak{so}(6)$\\

& $K_6$ & $\mathfrak{so}(16)$ &
	     $\mathfrak{so}(10)\oplus \mathfrak{so}(6)$ \\

& $K_8$ & $\mathfrak{e}_7\oplus\mathfrak{su}(2)$  &
 $\mathfrak{so}(10)\oplus(\mathfrak{so}(4) + \mathfrak{so}(2))$ \\

& $K_1 + K_8$ &  $\mathfrak{e}_7\oplus\mathfrak{su}(2)$  &
 $(\mathfrak{so}(8) + \mathfrak{so}(2))\oplus(\mathfrak{so}(4)+
 \mathfrak{so}(2))$ \\ 

& $K_2 + K_7$ & $\mathfrak{e}_7\oplus\mathfrak{su}(2)$ &
	     $\mathfrak{u}(3)\oplus\mathfrak{u}(5)$ \\

& $K_3 + K_8$ &  $\mathfrak{so}(16)$  & $(\mathfrak{so}(6) +
	     \mathfrak{so}(4))\oplus(\mathfrak{so}(4) + \mathfrak{so}(2))$ \\

& $K_2 + K_6 + K_7$ & $\mathfrak{so}(16)$ &
	     $\mathfrak{u}(3)\oplus\mathfrak{u}(5)$ \\

\hline

$(\mathfrak{e}_7,
 \mathfrak{so}(6)\oplus\mathfrak{so}(6)\oplus\mathfrak{su}(2), K_4)$ 
& $K_1$ & $\mathfrak{so}(12)\oplus\mathfrak{su}(2)$  &
	 $\mathfrak{so}(6)\oplus(\mathfrak{so}(4)+\mathfrak{so}(2))
 \oplus \mathfrak{su}(2)$\\

& $K_2$ & $\mathfrak{su}(8)$  &
	     $\mathfrak{so}(6)\oplus\mathfrak{so}(6)\oplus\mathfrak{so}(2)$ \\  

& $K_4$ & $\mathfrak{so}(12)\oplus\mathfrak{su}(2)$ &
	     $\mathfrak{so}(6)\oplus\mathfrak{so}(6)\oplus\mathfrak{su}(2)$ \\

& $K_1 + K_2$ & $\mathfrak{e}_6\oplus\mathbb{R}$  &
	 $\mathfrak{so}(6)\oplus(\mathfrak{so}(4)+\mathfrak{so}(2))
          \oplus \mathfrak{so}(2)$ \\

& $K_1 + K_6$ & $\mathfrak{so}(12)\oplus\mathfrak{su}(2)$ &
 $\begin{array}{l}\!\!\!(\mathfrak{so}(4)+\mathfrak{so}(2))\\ 
\hspace{14pt}\oplus(\mathfrak{so}(4)+\mathfrak{so}(2))\oplus\mathfrak{su}(2)
  \end{array}$\\ 

& $K_3 + K_7$ & $\mathfrak{su}(8)$ &
 $\mathfrak{u}(3)\oplus\mathfrak{u}(3)\oplus\mathfrak{su}(2)$\\ 

& $K_1 + K_2 + K_6$ & $\mathfrak{su}(8)$ &
 $\begin{array}{l}
\!\!\!(\mathfrak{so}(4)+\mathfrak{so}(2))\\
\hspace{14pt}\oplus(\mathfrak{so}(4)+\mathfrak{so}(2))\oplus\mathfrak{so}(2)
\end{array}$ \\ 

& $K_2 + K_3 + K_7$ & $\mathfrak{so}(12)\oplus\mathfrak{su}(2)$ &
 $\mathfrak{u}(3)\oplus\mathfrak{u}(3)\oplus\mathfrak{so}(2)$\\ 

& $K_3 + K_4 + K_7$ & $\mathfrak{e}_6\oplus\mathbb{R}$ &
 $\mathfrak{u}(3)\oplus\mathfrak{u}(3)\oplus\mathfrak{su}(2)$\\ 

\hline

$(\mathfrak{f}_4, \mathfrak{so}(6)\oplus\mathfrak{so}(3), K_3)$ &
 $K_1$ & $\mathfrak{sp}(3)\oplus\mathfrak{su}(2)$  &
	 $(\mathfrak{so}(4)+\mathfrak{so}(2))\oplus\mathfrak{so}(3)$\\    

& $K_3$ & $\mathfrak{so}(9)$  &
	     $\mathfrak{so}(6)\oplus\mathfrak{so}(3)$ \\  

& $K_4$ & $\mathfrak{so}(9)$ &
	     $\mathfrak{so}(6)\oplus\mathfrak{so}(2)$ \\

& $K_1 + K_4$ & $\mathfrak{sp}(3)\oplus\mathfrak{su}(2)$  &
	 $\mathfrak{so}(6)\oplus\mathfrak{so}(2)$\\
\hline
\hline
\multicolumn{1}{c|}{$(\mathfrak{g}, \mathfrak{h}, H)$} &
 \multicolumn{1}{|c|}{$\tau$}
 &\multicolumn{1}{|c|}{$\mathfrak k$} & 
 \multicolumn{1}{|c}{$\mathfrak {h\cap k}$}  
 \\
\hline
\hline
$(\mathfrak{e}_7,
 \mathfrak{so}(6)\oplus\mathfrak{so}(6)\oplus\mathfrak{su}(2), K_4)$ 
& $\varphi$ & $\mathfrak{e}_6\oplus\mathbb{R}$ &
	 $\mathfrak{so}(16)\oplus\mathfrak{sp}(1)$ \\ 

& $\varphi\circ\tau_{K_2}$ & $\mathfrak{su}(8)$ &
	 $\mathfrak{so}(16)\oplus\mathfrak{so}(2)$  \\ 

& $\varphi\circ\tau_{K_4}$ & $\mathfrak{su}(8)$ &
	 $\mathfrak{so}(16)\oplus\mathfrak{sp}(1)$  \\
\hline
\multicolumn{4}{l}{$\varphi\ : \ E_{\alpha_1}\mapsto E_{\alpha_6}, \
 E_{\alpha_2}\mapsto E_{\alpha_2}, \ E_{\alpha_3}\mapsto
 E_{\alpha_5}, \ E_{\alpha_4}\mapsto E_{\alpha_4}, \ E_{\alpha_7}\mapsto
 E_{\alpha_0}$} 

\end{longtable}
\end{center}
}

{\small
\begin{center}
\renewcommand{\arraystretch}{1.2}
\begin{longtable}{l|l|l|l}
\caption{$\dim\mathfrak{z}=1, \ \tau\circ\sigma = \sigma\circ\tau$, 
$\sigma=\tau_{(1/2)H}$ and $\mathfrak{k=g}^{\tau}$.}\\
\hline
\multicolumn{1}{c|}{$(\mathfrak{g}, \mathfrak{h}, H)$} &
 \multicolumn{1}{|c|}{$h \ (\tau=\tau_h)$}
 &\multicolumn{1}{|c|}{$\mathfrak k$} & 
 \multicolumn{1}{|c}{$\mathfrak {h\cap k}$}  
 \\ 
\hline
\hline

$(\mathfrak{e}_8,
 \mathfrak{su}(8)\oplus\mathbb{R}, K_2)$ & $K_1$ & $\mathfrak{so}(16)$  & $
 \mathfrak{s(u}(7) + \mathfrak u(1))\oplus\mathbb{R}$ \\

 & $K_2$ & $\mathfrak{so}(16)$ & $\mathfrak{su}(8)\oplus\mathbb{R} $ \\

 & $K_3$ & $\mathfrak{e} _7\oplus\mathfrak{su}(2)$  &
	     $\mathfrak{s(u}(6) + \mathfrak u(2))\oplus\mathbb{R}$ \\   

 & $K_4$ &  $\mathfrak{e}_7\oplus\mathfrak{su}(2)$  &
	     $\mathfrak{s(u}(5) + \mathfrak u(3))\oplus\mathbb{R}$ \\ 

 & $K_5$ &  $\mathfrak{so}(16)$  &
	     $\mathfrak{s(u}(4) + \mathfrak u(4))\oplus\mathbb{R}$ \\ 

 & $K_8$ & $\mathfrak{e}_7\oplus\mathfrak{su}(2)$  &
	     $\mathfrak{s(u}(7) + \mathfrak u(1))\oplus\mathbb{R}$ \\

 & $K_2 + K_3$ & $\mathfrak{e}_7\oplus\mathfrak{su}(2)$  &
	     $\mathfrak{s(u}(6) + \mathfrak u(2))\oplus\mathbb{R}$ \\  

& $K_2 + K_4$ & $\mathfrak{so}(16)$  &
	     $\mathfrak{s(u}(5) + \mathfrak u(3))\oplus\mathbb{R}$ \\ 

\hline
$(\mathfrak{e}_8,
 \mathfrak{e}_6\oplus\mathfrak{su}(2)\oplus\mathbb{R}, K_7)$ & $K_1$ &
	 $\mathfrak{so}(16)$  &
	     $(\mathfrak{so}(10) + \mathbb{R}) \oplus\mathfrak{su}(2)
 \oplus\mathbb{R}$ \\ 

 & $K_2$ & $\mathfrak{so}(16)$  &
	     $(\mathfrak{su}(6) +
	     \mathfrak{su}(2))\oplus\mathfrak{su}(2)\oplus\mathbb{R}$ \\  

 & $K_7$ & $\mathfrak{e}_7\oplus\mathfrak{su}(2)$ &
	     $\mathfrak{e} _6\oplus\mathfrak{su}(2)\oplus\mathbb{R} $ \\ 

& $K_8$ & $\mathfrak{e}_7\oplus\mathfrak{su}(2)$  &
	     $\mathfrak{e}_6\oplus\mathfrak{so}(2)\oplus\mathbb{R}$ \\ 

 & $K_1 + K_7$ & $\mathfrak{e}_7\oplus\mathfrak{su}(2)$  &
	     $(\mathfrak{so}(10) + \mathbb
	     R)\oplus\mathfrak{su}(2)\oplus\mathbb{R}$ \\

 & $K_1 + K_8$ &  $\mathfrak{e}_7\oplus\mathfrak{su}(2)$  &
	     $(\mathfrak{so}(10) + \mathbb
	     R)\oplus\mathfrak{so}(2)\oplus\mathbb{R}$ \\ 

 & $K_2 + K_7$ & $\mathfrak{e}_7\oplus\mathfrak{su}(2)$  &
	     $(\mathfrak{su}(6) +
	     \mathfrak{su}(2))\oplus\mathfrak{su}(2)\oplus\mathbb{R}$ \\  

 & $K_2 + K_8$ &  $\mathfrak{so}(16)$  &
	     $(\mathfrak{su}(6) +
	     \mathfrak{su}(2))\oplus\mathfrak{so}(2)\oplus\mathbb{R}$ \\ 

\hline

$(\mathfrak{e}_7,
\mathfrak{su}(6)\oplus\mathfrak{su}(2)\oplus\mathbb{R}, K_3)$ & $K_1$ &
	     $\mathfrak{so}(12)\oplus\mathfrak{su}(2)$ & $\mathfrak{su}(6)
             \oplus \mathfrak{so}(2) \oplus \mathbb{R} $ \\ 

 & $K_2$ & $\mathfrak{su}(8)$ & $\mathfrak{s(u}(5) + \mathfrak u(1))
             \oplus \mathfrak{su}(2)\oplus\mathbb{R}$ \\

 & $K_3$ & $\mathfrak{so}(12) \oplus\mathfrak{su}(2)$ &
	     $\mathfrak{su}(6)\oplus\mathfrak{su}(2)\oplus\mathbb{R} $ \\

 & $K_4$ & $\mathfrak{so}(12)\oplus\mathfrak{su}(2)$  &
	     $\mathfrak{s(u} (4) + \mathfrak u(2)) \oplus
             \mathfrak{su}(2) \oplus\mathbb{R}$ \\ 

& $K_5$ & $\mathfrak{su}(8)$  &
	     $\mathfrak{s(u}(3) + \mathfrak u(3)) \oplus
             \mathfrak{su}(2) \oplus\mathbb{R}$ \\

& $K_7$ & $\mathfrak{e} _6 \oplus\mathbb{R}$  &
	     $\mathfrak{s(u}(5) + \mathfrak u(1)) \oplus
             \mathfrak{su}(2) \oplus\mathbb{R}$ \\

& $K_1 + K_2$ & $\mathfrak{e}_6\oplus\mathbb{R}$  &
	     $\mathfrak {s(u}(5) + \mathfrak u(1)) \oplus
             \mathfrak{so}(2) \oplus\mathbb{R}$ \\ 

 & $K_1 + K_4$ &  $\mathfrak{so}(12)\oplus\mathfrak{su}(2)$  &
	     $\mathfrak {s(u} (4) + \mathfrak u(2))\oplus 
             \mathfrak{so}(2)\oplus\mathbb{R}$ \\ 

 & $K_1 + K_5$ &  $\mathfrak{su}(8)$  &
	     $\mathfrak{s(u}(3) + \mathfrak u(3))\oplus
             \mathfrak{so}(2)\oplus\mathbb{R}$ \\ 

 & $K_3 + K_4$ &  $\mathfrak{so}(12)\oplus\mathfrak{su}(2)$  &
	     $\mathfrak {s(u} (4) + \mathfrak u(2))\oplus
             \mathfrak{su}(2)\oplus\mathbb{R}$ \\ 

 & $K_3 + K_5$ &  $\mathfrak{e}_6\oplus\mathbb{R}$  &
	     $\mathfrak {s(u}(3) + \mathfrak
	     u(3))\oplus\mathfrak{su}(2)\oplus\mathbb{R}$ \\ 

\hline

$(\mathfrak{e}_7,
 \mathfrak{su}(5)\oplus\mathfrak{su}(3)\oplus\mathbb{R}, K_5)$ & $K_1$ &
	 $\mathfrak{so}(12)\oplus\mathfrak{su}(2)$ & $\mathfrak{s(u}(4)
	     + \mathfrak u(1))\oplus\mathfrak{su}(3)\oplus\mathbb{R}$ \\  

 & $K_3$ &  $\mathfrak{so}(12)\oplus\mathfrak{su}(2)$ &
	     $\mathfrak{s(u}(3) + \mathfrak u(2))\oplus
             \mathfrak{su}(3)\oplus\mathbb{R}$ \\ 

 & $K_5$ & $\mathfrak{su}(8)$ &
	     $\mathfrak{su}(5)\oplus\mathfrak{su}(3)\oplus\mathbb{R} $ \\

 & $K_6$ & $\mathfrak{so}(12)\oplus\mathfrak{su}(2)$ & $\mathfrak{su}(5) 
           \oplus\mathfrak{s(u}(2) + \mathfrak u(1))\oplus\mathbb{R} $ \\ 

 & $K_7$ & $\mathfrak{e}_6\oplus\mathbb{R}$ & $\mathfrak{su}(5)\oplus 
           \mathfrak{s(u}(2) + \mathfrak u(1))\oplus\mathbb{R}$ \\ 

 & $K_1 + K_5$ & $\mathfrak{su}(8)$ & $\mathfrak{s(u}(4) + \mathfrak u(1))
           \oplus\mathfrak {su}(3)\oplus\mathbb{R}$ \\ 

 & $K_1 + K_6$ & $\mathfrak{so}(12)\oplus\mathfrak{su}(2)$ &
	     $\begin{array}{l}\!\!\!\mathfrak{s(u}(4) + \mathfrak u(1))\\
             \hspace{14pt}\oplus\mathfrak{s(u}(2) + 
             \mathfrak u(1)) \oplus \mathbb{R}\end{array}$ \\ 

 & $K_1 + K_7$ & $\mathfrak{e}_6\oplus\mathbb{R}$ &
	     $\begin{array}{l}\!\!\!\mathfrak{s(u}(4) + \mathfrak u(1))\\
             \hspace{14pt}\oplus\mathfrak{s(u}(2) + \mathfrak u(1))\oplus
             \mathbb{R}\end{array}$\\

 & $K_3 + K_5$ & $\mathfrak{e}_6 \oplus\mathbb{R}$  &
	     $\mathfrak{s(u}(3) + \mathfrak u(2))\oplus\mathfrak{su}(3)
             \oplus\mathbb{R}$ \\

 & $K_3 + K_6$ &  $\mathfrak{so}(12)\oplus\mathfrak{su}(2)$  &
	     $\begin{array}{l}\!\!\!\mathfrak{s(u}(3) + \mathfrak u(2))\\
             \hspace{14pt}\oplus\mathfrak{su}(2) + \mathfrak u(1))\oplus
             \mathbb{R}\end{array}$ \\  

 & $K_3 + K_7$ &  $\mathfrak{su}(8)$ &
	     $\begin{array}{l}\!\!\!\mathfrak{s(u}(3) + \mathfrak u(2))\\
             \hspace{14pt}\oplus\mathfrak{s(u}(2) + \mathfrak u(1))
             \oplus\mathbb{R}\end{array}$ \\

\hline

$(\mathfrak{e}_6,\mathfrak{su}(3)\oplus\mathfrak{su}(3)\oplus\mathfrak{su}(2)
 \oplus\mathbb{R}, K_4)$ & $K_1$ & $\mathfrak{so}(10)\oplus\mathbb{R}$ & $
\begin{array}{l}\!\!\!\mathfrak{s(u}(2) + \mathfrak u(1))\\
\hspace{14pt}\oplus \mathfrak{su}(3) \oplus \mathfrak{su}(2) \oplus
 \mathbb{R}\end{array}$ \\

 & $K_4$ &  $ \mathfrak{su}(6) \oplus\mathfrak{su}(2)$ &  
 $\mathfrak{su}(3)\oplus\mathfrak{su}(3)\oplus\mathfrak{su}(2)
 \oplus\mathbb{R}$ \\ 

 & $K_5$ & $\mathfrak{su}(6)\oplus\mathfrak{su}(2)$ &
	     $\begin{array}{l}\!\!\!\mathfrak{s(u}(2) + \mathfrak u(1))\\
             \hspace{14pt}\oplus \mathfrak{su}(3) \oplus
	      \mathfrak{su}(2) \oplus \mathbb{R}\end{array}$ \\

& $K_1 + K_2$ & $\mathfrak{so}(10)\oplus\mathbb{R}$ & 
              $\begin{array}{l}\!\!\!\mathfrak{s(u}(2) + \mathfrak u(1))\\
              \hspace{10pt}\oplus\mathfrak{su}(3)\oplus\mathfrak{so}(2)
	       \oplus \mathbb{R}\end{array}$ \\ 

& $K_1 + K_5$ & $\mathfrak{su}(6)\oplus\mathfrak {su} (2)$  &
	     $\begin{array}{l}\!\!\!\mathfrak{s(u}(2) + \mathfrak u(1))\\
              \oplus \mathfrak{s(u}(2) + \mathfrak
	      u(1))\oplus\mathfrak{su}(2) \oplus \mathbb{R}\end{array}$ \\

 & $K_2 + K_4$ & $\mathfrak{su}(6)\oplus\mathfrak{su}(2)$ &
 $\mathfrak{su}(3) \oplus \mathfrak{su}(3) \oplus \mathfrak{so}(2)
 \oplus \mathbb{R}$\\

 & $K_1 + K_2+ K_5$ & $\mathfrak{su}(6)\oplus\mathfrak{su}(2)$ &
	     $\begin{array}{l}\!\!\!\mathfrak{s(u}(2) + \mathfrak
	      u(1))\\
             \oplus\mathfrak{s(u}(2) + \mathfrak u(1))
	      \oplus \mathfrak{so}(2)\oplus\mathbb{R}\end{array}$ \\  

 & $K_1 + K_4 + K_5$ & $\mathfrak{so}(10)\oplus\mathbb{R}$  &
	     $\begin{array}{l}\!\!\!\mathfrak{s(u}(2) + \mathfrak u(1))\\
             \oplus\mathfrak{s(u}(2) + \mathfrak u(1))
	      \oplus \mathfrak{su}(2)\oplus\mathbb{R}\end{array}$ \\ 

\hline

$(\mathfrak{f}_4,
 \mathfrak{su}(3)\oplus\mathfrak{su}(2)\oplus\mathbb{R}, K_2)$ & $K_1$ & $
	 \mathfrak{sp}(3)\oplus\mathfrak{su}(2)$ &
	     $\mathfrak{su}(3)\oplus\mathfrak{so}(2)\oplus\mathbb{R} $ \\ 
	    
 & $K_2$ & $\mathfrak{sp}(3)\oplus\mathfrak{su}(2)$  &
	     $\mathfrak{su}(3)\oplus\mathfrak{su}(2)\oplus\mathbb{R}$ \\

 & $K_4$ & $\mathfrak{so}(9)$  & $\mathfrak{s(u}(2) + \mathfrak
	     u(1))\oplus\mathfrak{su}(2)\oplus\mathbb{R} $ \\ 

 & $K_1 + K_3$ & $\mathfrak{sp}(3)\oplus\mathfrak{su}(2)$ &
	     $\mathfrak{s(u}(2) +
	     \mathfrak u(1))\oplus\mathfrak{so}(2)\oplus\mathbb{R}$ \\  

 & $K_2 + K_4$ & $\mathfrak{sp}(3)\oplus\mathfrak{su}(2)$ &
	     $\mathfrak{s(u}(2) + \mathfrak
	     u(1))\oplus\mathfrak{su}(2)\oplus\mathbb{R}$ \\ 

\hline 
 $(\mathfrak{g}_2,
 \mathfrak{su}(2)\oplus\mathbb{R}, K_1)$ & $K_1$ &
	 $\mathfrak{su}(2) \oplus \mathfrak{su}(2)$ & $\mathfrak{su}(2)
          \oplus \mathbb{R}$ \\

& $K_2$ & $\mathfrak{su}(2)\oplus\mathfrak{su}(2)$  &
	     $\mathfrak{so}(2)\oplus\mathbb{R}$ \\
\hline
\hline
\multicolumn{1}{c|}{$(\mathfrak{g}, \mathfrak{h}, H)$} &
 \multicolumn{1}{|c|}{$\tau$}
 &\multicolumn{1}{|c|}{$\mathfrak k$} & 
 \multicolumn{1}{|c}{$\mathfrak {h\cap k}$}\\
\hline
\hline
$(\mathfrak{e}_6,\mathfrak{su}(3) \oplus\mathfrak{su}(3)
 \oplus\mathfrak{su}(2)\oplus\mathbb{R}, K_4 )$ & $\psi$ & $\mathfrak{f} _4$  &
 $\mathfrak{su}(3) \oplus\mathfrak{su}(2) \oplus\mathfrak{sp}(1)
 \oplus\mathbb{R}$ \\ 

&  $\psi\circ\tau_{K_2}$ &  $ \mathfrak{sp}(4)$ &
 $\mathfrak{su}(3) \oplus\mathfrak{su}(2) \oplus\mathfrak{so}(2)
 \oplus\mathbb{R}$ \\ 

&  $\psi\circ\tau_{K_4}$ &  $ \mathfrak{sp}(4)$ &
 $\mathfrak{su}(3)\oplus\mathfrak{su}(2)\oplus\mathfrak{sp}(1)
 \oplus\mathbb{R}$ \\  
\hline
\multicolumn{4}{l}{$\psi\ : \ E_{\alpha_1}\mapsto E_{\alpha_6}, \
 E_{\alpha_2}\mapsto E_{\alpha_2}, \ E_{\alpha_3}\mapsto E_{\alpha_5}, \
 E_{\alpha_4}\mapsto E_{\alpha_4}$}
\end{longtable}
\end{center}
}

\bigskip\par
\begin{flushright}
\begin{tabular}{l}
Department of Computer and Media Science \\
Saitama Junior College \\
Hanasaki-ebashi, Kazo, Saitama 347-8503 \\
Japan \\
e-mail : kurihara@sjc.ac.jp \vspace{1cm}\\

Department of Mathematics \\
Chiba Institute of Technology \\
Shibazono, Narashino, Chiba 275-0023 \\
Japan \\
e-mail : tojo.koji@it-chiba.ac.jp
\end{tabular}
\end{flushright}

\end{document}